\documentclass{article}
\usepackage{graphicx} % Required for inserting images
\usepackage[letterpaper,top=2cm,bottom=2cm,left=3cm,right=3cm,marginparwidth=1.75cm]{geometry}
 \usepackage{algorithm} 
\usepackage{algorithmic}
\usepackage[algo2e,ruled,vlined]{algorithm2e}
\usepackage[noadjust]{cite}
\usepackage{amsmath, amsfonts, amsthm, amssymb, bm}
\usepackage{graphicx,soul}
\usepackage[colorlinks=true, allcolors=blue]{hyperref}
\usepackage{booktabs}
\usepackage{amsfonts}
\usepackage{color,xcolor}
\usepackage{multirow}
\usepackage{soul}

\newcommand{\NSR}{\operatorname{NSR}}

\usepackage{scalerel}[2014/03/10]
\usepackage[usestackEOL]{stackengine}
\def\intavg{\,\ThisStyle{\ensurestackMath{%
    \stackinset{c}{0\LMpt}{c}{0\LMpt}{\SavedStyle-}{\SavedStyle\phantom{\int}}}%
    \setbox0=\hbox{$\SavedStyle\int\,$}\kern-\wd0}\int}
\DeclareMathOperator*{\argmin}{\arg\min}

 \newcommand{\bbE}{\mathbb E}

\newcommand{\bbR}{\mathbb R}

 \newcommand{\bb}{\mathbf b}
\newcommand{\bc}{\mathbf c}

\newcommand{\bu}{\mathbf u}  
  
\newcommand{\by}{\mathbf y}

 \newcommand{\bF}{\mathbf F}

 \newcommand{\bV}{\mathbf V}

\newcommand{\cA}{\mathcal A} 
 \newcommand{\cD}{\mathcal D} 
 \newcommand{\cF}{\mathcal F}
 
\newcommand{\cI}{\mathcal I} 
 \newcommand{\cL}{\mathcal L}
\newcommand{\cM}{\mathcal M} 
\newcommand{\cO}{\mathcal O} \newcommand{\cP}{\mathcal P} 
 
\newcommand{\cS}{\mathcal S} \newcommand{\cT}{\mathcal T}
\newcommand{\cU}{\mathcal U} 
\newcommand{\cW}{\mathcal W} \newcommand{\cX}{\mathcal X} 
\newcommand{\cY}{\mathcal Y} 

\usepackage{booktabs} 
\usepackage{siunitx}  
\usepackage{multirow}
\usepackage{caption}
\usepackage{subcaption}
\usepackage{adjustbox}
\usepackage{authblk}
\newcommand{\Nd}{N_{\operatorname{dict}}}

\title{PDE Identification Using Noise Adaptive  Differentiation in Strong Form (S-IDENT)} 

\author[1]{Roy Y. He}
\author[2]{Sung Ha Kang}
\affil[1]{Department of Mathematics, City University of Hong Kong, Hong Kong}
\affil[2]{School of Mathematics, Georgia Institute of Technology, GA, USA}
\date{}

\begin{document}
\maketitle

\begin{abstract}
We explore identifying partial differential equations (PDEs) from noisy observations of single time-space trajectories. Recent developments show the benefits of identifying PDEs in their weak forms. We investigate the use of differential Strong-form dictionaries for PDE IDENTification (S-IDENT), which enables finding more general linear and nonlinear PDEs. Building on an extensive exploration of integral-type denoised differentiation approaches, we propose to use Savitzky--Golay (SG) differentiation with an adaptive window length chosen based on Stein's Unbiased Risk Estimate (SURE). This offers a guaranteed order of accuracy while producing estimators with minimal variance. The identification process is further refined and stabilized through trimming and reduction-in-residual model selection. Numerical evidence shows that S-IDENT can successfully identify nonlinear PDEs at higher levels of noise than existing strong-form methods, while also yielding results comparable to weak-form approaches. We further verify the effectiveness of S-IDENT through comparisons with various strategies to approximate differential features. We provide numerical evidence that general differential-form dictionaries are larger and more ill-conditioned than those used for weak-form identification, yet S-IDENT does not significantly suffer from this combinatorial increase in dictionary size.
\end{abstract}

\section{Introduction}

Given a collection of time-dependent data, that records the evolution of quantities of interest whether scalar or vector fields, it is valuable to automatically discover a mathematical model that characterizes their dynamics, so that analysis, simulation, prediction, and control can be carried out in a systematic manner.
We explore identifying partial differential equations (PDEs) from single noisy trajectories. Assuming that the underlying PDE is of evolutionary type, we focus on the case where the velocity of the observed dynamics is a linear combination of linear or nonlinear terms involving spatial derivatives. 
One can formulate the problem as a simple linear system using the given data. Let $u$ be a function approximating the observed data, and let $\cF_k$ be the feature terms, which are monomials of partial derivatives of $u$ in space, e.g., $u$, $u_x$, $u_{xx}^2$; then one can represent the underlying PDE
\begin{equation}\label{eq:PDE_general}
 u_t = \sum_{k=1}^K c_k \cF_k \quad \text{ as } \quad \bb = \bF \bc. 
 \end{equation}
Here $\bb$ is a numerical approximation of $u_t$, $\bF$ is the feature matrix where each column $\cF_k$ is a numerical estimation of a feature, and $\bc$ is the coefficient vector to be identified. Adding a sparsity constraint on $\bc$ is a natural choice to find a meaningful equation. There are many recent sparse-regression-based approaches such as~\cite{brunton2016discovering,cheng2025Dyn,kang2021ident,rudy2017data,messenger2021weak,he2022asymptotic,he2022robust,cui2025stoch,tang2025wg,tang2023fourier,tang2023weakident,he2023group,he2024much,cohen2024physics}; for example, SINDy~\cite{rudy2017data} and IDENT~\cite{kang2021ident} utilize different sparsity-inducing mechanisms, e.g., $\ell_0$ and $\ell_1$ regularization, to find candidate models. See \cite{he2025ident} for a unifying explanation of the IDENT approach.

\begin{figure}
\centering
\includegraphics[width=0.98\textwidth]{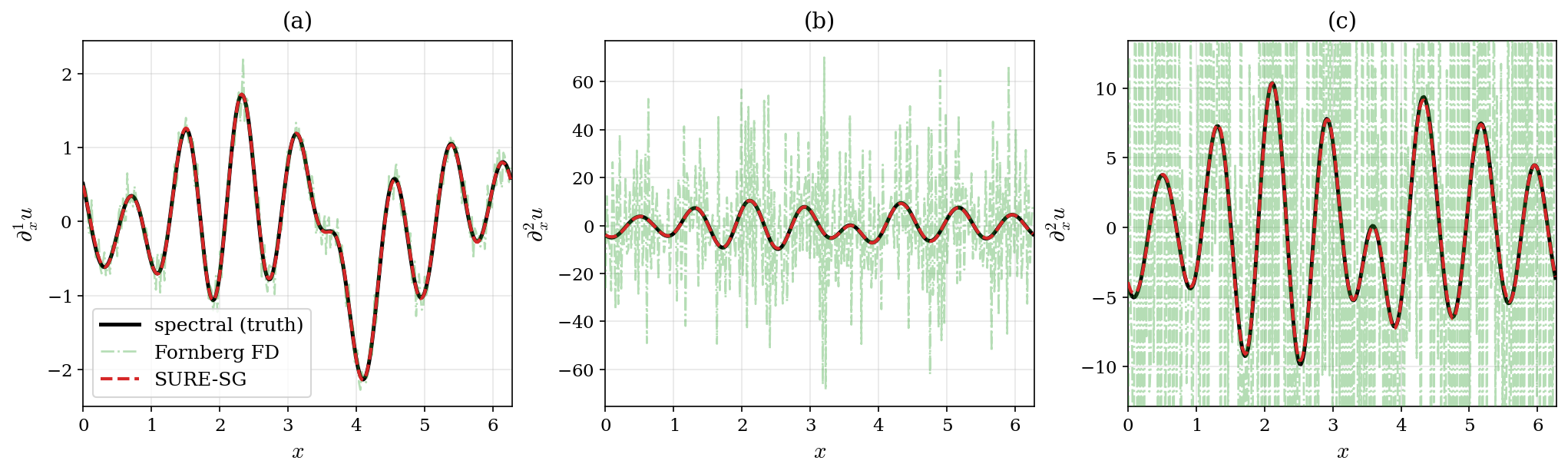}
\caption{Noise amplification due to numerical differentiation. Estimation of the (a) first-order derivative, (b) second-order derivative, and (c) a zoom-in of (b). The small noise in the data (Gaussian with $1\%$ noise-to-signal ratio~\eqref{eq_nsr_energy}) gets amplified every time a finite difference (Fornberg's filter; see Subsection~\ref{sec:fornberg}) is used; see the green dashed lines. In contrast, the adaptive method used in this paper, SURE-SG (the red solid curves; see Subsection~\ref{sec:sure_sg}), yields more stable estimates for the differential features.}
\label{fig:noiseAmp}
\end{figure}

Since both $\bF$ and $\bb$ in the feature system~\eqref{eq:PDE_general} are numerically computed from given data, observational noise causes unstable estimations and adds challenges to the problem. Figure~\ref{fig:noiseAmp} illustrates the noise amplification due to numerical differentiation: a small amount of noise is quickly amplified in (a), (b), and (c) as we approximate first- and second-order derivatives.
In IDENT \cite{kang2021ident}, a moving least-squares (MLS) method is used to denoise the given data. To approximate higher-order-derivative feature terms, Robust-IDENT~\cite{he2022robust} proposes a Successive Denoised Differentiation (SDD) method, where a local weighted MLS denoising is applied whenever one order of differentiation is applied.
Weak-form-based approaches have shown stable state-of-the-art results, such as Weak-SINDy \cite{messenger2021weak} and Weak-IDENT \cite{tang2023weakident}. From the given noisy data, a particular test function is constructed, separating signal from noise in the Fourier domain. Using the weak form, differentiation in feature terms is not directly computed but it is  moved to the smooth test function via integration by parts giving stability in computation.
Such a weak formulation requires specific types of features, in that they should be written as linear differential operators, e.g., $\partial^p_x(u^q)$ or $\partial^p(u^{m}v^n)$ for some positive integers $p, q, m$, and $n$. Tang et al.~\cite{tang2026priorident} showed that prior physical knowledge helps to restrain the feature dictionaries, and that they can be tailored to bias different physical laws, e.g., conservative or dissipative.

For modeling complicated unknown behaviors, more general nonlinear features must be included. For example, when the system has interaction or feedback with its state, such as deformation where the resistance depends on the rate of transport, shear can become thinning, as in blood flow, or thickening, as in wet suspensions. Higher-order terms naturally arise: surface-tension variation depending on the curvature necessarily involves flux proportional to the gradient of the pressure, which can rely on the surface curvature. For example, the thin film equation~\cite{bernis1990higher}
\begin{equation}\label{eq:thin-film-formula}
u_t = -\partial_x(u^2\,u_{xxx}) = -2uu_xu_{xxx} -u^2u_{xxxx}
\end{equation}
contains the fourth-order derivative and products of three terms.
Cross-coupling often occurs in a multi-component system~\cite{drinfeld1981equations}, giving rise to terms like $v\partial_x u$.
These terms cannot always exploit the denoising advantages of the weak formulation. In contrast, frameworks such as IDENT~\cite{kang2021ident} and Robust-IDENT~\cite{he2022robust} can include general terms in the feature dictionary, and we further develop this direction for identifying more complex PDEs with higher levels of noise.

We consider more comprehensive families of dictionaries parameterized by two integers: the highest order of differentiation $p\geq 0$ and the maximal number of multiplications $q\geq 1$.
For any integer $p\geq 0$ and a set of functions $ \cA=\{f_1,\dots,f_N\}\subseteq C^p(\Omega)$ with $\Omega\subseteq \mathbb{R}^D$ for some $D\geq 1$, we define the set of all partial derivatives of functions in $\cA$ up to the $p$-th order as
\begin{equation}\label{eq:DpA}
\mathcal{D}^p(\cA) := \{\partial^\alpha f_n: |\alpha|\leq p,  n=1,\dots,N\},
\end{equation}
where $\alpha=(\alpha_1,\dots,\alpha_D)$ is a multi-index, $|\alpha|:=\sum_{d=1}^D\alpha_d$, and $\partial^\alpha f_n=\partial^{\alpha_1}_1\cdots\partial_D^{\alpha_D} f_n$. The number of elements in $\mathcal{D}^p(\cA)$ is at most $N\cdot{p+D\choose D}$.
For any integer $q\geq 1$, we define the set of all products of up to $q$ functions from $\cA$ as
\begin{equation}\label{eq:MqA}
\mathcal{M}^q(\cA) = \{ f^{\beta_1}_{1}\cdots f^{\beta_N}_{N}:~ \beta_n\geq 0,~\sum_{n=1}^N\beta_n\leq q\}\;.
\end{equation} 
The number of elements in $\mathcal{M}^q(\cA)$ is at most ${q+N\choose N}$.
Using notations in~\eqref{eq:DpA} and~\eqref{eq:MqA}, the strong-form dictionaries considered in differential-form PDE identification such as IDENT~\cite{kang2021ident} and Robust-IDENT \cite{he2022robust} can be expressed as
\begin{equation}\label{eq:typeS}
\text{\textbf{Type-S} $(p,q)$ dictionary :} \quad \cM^q\left(\cD^p(\cA)\right).
\end{equation}
This set consists of features that are defined by taking differentiation first, then multiplication, and it covers any PDEs of the form~\eqref{eq:PDE_general}. Alternatively, the dictionaries exploited in weak-form PDE identification such as \cite{messenger2021weak,tang2023weakident} can be expressed as 
\begin{equation}\label{eq:typeW}
\text{\textbf{Type-W} $(p,q)$ dictionary :} \quad \cD^p\left(\cM^q(\cA)\right).
\end{equation}
 The order of operations is important, since this determines whether the differentiation can be moved to the smooth test function via integration by parts.
By Leibniz's rule, given any $(p,q)$, a Type-W $(p,q)$ dictionary can always be represented as a linear combination of features from the Type-S $(p,q)$ dictionary. Table~\ref{tab_dictionary_sizes} compares the sizes of Type-S and Type-W dictionaries, showing the challenge of identifying PDEs using Type-S dictionaries.
There are additional intrinsic difficulties in identification with more general dictionaries.
We observe that Type-S dictionaries' spectra decay faster than those of Type-W ones, even when their sizes are comparable, indicating significant information redundancy. We explore this in more detail in Appendix \ref{sec:S-correlation}.

\begin{table}
\centering
\renewcommand{\arraystretch}{1.4}
\setlength{\tabcolsep}{6pt}
\begin{tabular}{lcccccccc}
\toprule
\multirow{2}{*}{\textbf{Type}} & \multirow{2}{*}{\textbf{Size formula}} & \multicolumn{3}{c}{$\boldsymbol{D=1,\ N=1}$} & & \multicolumn{3}{c}{$\boldsymbol{D=2,\ N=2}$} \\
\cmidrule(lr){3-5} \cmidrule(lr){7-9}
& $(p,q)$& $(4,3)$ & $(6,4)$ & $(6,6)$ & & $(4,3)$ & $(6,4)$ & $(6,6)$ \\
\midrule
\textbf{Type-S} & $\displaystyle \binom{q + N\binom{p+D}{D}}{q}$ & $56$ & $330$ & $1{,}716$ & & $5{,}456$ & $487{,}635$ & $61{,}474{,}519$ \\[4pt]
\textbf{Type-W} & $\displaystyle 1 + \left[\binom{q+N}{N} - 1\right]\binom{p+D}{D}$ & $16$ & $29$ & $43$ & & $121$ & $365$ & $729$ \\
\bottomrule
\end{tabular}
\caption{Type-S~\eqref{eq:typeS} and Type-W~\eqref{eq:typeW} dictionaries comparison. 
Each dictionary is determined by  the maximal order of derivatives $p$, the maximal number of terms multiplied  $q$, the space dimension $D$ , and the number of variables in the system $N$. The size of Type-S dictionary is significantly bigger than Type-W, indicating the complexity of the problem.  We explore more details in Appendix \ref{sec:S-correlation}.}
\label{tab_dictionary_sizes}
\end{table}

We propose a new method that identifies PDEs consisting of Type-S features (\textbf{S-IDENT}). Since Type-W dictionaries are special cases of Type-S dictionaries, S-IDENT can also find PDEs with Type-W features. Our method integrates an effective noise-adaptive differentiation scheme with powerful IDENT modules developed from a series of works~\cite{kang2021ident,he2022robust,tang2023weakident,he2023group,tang2025wg}. Building on a comprehensive exploration of denoising differentiation, we employ the Savitzky--Golay (SG) filter~\cite{savitzky1964smoothing}, with its window length adaptively determined by minimizing Stein's Unbiased Risk Estimate (SURE)~\cite{stein1981estimation}. This enables automatic derivative estimation without manually tuning the filter to the noise level, yielding more reliable features; see Figure~\ref{fig:adaptive}.
Numerical experiments show that S-IDENT not only handles nonlinear PDEs such as~\eqref{eq:thin-film-formula} which are representable by Type-S features, but also achieves performance comparable to Weak-SINDy~\cite{messenger2021weak} and Weak-IDENT~\cite{tang2023weakident}, which are tailored for Type-W dictionaries and produce integral-form PDEs.
Main contributions of this paper include the following:
\begin{enumerate}
\item We propose a new method, S-IDENT, a PDE identification method with dictionaries containing more general differential and nonlinear Type-S features.
By integrating a noise-adaptive differentiation scheme with powerful modules from IDENT variants, we enhance both identification accuracy and robustness for PDE identification from noisy observations of single trajectories.
\item We conduct a comprehensive study of denoised differentiation techniques based on weighted integration, fostering a clearer understanding of the role of denoising and its relation to the order and nonlinearity of the underlying PDEs. We propose Savitzky--Golay (SG) differentiation with Stein's Unbiased Risk Estimate (SURE)~\cite{stein1981estimation} to achieve automatic window-length selection that adapts to the estimated noise level.
\item We provide various comparisons to validate the proposed model. While uncovering intrinsic difficulties in PDE identification with more general feature dictionaries, we show that S-IDENT remains effective.  S-IDENT can identify more general differential equations, with comparable results to the state-of-the-art Weak form based methods, with much larger size of dictionary. 
\end{enumerate}

This paper is organized as follows. In Section~\ref{sec_core_ideas}, we present our proposed method, S-IDENT, for identifying PDEs from a single noisy observation.
In Section~\ref{sec_review_diff}, we provide a comprehensive overview of derivative estimation methods from discrete data, with a focus on the weighted-integration type~\eqref{eq_general_filter}, which justifies the proposed choice of approximating various  features.
In Section~\ref{sec_numerical}, we present numerical experiments to validate S-IDENT via general performance and comparative studies. We conclude this paper in Section~\ref{sec:conclude}. More experimental details and results including the study of feature correlation are collected in the Appendix.

\section{The Proposed Method: S-IDENT}\label{sec_core_ideas}

\begin{figure}
\centering
\includegraphics[width=0.8\textwidth]{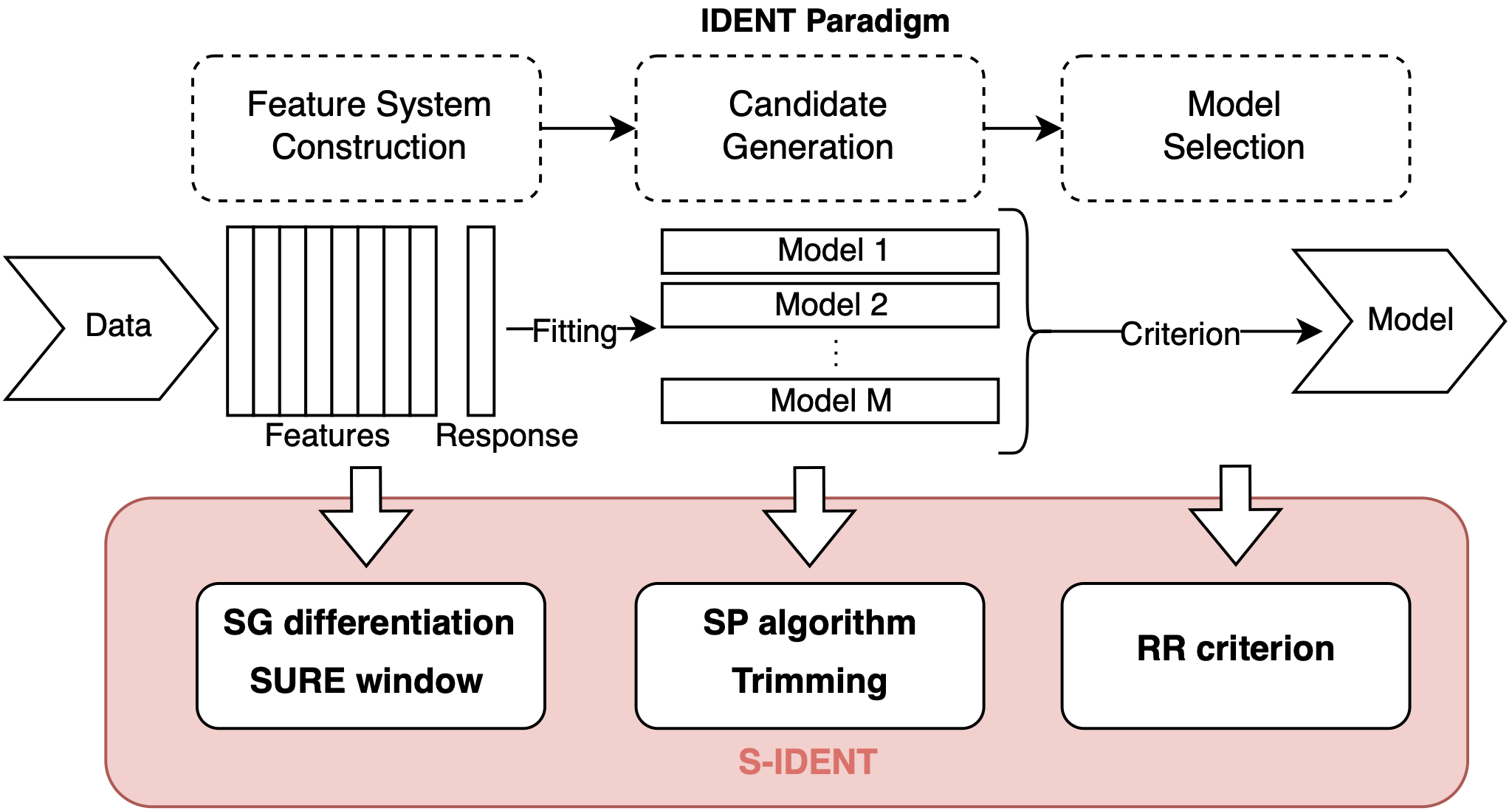}
\caption{The proposed S-IDENT framework (Subsection~\ref{sec:overview}): (a) Savitzky--Golay (SG) differentiation with Stein's Unbiased Risk Estimate (SURE) adaptive window (Subsection~\ref{sec:sure_sg}), (b) Subspace Pursuit (SP) greedy algorithm for candidate generation and feature trimming, and (c) Reduction in Residual (RR) criterion (Subsection~\ref{sec:identification}). Table~\ref{tab:table_summary} relates S-IDENT to other methods.}\label{fig:demo}
\end{figure}

In this section, we present our method, S-IDENT, for PDE identification with a Type-S dictionary from noisy observations of single trajectories. We aim to identify PDEs in their differential forms with enhanced robustness compared to IDENT \cite{kang2021ident} and Robust-IDENT \cite{he2022robust}. Since Type-W is a special case of Type-S, we also aim to achieve performance comparable to that of state-of-the-art methods such as Weak-SINDy \cite{messenger2021weak} and Weak-IDENT \cite{tang2023weakident}.
Figure~\ref{fig:demo} demonstrates the paradigm of IDENT approaches, while highlighting the improved modules in the proposed S-IDENT.

\subsection{General framework of IDENT approach}\label{sec:overview}

Let $\cU = \{(t_n, x_i, U_{i}^n)\;, n=1,\dots,N, i=1,\dots,I\}$ be the given data on a regular time-space grid $\Gamma:=\{(t_n,x_i): n=1,\dots,N,i=1,\dots, I\}\subset [0,T]\times \Omega\subset\mathbb{R}\times\mathbb{R}^D$ with some finite $T>0$ and a bounded domain $\Omega$; here $D\geq 1$ denotes the space dimension. For an integer $K\geq 1$, we consider $K$ candidate feature operators $\cF_k: \cX \to\cY$, $k=1,\dots,K$, mapping between Hilbert spaces $\cX\subset \cap_{k=1}^{K}\text{dom}(\cF_k)$ and $\cY$. Assuming that $\cX\neq \varnothing$, we model the observed data as
\[
U_i^n = \widetilde{u}(t_n,x_i) := u(t_n, x_i) + \varepsilon_{i}^n, \quad  n=1,\dots,N, \quad i=1,\dots,I,
\]
for a function $u\in \cX$ and independent and identically distributed (i.i.d.) perturbations $\varepsilon_i^n$ with mean $\bbE(\varepsilon_i^n)=0$ and variance $\sigma^2:=\text{Var}(\varepsilon_i^n)<+\infty$. We propose to identify an evolutionary-type PDE~\eqref{eq:PDE_general} based on the observed noisy trajectory $\cU$. The coefficients $\{c_k:k=1,\dots,K\}$ are unknown, and most of them are assumed to be zero for interpretability. The framework of IDENT \cite{kang2021ident} and its variants (see \cite{he2025ident} for a review) finds the non-zero coefficients via the following steps:

\begin{enumerate}
\item \textbf{Feature system construction.} From the observed trajectory data $\cU$, approximate the feature values $\cF_k$ at the grid points in $\Gamma$ for $k = 1, \dots, K$. This step yields a feature matrix $\bF$ with $K$ columns, each $\cF_k$ approximating  a feature, together with a feature response vector $\bb$ associated with $u_t$. In the differential-form IDENT~\cite{kang2021ident,he2022robust,he2023group}, each row corresponds to the evaluation of features at a space-time grid point, whereas in the weak-form IDENT~\cite{tang2023weakident,tang2025wg}, each row corresponds to a test function localized at a space-time grid point. We propose a new noise-adaptive strategy for S-IDENT, which is also a differential-form, in Subsection~\ref{sec:sure_sg}.
\item \textbf{Candidate model generation.} Sparse solutions are considered as candidate differential equations. In IDENT~\cite{kang2021ident}, $\ell_1$-regularization is considered and LASSO \cite{tibshirani1996regression} is used, while Robust-IDENT~\cite{he2022robust} and subsequent IDENT variants~\cite{he2023group,he2024much,tang2023weakident,tang2025wg,tang2023fourier,tang2026priorident} adopt $\ell_0$-constrained optimization and greedy algorithms. Candidate differential equations are computed for each level of sparsity. In S-IDENT, we also use a greedy algorithm, Subspace Pursuit (SP)~\cite{dai2009subspace}, with additional trimming introduced in Weak-IDENT \cite{tang2023weakident} for more refined identification.
\item \textbf{Optimal model selection.} From a list of candidate models generated in the previous step for each sparsity level, we use model validation to choose an optimal solution. Model validation methods include the Time Evolution Error (TEE)~\cite{kang2021ident} and the Cross-validation Estimation Error (CEE)~\cite{he2022robust}, and we use
the Residual Reduction (RR)~\cite{he2023group}  based on regression efficiency. Identification improvements for S-IDENT are presented in Subsection~\ref{sec:identification}.
\end{enumerate}

The first step is to set up the feature linear system, and the second and third steps constitute the identification of the PDE. The proposed S-IDENT improves both parts by introducing a noise-adaptive differentiation method and integrating modules for more accurate and robust identification.

\subsection{Savitzky--Golay (SG) filter and Stein’s Unbiased Risk Estimate (SURE) minimization}\label{sec:sure_sg}

We apply Savitzky--Golay (SG) differentiation~\cite{savitzky1964smoothing} to approximate the partial differential features from noisy data. This choice is motivated by its theoretical and computational advantages: it is known that, for white Gaussian noise, the least-squares polynomial fit is the BLUE (best linear unbiased estimator) by Gauss--Markov~\cite[Chapter 4a]{rao1973linear}. This minimizes the estimator's variance subject to exactly reproducing polynomials up to degree $d$ within the window. As detailed in Subsection~\ref{sec:sg}, the local least-squares fitting~\eqref{eq_SG_fitting} is equivalent to linear filtering by exploiting an orthogonal basis, and denoising can be applied efficiently via convolution.
Savitzky--Golay differentiation is specified by two parameters: the degree of the local polynomial and the window length. In general, increasing the polynomial degree adds model complexity and thus risks overfitting, while a larger window length increases bias and a smaller window length increases variance.

We propose to automatically determine the window length parameter by minimizing Stein's Unbiased Risk Estimate (SURE)~\cite{stein1981estimation}, while fixing the degree of the polynomial to be greater than the highest order of derivatives in the given dictionary.
The SURE offers a convenient and systematic method to select hyper-parameters of least-squares estimators. Suppose $x_1,\dots,x_N$ are uniform grid points on the interval $\cI := [x_1, x_N]$ with spacing $\Delta x = x_2-x_1$, $u:\cI\to\mathbb{R}$ is the underlying function extended periodically (so that $u(x_{N+1}) = u(x_1)$), $\bu:=(u(x_1),\dots,u(x_N))$, and $\by=(y_1,\dots,y_N)\in\mathbb{R}^N$ with $y_i = u(x_i)+\varepsilon_i$, where $\varepsilon_i$ is independently sampled from a Gaussian distribution with mean $0$ and standard deviation $\sigma>0$. Let $\widehat{u}_L$ be the SG estimator obtained from $\by$ with window size $L$ greater than the polynomial degree, and define $\widehat{\bu}_L:=(\widehat{u}_L(x_1),\dots,\widehat{u}_L(x_N))$; then the estimator's mean squared error (MSE) is
 \begin{equation}\label{eq_MSE}
\bbE\|\bu - \widehat{\bu}_L\|_2^2 = -N\sigma^2 + \bbE\|\by-\widehat{\bu}_L\|_2^2 + 2\sum_{i=1}^N\operatorname{Cov}(y_i,\widehat{u}_L(x_i))\;.
 \end{equation}
 By Stein's lemma~\cite{stein1981estimation}, the sample-based quantity
  \begin{equation}\label{eq_SURE_estimate}
\widehat{R}_L(\sigma) := -N\sigma^2 + \|\by-\widehat{\bu}_L\|_2^2 + 2\sigma^2\sum_{i=1}^N\frac{\partial \widehat{u}_L(x_i)}{\partial y_i}
 \end{equation}
is an unbiased estimator for~\eqref{eq_MSE}, and $\widehat{R}_L(\sigma)$ in \eqref{eq_SURE_estimate} is called Stein's Unbiased Risk Estimate (SURE). For SG estimation and general filter-based methods, as reviewed in Section~\ref{sec_review_diff}, the divergence $\partial \widehat{u}_L(x_i)/\partial y_i$ reduces to the central weight. We select the window length $L^*$ by minimizing the SURE and define
\begin{equation}\label{eq_SURE_window}
L^* = \operatorname*{arg\,min}_{L\in \Lambda} \widehat{R}_{L}(\sigma),
\end{equation}
where $\Lambda\subset\mathbb{N}$ is a finite set of candidate window lengths.
To apply~\eqref{eq_SURE_window}, we estimate the noise standard deviation $\sigma$ by
\begin{equation*}
\widehat{\sigma}_N=\sqrt{\frac{1}{6N}\sum_{i}\!\left(\delta^2 y_i - \frac{1}{N}\sum_j \delta^2 y_j\right)^{\!2}},
\end{equation*}
where $\delta^2 y_i := y_{i+1} + y_{i-1} - 2y_i$ for $i = 1, \dots, N$, with indices taken modulo $N$. If $|\cI|$ is fixed and $u \in C^6(\cI)$, then $\widehat{\sigma}_N^2$ is an asymptotically unbiased estimator of $\sigma^2$ as $N \to \infty$, with
\[
\mathbb{E}[\widehat{\sigma}_N^2]
= \left(1 - \frac{2}{3N^2}\right)\sigma^2
+ \underbrace{\frac{1}{6N}\sum_i\left(\delta^2u(x_i)-\frac{1}{N}\sum_j\delta^2u(x_j)\right)^2}_{:=V_N(u)}.
\]
The term $V_N$ is related to the local variation of the function $u$, and it can introduce additional variability into the estimator when $u$ is oscillatory.

\begin{figure}
\centering
\includegraphics[width=0.98\textwidth]{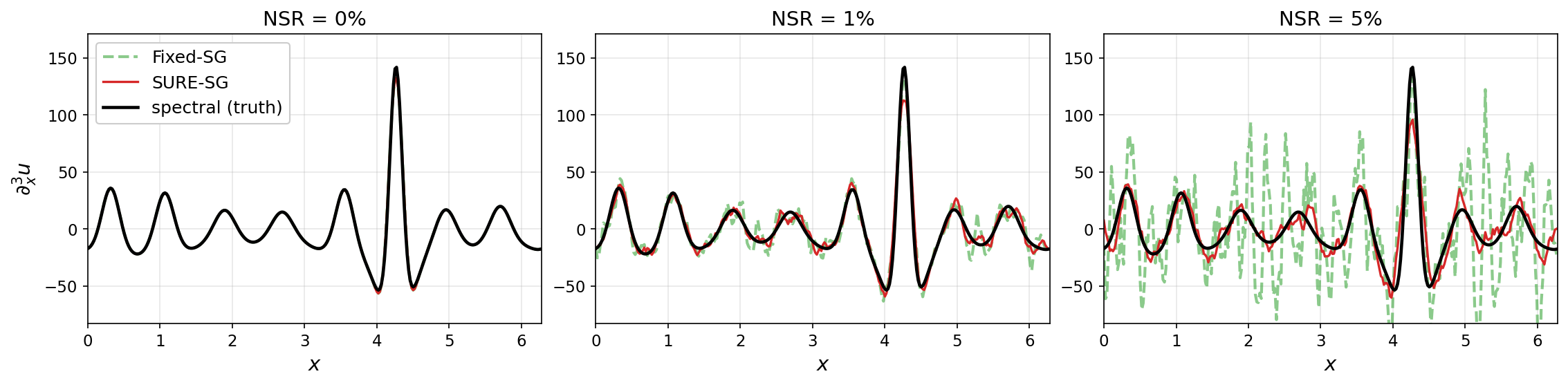}
\caption{SURE-SG: Stein’s Unbiased Risk Estimate (SURE) to automatically choose adaptive window length (Subsection~\ref{sec:sure_sg}) for 
Savitzky--Golay differentiation (Subsection~\ref{sec:sg}). 
Green dashed curves use a fixed window size of $39$; red solid curves use the proposed SURE-SG, giving window sizes $39$, $47$, and $61$ for NSRs of $0\%$, $1\%$, and $5\%$, respectively, when approximating the third-order derivative.}\label{fig:adaptive}
\end{figure}

Figure~\ref{fig:adaptive} shows an example of approximating the third-order
derivative. Green dashed curves use a fixed window length, while red solid curves use the proposed Savitzky--Golay differentiation with an adaptively chosen window length via SURE minimization. The adaptive selection yields window sizes $39$, $47$, and $61$ for NSRs of $0\%$, $1\%$, and $5\%$, respectively. The proposed denoising gives a stable and accurate approximation.

When the data is not periodic, we discard the estimates within $\lfloor L/2 \rfloor$ points of each boundary, where the SG window extends beyond the data. For higher-dimensional data, we apply the tensor product of the one-dimensional SG filters along each dimension, with the window size for each dimension selected by minimizing the SURE aggregated over all one-dimensional slices along that dimension, weighted by interior length.

\subsection{Improved model identification for S-IDENT}\label{sec:identification}

From the given data, Savitzky--Golay (SG) differentiation with SURE, abbreviated as SURE-SG, is used to construct the feature linear system $(\bF,\bb)$.
As we require that the candidate vectors are sparse, i.e., many entries are zero, we obtain candidates by considering the minimization problems
\begin{equation}\label{eq_sp_sparse}
\min_{\bc}\|\bF\bc-\bb\|_2^2\;,~\text{s.t.}~\|\bc\|_0 = m.
\end{equation}
Here $\|\bc\|_0\in\{1,\dots,K\}$ denotes the number of non-zero entries in the vector $\bc$; hence, a solution of~\eqref{eq_sp_sparse} yields a candidate model with exactly $m$ active features. 
We find a solution of~\eqref{eq_sp_sparse} for each sparsity level $m=1,\dots,M$, and define a list of $M$ candidate coefficient vectors as $\mathcal{L}:=\{\bc_1,\bc_2,\dots,\bc_M\}\subset\mathbb{R}^K$, such that $\bF\bc_m$ approximates $\bb$ for $m=1,\dots,M$; here $M\leq K$ denotes the total number of generated candidates. It was found in~\cite{he2022robust} that Subspace Pursuit (SP)~\cite{dai2009subspace} effectively solves the problem~\eqref{eq_sp_sparse}. 

We also adopt \textit{trimming}~\cite{tang2023weakident}, which was introduced to improve the identification accuracy for the weak formulation. Specifically, for each identified $\bc=(c_1,\dots,c_K)\in\cL$, we set the $k$-th entry to zero if
\begin{equation}\label{eq:trimming}
\frac{|c_k|\,\|\mathbf{f}_k\|}{\max_{\ell=1,\dots,K}\{|c_\ell|\,\|\mathbf{f}_\ell\|\}} < \tau,
\end{equation}
where $\mathbf{f}_k$ is the $k$-th column of $\bF$, and $\tau\in(0,1)$ is a threshold parameter. After trimming, we refit the data to estimate the non-zero entries. Denote the resulting modified candidate coefficients by $\{\bc_k':k=1,\dots,M\}$. If $\bc'_k=\bc'_{j}$ for some $k<j$, we discard $\bc'_j$; we denote the resulting set of vectors as $\cL'$, with $M'=|\cL'|$ its size. For more details about SP and trimming, we refer the reader to~\cite{he2022robust,tang2023weakident}.

Among candidates in $\cL'$, we determine the optimal one by the Reduction in Residual (RR) criterion proposed in~\cite{he2023group}. First, we compute the residual sum of squares
$
R_m=\|\bF\bc'_m-\bb\|_2^2,
$
for $m=1,\dots,M'$. 
Let $L\geq 1$ be a fixed integer. If $M' \leq L$, we take the optimal
\begin{equation*}
\bc^*:=\arg\min_{\bc'_m\in\cL'} R_m.
\end{equation*}
Otherwise, for $m=1,\dots, M'-L$, we compute the Reduction in Residual (RR) given by 
 \begin{align}
 s_m = \frac{R_{m}-R_{m+L}}{LR_1}\;,\;\; m=1,\dots, M'-L,  \label{eq:cand_score}
 \end{align}
and choose the optimal candidate as $\bc^* = \bc'_{m^*}$, where
\begin{align}
m^*=\min\{m:1\leq m\leq M'-L,\ s_m<\rho\}, \label{eq_thresh}
\end{align}
for some threshold parameter $\rho>0$. It is the smallest sparsity index $m$ for which $s_m$ is below $\rho$. As reported in~\cite{he2023group}, the identification is not sensitive to the choice of $L$ and $\rho$.

For identification results, we pay attention to two aspects of the coefficient vector $\bc$: the support of $\bc$ indicated by $S$, and the coefficient values of $\bc$.    We describe the evaluation  metrics in Section \ref{sec_numerical}.

\subsection{Relations to existing methods} \label{sec:summaryTable}

We present the relations between existing sparsity-based frameworks for PDE identification and the proposed S-IDENT. We summarize representative methods in Table~\ref{tab:table_summary}, and we consider the following aspects: (Type) the type of features contained in the dictionary; (Approx.) feature approximation methods; (Param.) methods for choosing the parameters of approximation methods; (Spars.) algorithms for generating sparse candidate models; (Trim.) whether feature trimming~\eqref{eq:trimming} is adopted; and (Selec.) the criterion for model selection.

In particular, S-IDENT is the first identification method combining SG with SURE to enable noise-adaptation.
In IDENT~\cite{kang2021ident} and Robust-IDENT~\cite{he2022robust}, the moving least squares (MLS) and MLS combined with Successive Denoised Differentiation (SDD) required manual tuning of the parameters.
The noise-adaptation in Weak-SINDy and Weak-IDENT is automatic but operates in the frequency domain. The trimming technique was introduced in Weak-IDENT~\cite{tang2023weakident} for Type-W dictionaries, and S-IDENT adopts it here for more general Type-S dictionaries. For the sparsity-based candidate generation, S-IDENT uses Subspace Pursuit (SP)~\cite{dai2009subspace} to address $\ell_0$-constrained problems, which differs from STRidge~\cite{rudy2017data} for $\ell_0$-regularization and Lasso~\cite{kang2021ident} for $\ell_1$-regularization, where continuous regularization parameters instead of integer parameters control the candidate sparsity. In GP-IDENT~\cite{he2023group} and WG-IDENT~\cite{tang2025wg}, which deal with varying coefficients, SP was extended to Group Projected Subspace Pursuit (GPSP).
As for the model selection methods, SINDy-PDE and Weak-SINDy both use statistical metrics: the Akaike information criterion (AIC) and the Bayesian information criterion (BIC); IDENT adopts the Time Evolution Error (TEE); and Robust-IDENT uses the Multishooting TEE (MTEE) as well as the Cross-validation Estimation Error (CEE); GP-IDENT and WG-IDENT introduced the Reduction in Residual (RR)~\eqref{eq:cand_score} when focusing on varying coefficients, and S-IDENT exploits this technique on PDEs with constant coefficients.

\begin{table}
\centering
\begin{tabular}{l|cccccc}
\toprule
Method & {Type} & {Approx.} & {Param.} &  {Spars.} & {Trim.} & {Selec.} \\
\midrule
SINDy-PDE \cite{rudy2017data} & S & SG & Manual  & STRidge &No  & AIC/BIC \\
{IDENT} \cite{kang2021ident} & S & MLS & Manual  & Lasso & No & TEE \\
{Robust-IDENT} \cite{he2022robust} & S & SDD+MLS & Manual  & SP &No  & MTEE/CEE \\
{GP-IDENT} \cite{he2023group} & S & SDD+SG & Manual & GPSP & No & RR \\
\textbf{S-IDENT} (Proposed) & \textbf{S} & \textbf{SG} & \textbf{Auto}  & \textbf{SP} & \textbf{Yes} & \textbf{RR} \\
{Weak-IDENT} \cite{tang2023weakident}& W  & Weak form & Auto & SP & Yes & CEE \\
{WG-IDENT} \cite{tang2025wg} & W & Weak form & Auto & GPSP & Yes & RR \\
{Weak-SINDy} \cite{messenger2021weak}& W  & Weak form & Auto & STRidge & No & AIC/BIC \\
\bottomrule
\end{tabular}
\caption{The proposed S-IDENT in relation to existing methods. Abbreviations are explained in the text. S-IDENT is first to adaptively choose approximation parameter for Strong form features, and combines benefits of Trimming and Reduction in Residual for more stable identification. Identification comparisons are presented in Section \ref{sec_numerical}. }\label{tab:table_summary}
\end{table}

Table~\ref{tab:table_summary} shows that S-IDENT advances identification by using the Type-S dictionary, which allows for a much wider range of feature terms. S-IDENT employs Type-S features, SG with automatic parameter selection via SURE, SP with trimming, and RR for model selection, which together represent several advances in the identification of PDEs.

\section{Denoised Differentiation from Discrete Data}\label{sec_review_diff}

We present the motivation behind the choice of the SURE-SG approach for denoised differentiation in feature approximation.
In particular, we analyze and provide a comprehensive review of the class of denoised differentiation methods based on integration. There is a balance one needs to strike between accuracy and stability.
After presenting some general analysis in Subsection~\ref{sec:bias}, we discuss methods deduced from different principles for derivative approximation: (1) Fornberg's finite difference~\cite{fornberg1988generation} based on local polynomial interpolation in Subsection~\ref{sec:fornberg}, (2) the Savitzky--Golay filter~\cite{savitzky1964smoothing} based on local polynomial fitting in Subsection~\ref{sec:sg}, (3) maximally flat differentiation~\cite{hosseini2017finite} based on frequency response specification in Subsection~\ref{sec:maxflat}, and (4) some other methods in Subsection~\ref{sec:other} to highlight the richness of the methods for denoised differentiation. We further discuss the relation between derivative approximation via integration and kernel convolution in Subsection~\ref{sec:WeakvsStrong}.

\subsection{Integration based denoised differentiation}\label{sec:bias}

Let $s\geq 0$ be an integer, let $\mu_s:[-1,1]\to\mathbb{R}$ be a function of bounded variation, and let $g\in L^1(\mathbb{R})$. For some $h\in(0,1)$, we consider the following weighted integral operator for approximating $s$-th order derivatives:
\begin{equation}\label{eq_general_filter}
\cT_{\mu_s}^h(g)(x):=h^{-s}\int_{-1}^1 g(x+th)\,d\mu_s(t).
\end{equation}
This generalizes linear filters in signal processing. For instance, if $\mu_s(t) = \sum_{i=1}^N w_{s}(i)H(t-t_i)$, where $w_{s}(i)\in\mathbb{R}$, $i=1,\dots,N$, $-1\leq t_1<\cdots<t_N\leq 1$, and $H$ is the Heaviside function, then the weighted integral form~\eqref{eq_general_filter} reduces to a weighted sum, approximating derivatives of $g$ as
\begin{equation}\label{eq_linear_filter}
  \cT_{\mu_s}^h(g)(x) = h^{-s}\sum_{i=1}^N w_{s}(i)g(x+t_ih),
 \end{equation}
which is a local weighted average of the discrete signal. If $t_1,\dots,t_N$ are equidistant with $\Delta t =t_2-t_1$ and $x\in \{nh\Delta t: n\in\mathbb{Z}\}$, then the weighted sum~\eqref{eq_linear_filter} is exactly a time-invariant linear filtering of the signal $g$ sampled with frequency $(h\Delta t)^{-1}$.
These include Fornberg's differentiation~\cite{fornberg1988generation}, Savitzky--Golay filtering~\cite{savitzky1964smoothing}, maximally flat differentiation~\cite{selesnick1998maximally,hosseini2017finite}, and the more general \emph{differentiation by integration}~\cite{lanczos1964evaluation,groetsch1998lanczos}. Such approaches are not only widely used in data analysis but also form a core component of many nonlocal models~\cite{du2019nonlocal}.

\begin{figure}
\centering
\begin{tabular}{c@{\vspace{2pt}}c@{\vspace{2pt}}c}
(a)&(b)&(c)\\
\includegraphics[width=0.33\textwidth,height=2in]{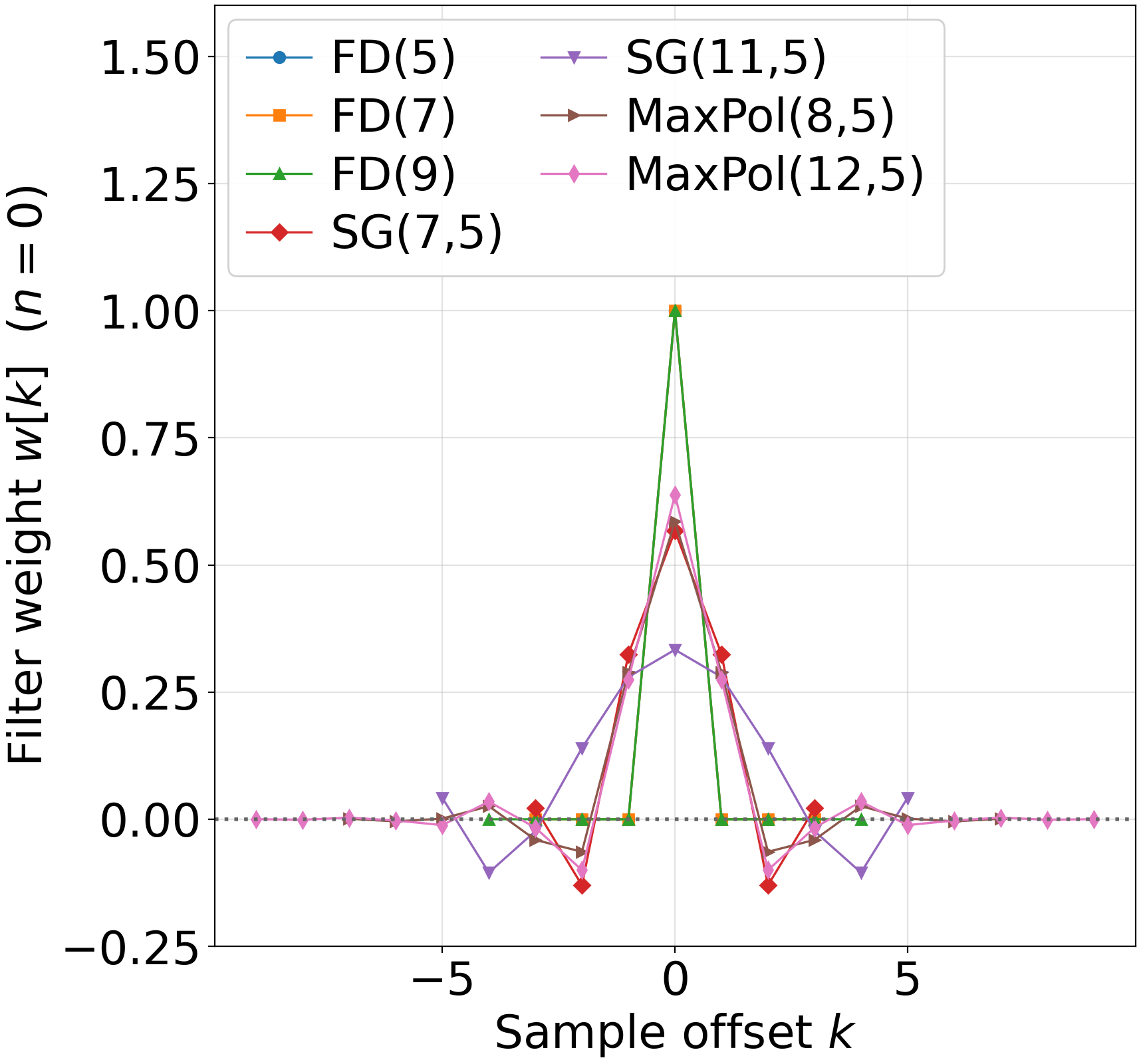}&
\includegraphics[width=0.33\textwidth,height=2in]{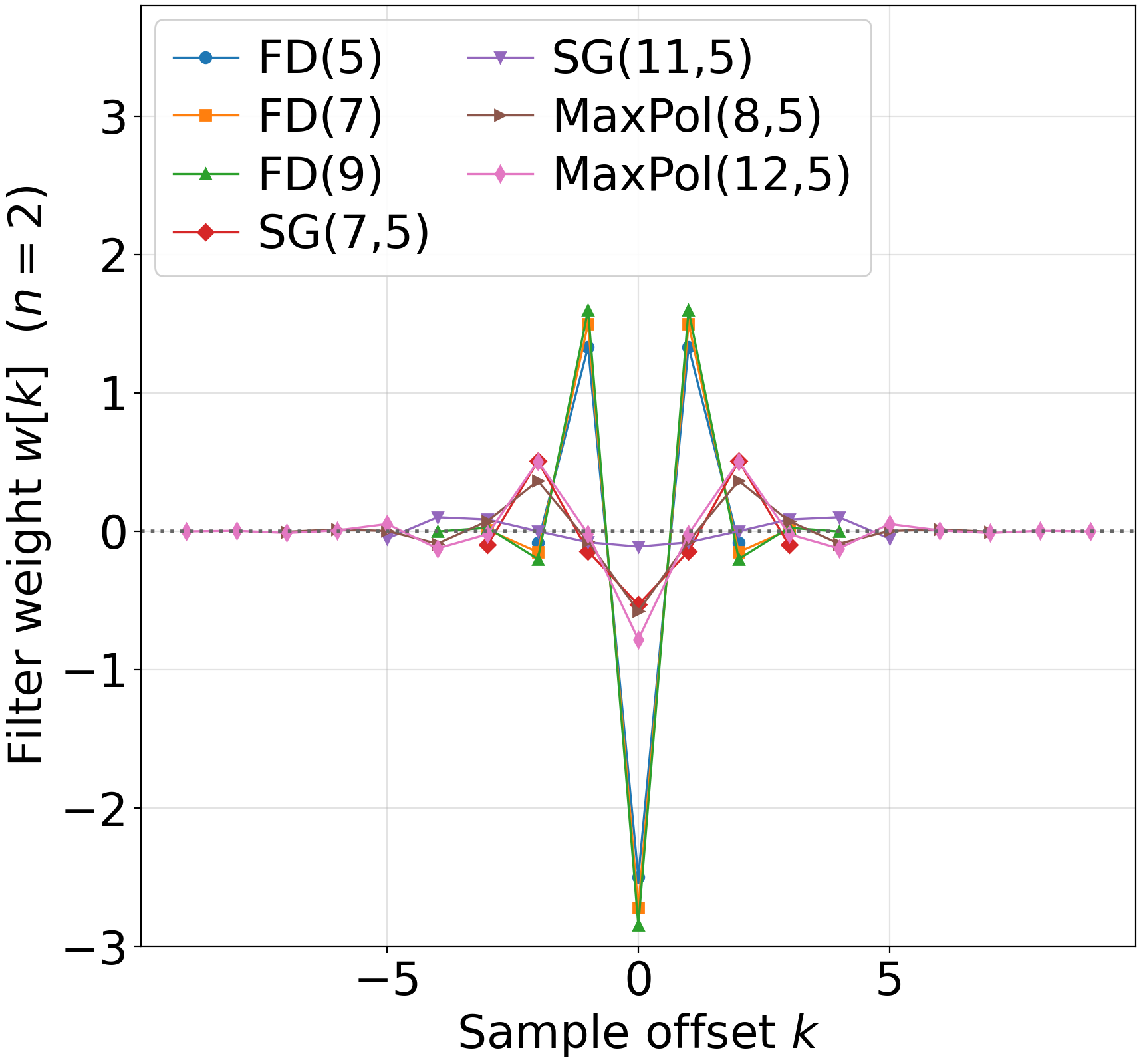}&
\includegraphics[width=0.33\textwidth,height=2in]{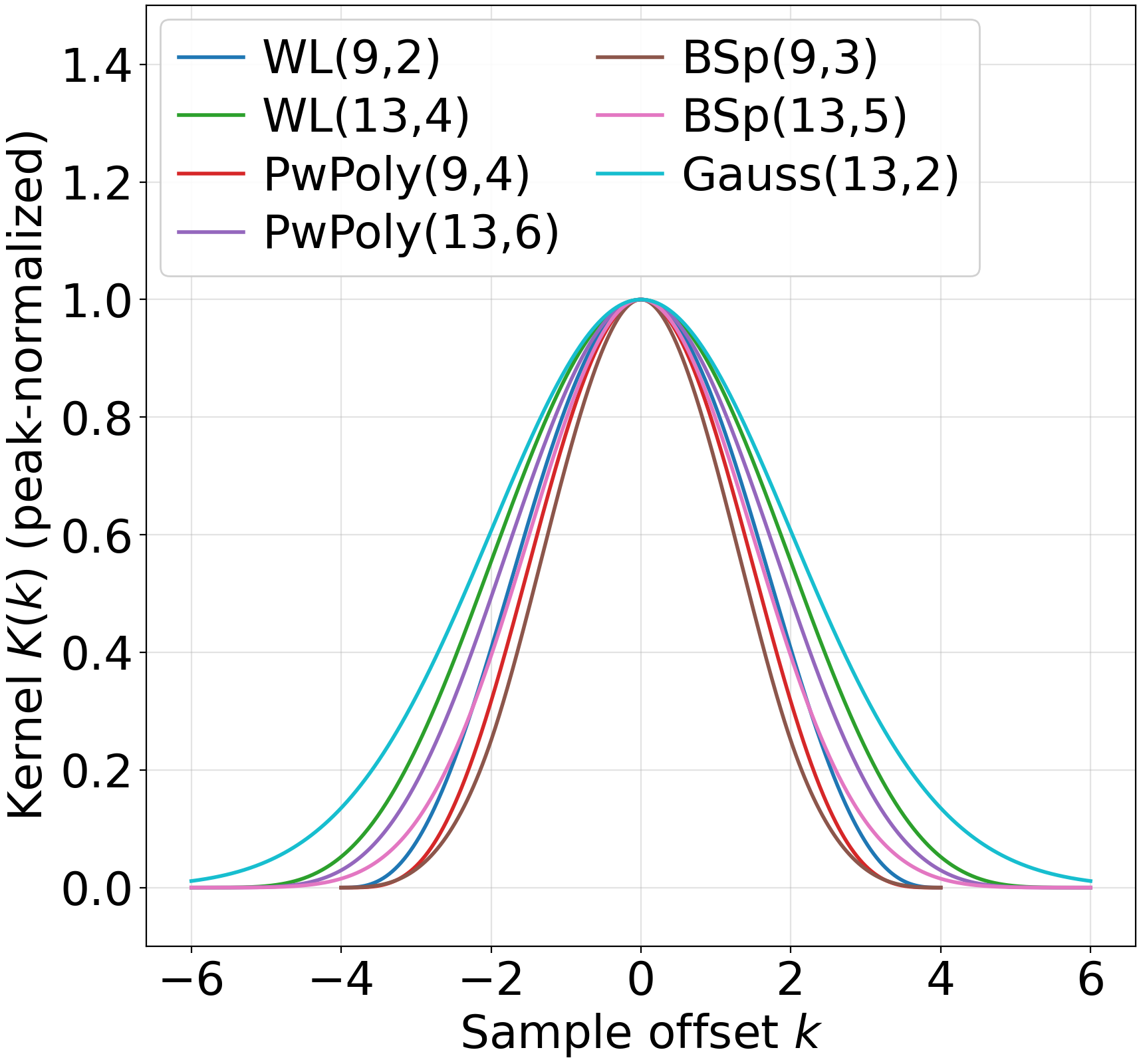}
\end{tabular}
\caption{Weight difference of various methods in terms of differentiation by integration.
In the form of the weighted sum~\eqref{eq_linear_filter}, (a) shows the weights when the methods are applied as a smoothing operator, and (b) shows  when the methods are applied to approximate the second-order derivative. (c) shows the kernel/test function in the weighted integral form~\eqref{eq_general_filter}.
(Abbreviations and details are explained in the text.) 
Figure~\ref{fig:mag-response} further shows their characterization in the frequency domain.}\label{fig:filter-shape}
\end{figure}

In Figure~\ref{fig:filter-shape} (a) and (b), we present the weights of various methods of the form~\eqref{eq_linear_filter} that can be candidates for approximating Type-S features. In particular, (a) shows the weights when the methods are applied as a smoothing operator ($0$-th order derivative), and (b) shows the weights when the methods are applied to approximate the second-order derivative.
In (a) and (b), Fornberg's finite difference~\eqref{eq_weight_Fornberg}, $\mathrm{FD}(N)$, of order $N=5, 7$, and $9$, the Savitzky--Golay method~\eqref{eq_weight_SG}, $\mathrm{SG}(p,q)$, with window size $p$ and polynomial degree $q$ for $(p,q)=(7, 5)$ and $(11,5)$, and the MaxPol method~\eqref{eq_weight_maximal_flat}, $\mathrm{MaxPol}(B,A)$, with smoothing order $B$ and order of accuracy $A$ for $(B,A)=(8,5)$ and $(12,5)$ are shown.
For the smoothing operator in (a), there are oscillations and negative weights for SG and MaxPol, while FD only has positive weights. When applied to differentiation, FD has weights with more significant oscillations, which are responsible for the amplified noise effects, as shown in Figure~\ref{fig:noiseAmp}. For comparison, we present the shape of test functions used for the weak formulation in (c): they include the Wendland kernel~\cite{wendland1995piecewise}, $\mathrm{WL}(p,q)$, with window size $p$ and degree $q$ for $(p,q)=(9,2)$ and $(13,4)$; the truncated polynomial kernel~\cite{messenger2021weak}, $\mathrm{PwPoly}(p,q)$, for $(p,q)=(9,4)$ and $(13,6)$; the B-spline kernel~\cite{tang2025wg}, $\mathrm{BSp}(p,q)$, for $(p,q)=(9,3)$ and $(13,5)$; and the Gaussian kernel with window size $h$ and standard deviation $s$ for $(h,s)=(13,2)$. We discuss the differences between derivative approximation and test function convolution in Subsection~\ref{sec:WeakvsStrong}.  
For these different methods, we consider  following analytical aspects to choose the best method for S-IDENT.

\vspace{0.2cm}
\textbf{Moment vanishing and derivative approximation:}
Assume that $g$ is $p$-times differentiable and $g^{(p)}$ is continuous for some integer $p > s$. By Taylor's theorem with the remainder in Lagrange form, we have an integral expansion of \eqref{eq_general_filter} as
\begin{align}
&\cT_{\mu_s}^{h}(g)(x)= \sum_{j=0}^{p-1}\left(\frac{h^{j-s}g^{(j)}(x)}{j!}\int_{-1}^1t^j\,d\mu_s(t)\right) + \frac{h^{p-s}}{p!}\int_{-1}^1 g^{(p)}(\zeta_{x,h})t^p \,d\mu_s(t),
\label{eq_error_decomposition_clean}
\end{align}
where $\zeta_{x,h}$ is between $x-h$ and $x+h$.
For any integer $r$ satisfying $p\geq r>s\geq 0$, it is clear from the integral expansion \eqref{eq_error_decomposition_clean} that if $\mu_s$ has the following \textit{moment vanishing} property:
\begin{equation}\label{eq_moment_vanishing}
\int_{-1}^1t^j\,d\mu_s(t) = \begin{cases}
 0,\quad &j \in\{0,1,\dots,r-1\}\setminus\{s\}\\
s!,\quad &j=s
\end{cases}\;,
\end{equation}
then we have
\begin{equation*}
g^{(s)}(x) = \cT^h_{\mu_s} (g) (x) +\cO(h^{r-s}),
\end{equation*} 
which means that the weighted integral form~\eqref{eq_general_filter} approximates the $s$-th derivative with a truncation error of order $h^{r-s}$.

\vspace{0.2cm}
\textbf{Effects of noise:} For noisy data, consider 
\begin{equation*}
  \widetilde{g}(x) := g(x) + \nu(x),
\end{equation*} 
where $\nu$ is an uncorrelated Gaussian process with mean $0$ and variance $\sigma^2>0$. It is easy to see that if $\mu_s$ is of bounded variation and satisfies the moment vanishing property~\eqref{eq_moment_vanishing}, then as $h\to0^+$, $\cT^h_{\mu_s} (\widetilde{g})(x)$ is an asymptotically unbiased estimator for $g^{(s)}(x)$ for any $x\in\mathbb{R}$. 
The mean squared error (MSE) is
\begin{equation}\label{eq_variance_diagonal_mass}
\bbE\left(\cT^h_{\mu_s} (\widetilde{g})(x)-g^{(s)}(x)\right)^2= \frac{1}{h^{2s}} \int_{-1}^1\int_{-1}^1\bbE\left(\nu(x+yh)\nu(x+zh)\right)\,d\mu_{s}(y)\,d\mu_{s}(z) + \cO(h^{2(r-s)}).
\end{equation}
 This shows that the weighted integral form \eqref{eq_general_filter} for the $s$-th order derivative amplifies the noise in the data by (i) the power of the scale parameter $h^{-2s}$ in the discrete case and $h^{-2s-1}$ in the continuous case, and (ii) the diagonal mass of the product measure. In the case of discrete filters, the integral in~\eqref{eq_variance_diagonal_mass} becomes the summation of the squared weights. 
In any case, bias decreases while the first term in~\eqref{eq_variance_diagonal_mass} increases as $h\to0^+$, reflecting the general principle of the trade-off between bias and variance.

\begin{figure}
\centering
\begin{tabular}{c@{\vspace{2pt}}c@{\vspace{2pt}}c}
(a)&(b)&(c)\\
\includegraphics[width=0.33\textwidth,height=2in]{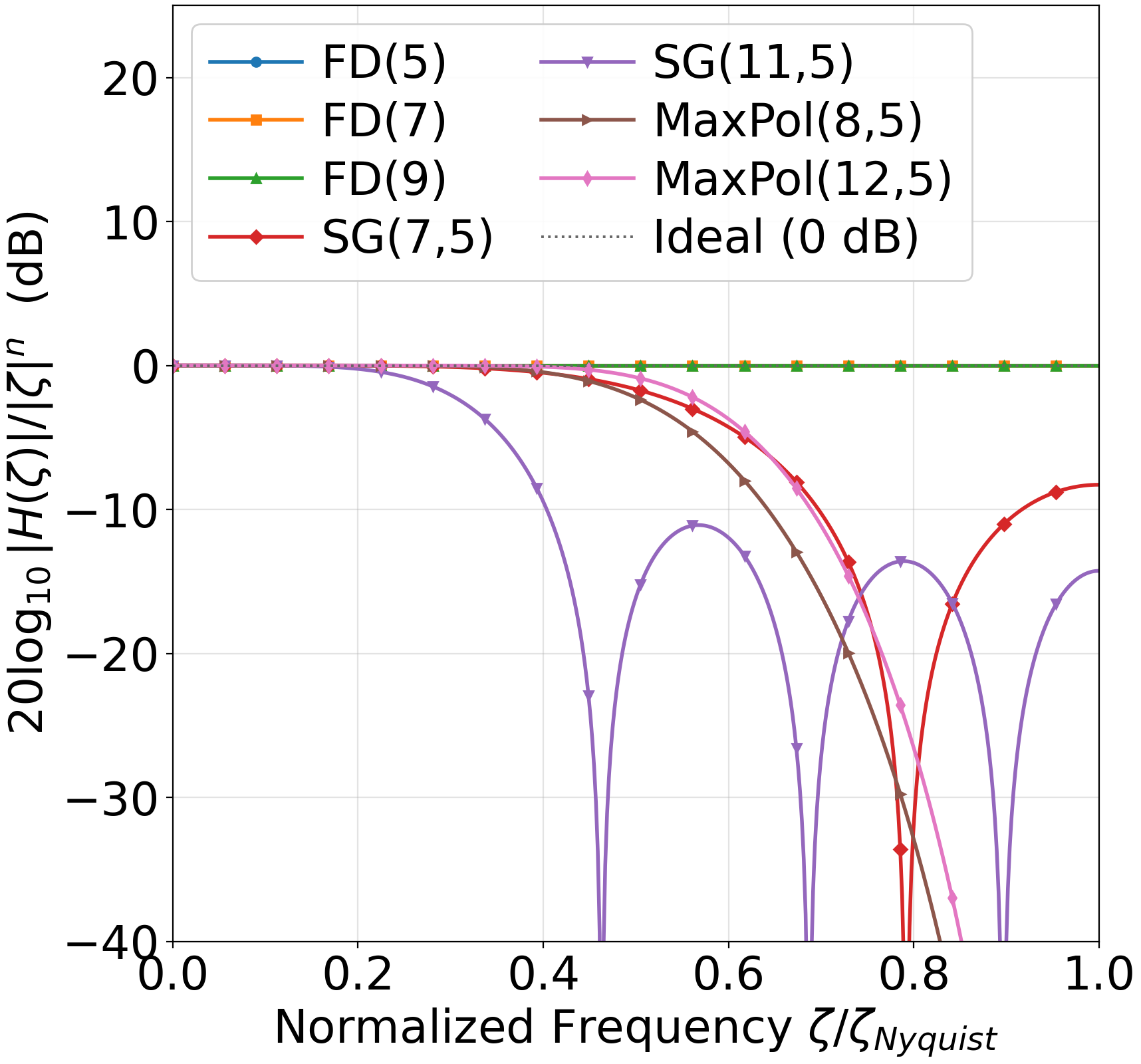}&
\includegraphics[width=0.33\textwidth,height=2in]{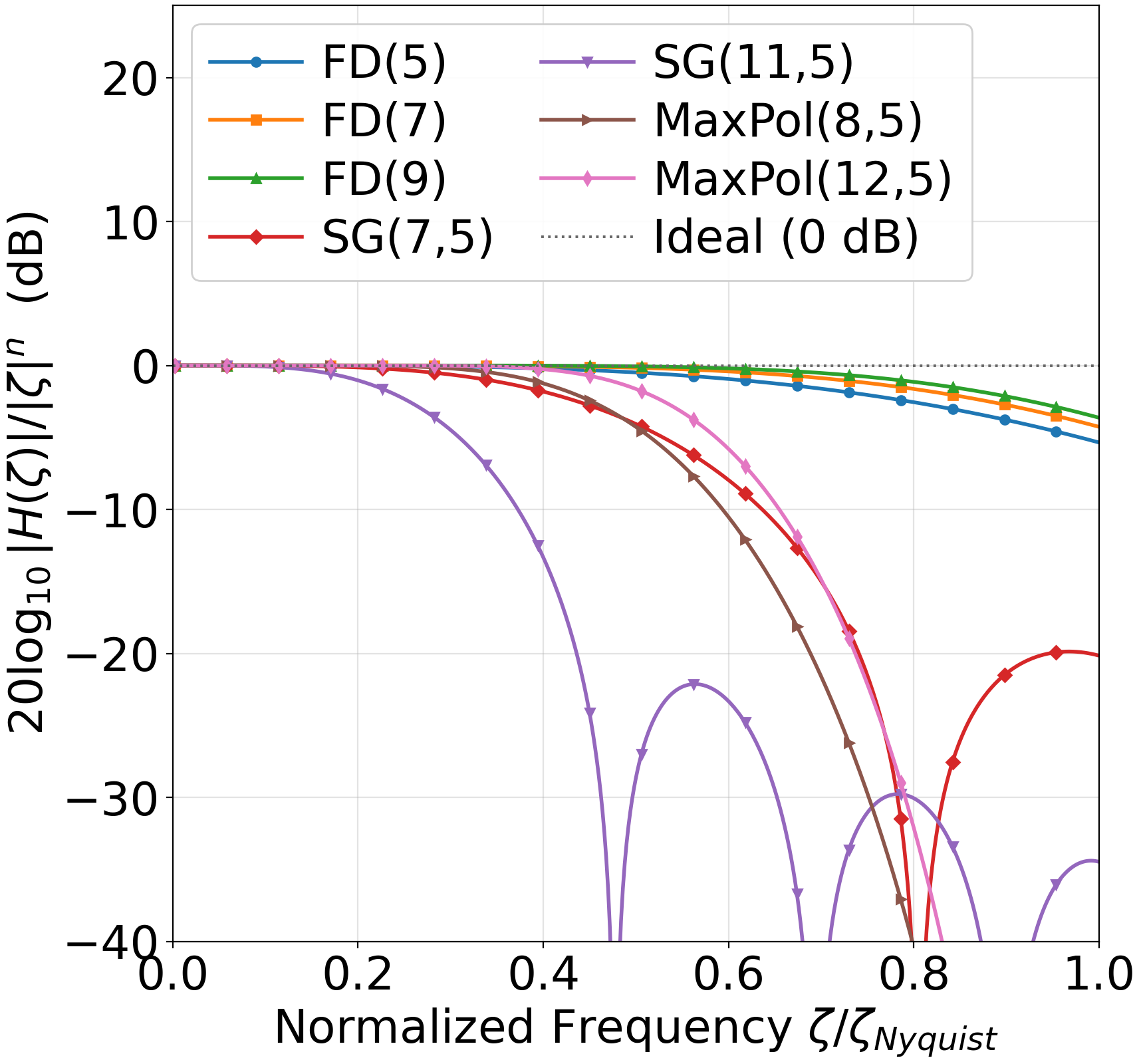}&
\includegraphics[width=0.33\textwidth,height=2in]{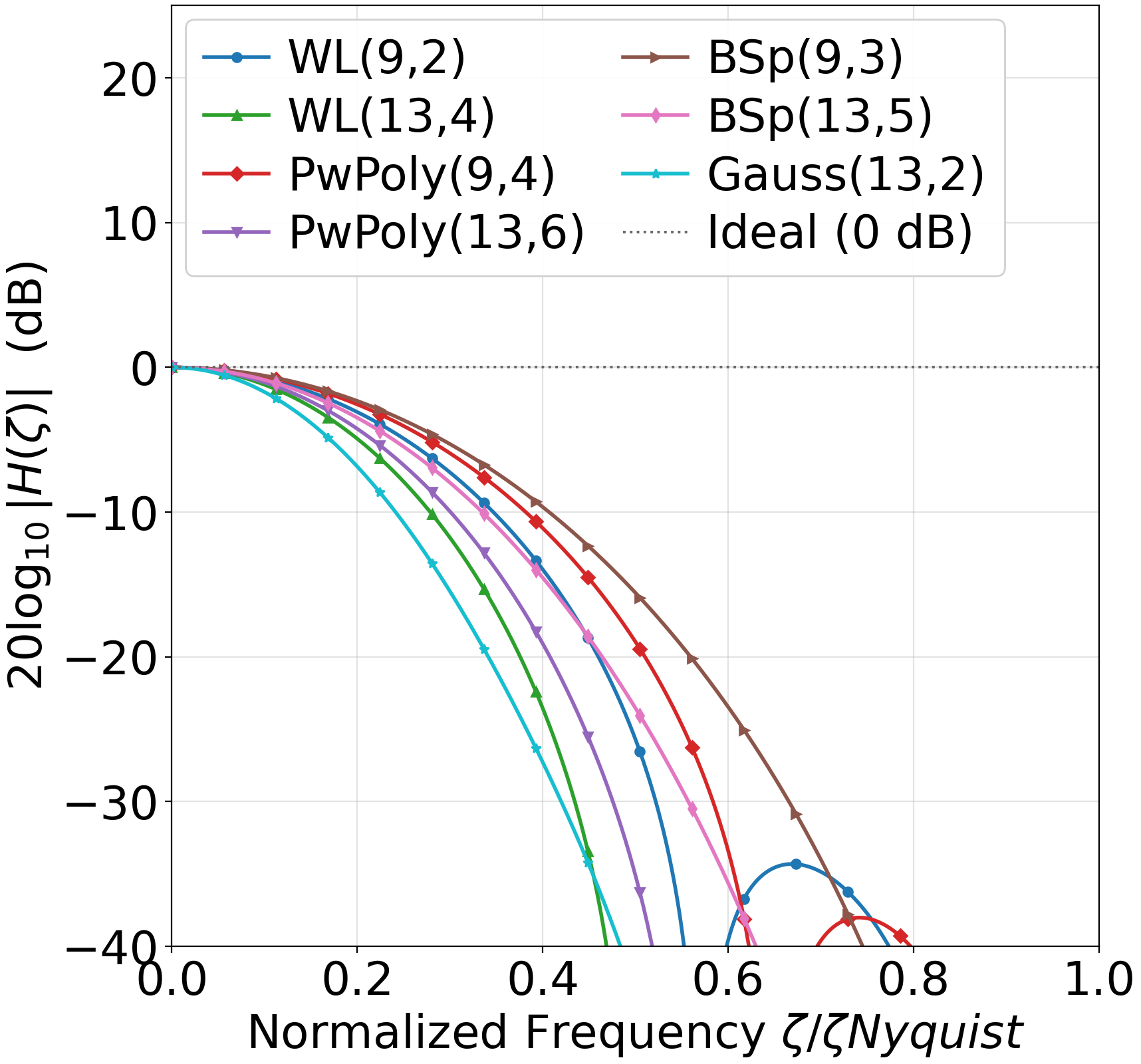}
\end{tabular}
\caption{Relative magnitude response \eqref{eq:magnitudeR} of the representative methods shown in Figure~\ref{fig:filter-shape}. (a) Effects of smoothing,  (b) 2nd-order differentiation, and (c) kernel test functions. The decay of the relative magnitude response at higher frequencies reflects the denoising effect, at the cost of reduced accuracy in approximating the ideal filters.  Here $\zeta_{\mathrm{Nyquist}}=1/(2h)$ is the Nyquist frequency. 
}\label{fig:mag-response}
\end{figure}

\vspace{0.2cm}
\textbf{Fourier domain:}
Examining the weighted integral \eqref{eq_general_filter} in the frequency domain provides insights into its global effects. For simplicity, we focus on the case where $g\in L^2(\mathbb{R})\cap C^p(\bbR)$, and apply the Fourier transform 
$\mathfrak{F}: L^2(\mathbb{R})\to L^2(\mathbb{R})$ to~\eqref{eq_general_filter}, 
\begin{equation*}
\mathfrak{F}(\cT_{\mu_s}^{h}(g))(\zeta) = \mathfrak{F}(g)(\zeta)\cdot h^{-s}\int_{-1}^1e^{2\pi i th\zeta}\,d\mu_s(t)\;,\quad \zeta\in\mathbb{R}.
\end{equation*}
This formula shows that the weighted integral~\eqref{eq_general_filter} modifies the spectral content of $g$ via point-wise multiplication by a factor dependent on the Fourier transform of the measure $\mu_s$. In particular, the factor
\begin{equation*}
H^h_{\mu_s}(\zeta):=h^{-s}\int_{-1}^1e^{2\pi i th\zeta}\,d\mu_s(t)\;,\quad\zeta\in\mathbb{R}\;,
\end{equation*}
is known as the \textit{frequency response} associated with the function $\mu_s$ and scale $h$, and its modulus signifies the amplification or attenuation of the frequency component of $g$ at $\zeta$. 
We further define the \textit{relative magnitude responses} on a log scale as 
\begin{equation}\label{eq:magnitudeR}
    M^h_{\mu_s}(\zeta):=20 \log_{10}\left(| H^h_{\mu_s}(\zeta) |/ |2\pi\zeta |^s\right)
\end{equation}
where $s$ is the order of the derivative. In Figure~\ref{fig:mag-response}, we show the relative magnitude responses of the methods shown in Figure~\ref{fig:filter-shape} at frequency $\zeta$.  Notice that $\zeta\mapsto \left(2\pi i\zeta\right)^s$ is the frequency response of the ideal $s$-th order differentiation, which clearly shows that oscillations are amplified during differentiation, with higher frequencies experiencing greater amplification.  When $M_{\mu_s}^h(\zeta)=0$, the integral operator~\eqref{eq_general_filter} matches the ideal differentiator's gain at frequency $\zeta$. When $M_{\mu_s}^h<0$, the $\zeta$-frequency is suppressed, and when $M_{\mu_s}^h>0$, it is amplified compared to the ideal filter.   
To remove oscillations caused by noise, the function $\mu_s$ can be constructed so that $H_{\mu_s}^h(\zeta)$ is close to $0$ or $M_{\mu_s}^h(\zeta)$ is close to $-\infty$ for large $\zeta$.

Interestingly, the moment vanishing property~\eqref{eq_moment_vanishing} translates, in the frequency domain, into derivative conditions on the frequency response at $\zeta=0$, as follows:
\begin{equation}\label{eq_frequency_condition}
\frac{d^j}{d\zeta^j}\bigg|_{\zeta=0}H_{\mu_s}^{h}(\zeta)=\begin{cases}
0&\quad j \in\{0,1,\dots,r-1\}\setminus\{s\}\\
(2\pi i)^{s} s!&\quad j=s
\end{cases}\;.
\end{equation} 
Note that the derivatives of the frequency response at $0$ dictate the consistency of the corresponding approximation for the function derivatives, while its values at frequencies away from $0$ can be designed for different purposes.

\subsection{Local polynomial interpolation: Fornberg's finite difference}\label{sec:fornberg}

A classical approach to approximate  derivatives from discrete data is computing the derivatives of locally interpolating polynomials. Given $N+1$ uniform interpolating points $\{x_0,x_1,\dots, x_N\}$ with grid size $\Delta x$ and noisy observations $y_n$ of $g(x_n)$ for $n=0,1,\dots,N$, the Lagrange interpolation polynomial is defined by
\begin{equation}\label{eq_fornberg}
\psi(x) = \sum_{n=0}^N \ell_{n}(x)y_n\;,\quad\text{where}
~\ell_n(x):= \frac{\omega_N(x)}{\omega'_N(x_n)(x-x_n)}
\end{equation}
and $\omega_N(x) := \prod_{n=0}^N(x-x_n)$. For any integer $s$ satisfying $0\leq s\leq N$, the $s$-th order derivative of $\psi$ is used to approximate that of the underlying function:
\begin{equation}\label{eq_weight_Fornberg}
g^{(s)}(x) \approx \psi^{(s)}(x) = \sum_{n=0}^N \ell^{(s)}_{n}(x)y_n\;,
\end{equation}
which is a form of a weighted sum~\eqref{eq_linear_filter}. 
Note that $\psi$ is a polynomial of degree $N$ uniquely determined by the provided $N+1$ data points. The functions $\{\ell_n:n=0,1,\dots,N\}$ form a  basis for polynomials of degree $\leq N$, and their derivatives at the interpolating points in~\eqref{eq_weight_Fornberg} can be efficiently computed via Fornberg's recursive algorithm~\cite{fornberg1988generation}. For any $g\in C^{N+1}(\bbR)$, the following error bound holds:
\begin{equation*}
\|g^{(s)} -\psi^{(s)}\|_{\infty}\leq \frac{\|\omega_N^{(s)}\|_{\infty}\|g^{(N+1)}\|_{\infty}}{s!(N+1-s)!}\;,
\end{equation*}
for $s=0,\dots,N$; see~\cite{howell1991derivative} for the proof. 
The factor $\|\omega_N^{(s)}\|_{\infty}= \cO(\Delta x^{N-s+1})$ reflects the influence of the locations of the interpolating points. 

As the interpolation polynomial~\eqref{eq_fornberg} passes through each given data point, the approximate differentiation~\eqref{eq_weight_Fornberg} is sensitive to noise; nevertheless, due to finite sampling, it does exhibit a certain level of smoothing. 
This is illustrated in Figure~\ref{fig:mag-response} (b), where the corresponding relative magnitude response~\eqref{eq:magnitudeR} drops below $0$ as the frequency increases.

\subsection{Local polynomial regression: Savitzky--Golay and SDD}\label{sec:sg}

Instead of local interpolating polynomials, the Savitzky-Golay (SG) method~\cite{savitzky1964smoothing} finds a degree-$d$ polynomial in the least squares sense:
\begin{equation}\label{eq:local_polynomial}
p_n(x) := \sum_{j=0}^d\widehat{b}_j(x_n)(x - x_n)^j
\end{equation}
for each sample point $x_n$ to approximate the function value $g(x_n)$, and the derivatives of $g$ at $x_n$ are approximated by the derivatives of $p_n$ at $x_n$. The local polynomial coefficients $\widehat{\bb}_d(x_n):=(\widehat{b}_0(x_n),\dots,\widehat{b}_d(x_n))$ of~\eqref{eq:local_polynomial} are computed from
\begin{equation}\label{eq_SG_fitting}
\widehat{\bb}_d(x_n)= \argmin_{\bb=(b_0,\dots,b_d)\in\mathbb{R}^{d+1}} \sum_{m=1}^N\left(\left(y_m - \sum_{j=0}^{d}b_j(x_m-x_n)^j\right)^2\cdot \Pi\left(\frac{x_m-x_n}{h}\right)\right).
\end{equation}
Here $\Pi$ is a rectangular window function that takes the value $1$ on $[-1,1]$ and $0$ otherwise; $h>0$ is called the window-length parameter, which controls the number of samples involved in the local polynomial parameter estimation. For every $x_n$, the problem~\eqref{eq_SG_fitting} admits a unique solution obtained via the pseudoinverse of a Vandermonde matrix. In particular, when the sampling is uniform, $x_n = x_0+n\Delta x$ for some $\Delta x>0$ and $h=L\Delta x$ for some integer $1\leq L < N/2$, the derivative of the polynomial $p^{(s)}_n(x_n)$ can be expressed in the form of a weighted sum as in \eqref{eq_linear_filter}: 
\begin{equation}\label{eq_weight_SG}
g^{(s)}(x_n)\approx p^{(s)}_n(x_n)= \sum_{\ell=-L}^{L}C_{d,\ell}^{(s)}y_{n+\ell}\;,~s=0,1,\dots,d,
\end{equation}
for every $L\leq n\leq N-L$. The weights $C_{d,\ell}^{(s)}\in\mathbb{R}$ can be computed efficiently by exploiting the discrete Chebyshev polynomials~\cite{persson2003smoothing}.

Importantly, the SG differentiation~\eqref{eq_weight_SG} satisfies the moment vanishing property~\eqref{eq_moment_vanishing} for $r=d+1$, yielding an asymptotically unbiased estimator of the derivative~\cite{fan1997local}. 
With independent additive noise, it is shown in~\cite{ruppert1994multivariate} that
\[
\bbE\left(p^{(s)}_n(x_n) -g^{(s)}(x_n) \right) = \cO\left(h^{d+1-s}\right)\;,~\operatorname{Var}\left(p^{(s)}_n(x_n)\right) = \cO\left(\Delta x h^{-2s-1}\sigma^2\right),
\]
for $L\leq n \leq N-L$. 
Thanks to the linear form~\eqref{eq_weight_SG}, the frequency properties of SG-based differentiation are well understood~\cite{luo2005properties}. For instance, its \textit{cutoff frequency}, beyond which the relative magnitude response \eqref{eq:magnitudeR} remains significantly low, is proportional to the polynomial degree $d$ and inversely proportional to the window length $h$. We also mention that for the points $x_n$ with $n<L$ or $n>N-L$, depending on the boundary conditions, different approaches can be used~\cite{gorry1990general}. Note that~\eqref{eq_weight_SG} requires an odd number of points ($2L+1$) in the window; in~\cite{luo2005savitzky}, variants of SG derivatives with an even number of data points are considered. 
As discussed in Subsection~\ref{sec:sure_sg}, the SG filter yields the minimum-variance approximation among all linear filters with a guaranteed order of accuracy.

\begin{figure}
\centering
\begin{tabular}{c@{\vspace{2pt}}c@{\vspace{2pt}}c}
(a)&(b)&(c)\\
\includegraphics[width=0.33\textwidth,height=2in]{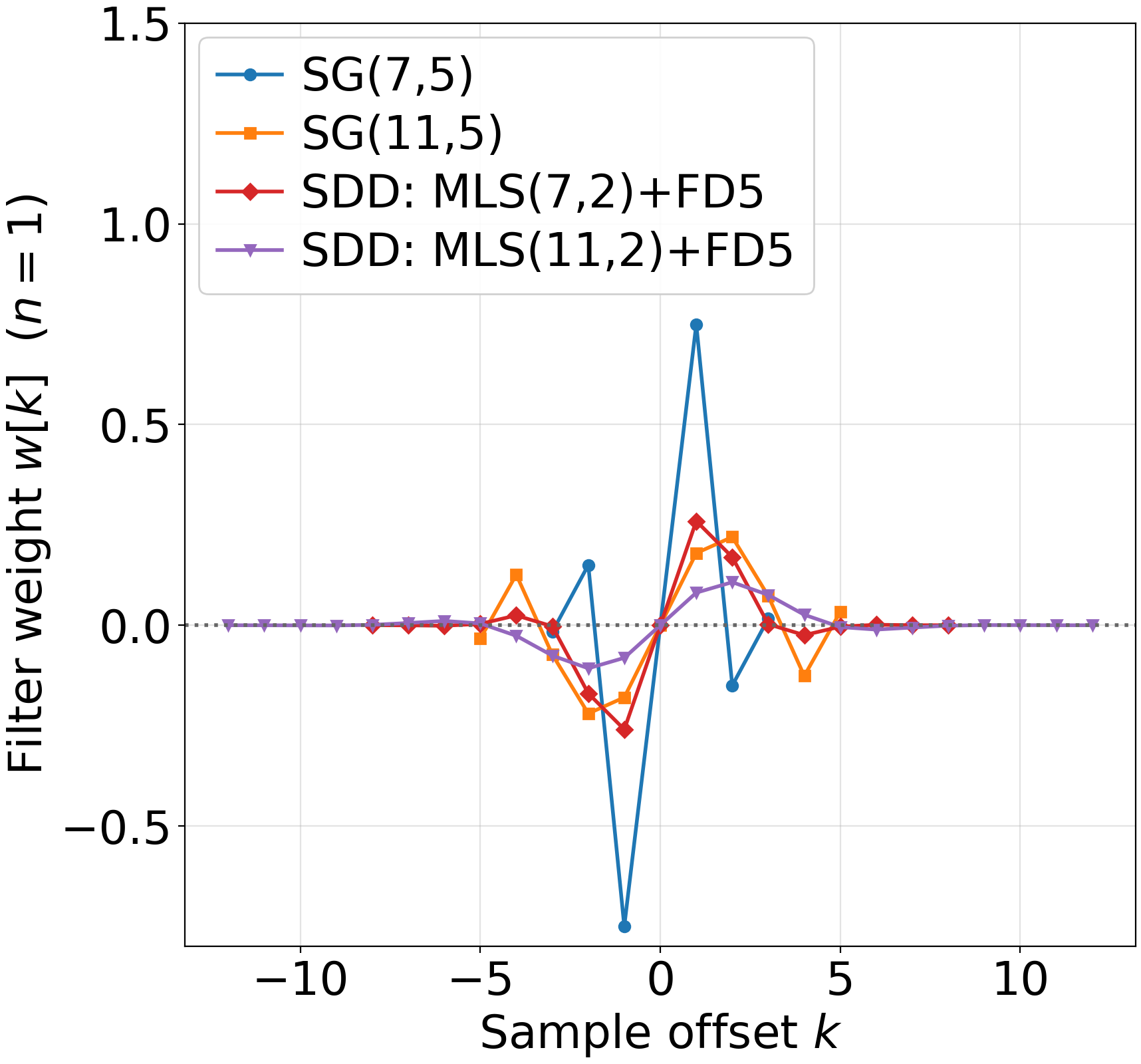}&
\includegraphics[width=0.33\textwidth,height=2in]{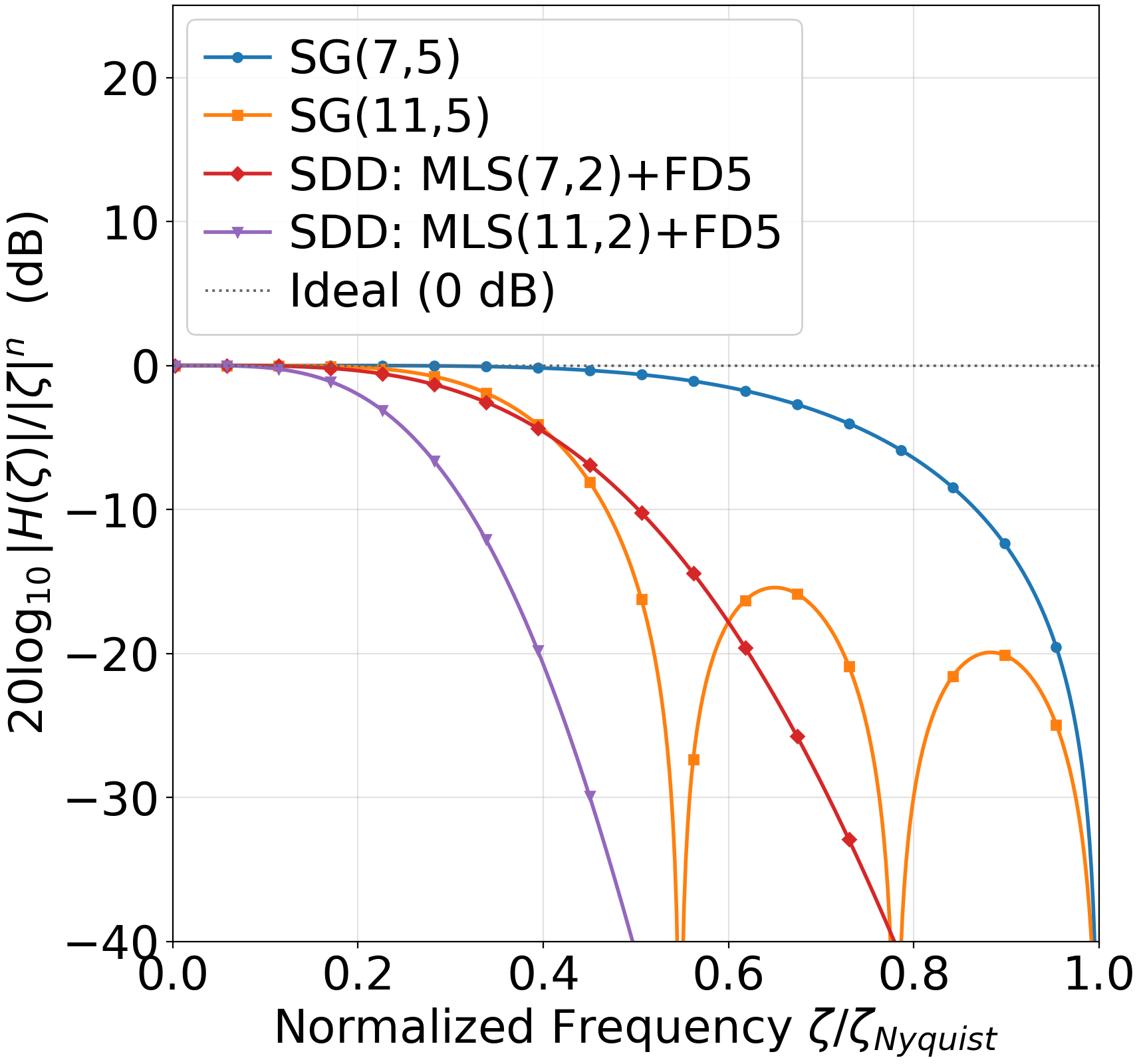}&
\includegraphics[width=0.33\textwidth,height=2in]{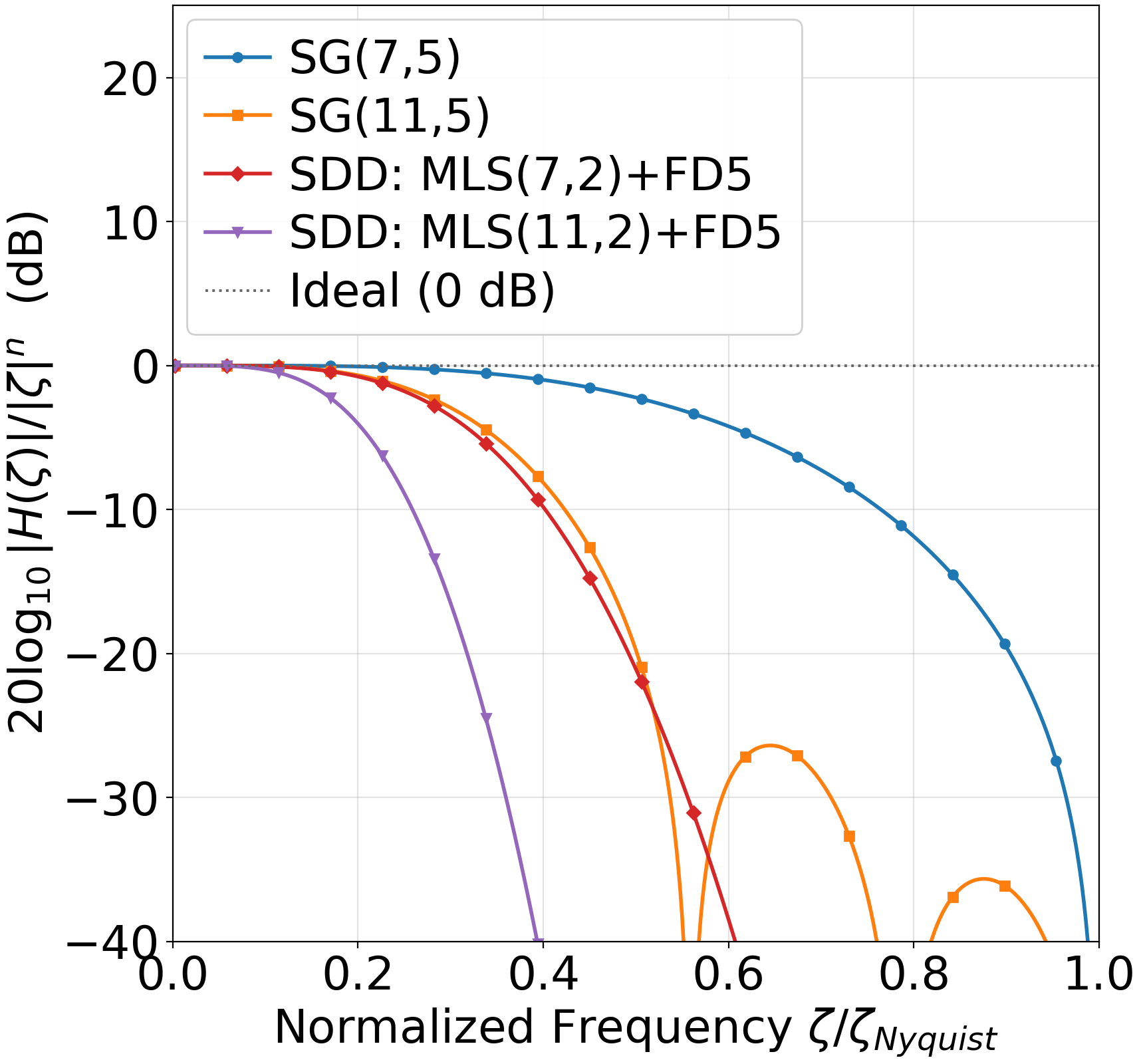}
\end{tabular}
\caption{Comparison between SG differentiation and Successive Denoised Differentiation (SDD)~\cite{he2022robust}. For SDD, MLS weighted by a Gaussian function is used for smoothing, and a $5$-point FD for differentiation. (a) The weights for smoothing in the physical domain; (b) and (c) are relative magnitude responses~\eqref{eq:magnitudeR} for (b) smoothing and (c) second-order derivative approximation.}\label{fig:SG_vs_SDD}
\end{figure}

In a more general context, equation \eqref{eq_SG_fitting} is a moving least squares (MLS) method weighted by a window function, and other weight functions are extensively studied~\cite{fan2018local}. For example, in Robust-IDENT~\cite{he2022robust}, a Gaussian weight function was explored for denoising, and a flexible paradigm called Successive Denoised Differentiation (SDD) was proposed. It  combines smoothing $\cS$, e.g., MLS, with a numerical differentiation scheme $\cD$, e.g., the essentially non-oscillatory (ENO) scheme~\cite{liu1994weighted}: 
SDD smooths the data with $\cS$ each time $\cD$ is applied to estimate the derivative one order higher. 
In Figure~\ref{fig:SG_vs_SDD}, we compare SDD with SG, in both the physical and frequency domains. For SDD, we use MLS with a Gaussian weight:
\begin{equation}\label{eq_SDD_mls}
\widehat{\bb}_d(x_n)= \argmin_{\bb=(b_0,\dots,b_d)\in\mathbb{R}^{d+1}} \sum_{m=1}^N\left(\left(y_m - \sum_{j=0}^{d}b_j(x_m-x_n)^j\right)^2\cdot \exp\left(-\frac{(x_m-x_n)^2}{\theta^2}\right)\right)\;
\end{equation}
for smoothing and a 5-point FD for differentiation, where $\theta>0$ is a hyper-parameter controlling the neighborhood weights. As in Figure~\ref{fig:filter-shape}, we present the weights for smoothing in the physical domain in (a). We compare SG with window size $p$ and polynomial degree $q$ for $(p,q)=(7,5)$ and $(11,5)$, and SDD with MLS~\eqref{eq_SDD_mls} using window size $N$, $\theta=(N-1)/4$, and polynomial degree $d$ with $(N,d)= (7,2)$ and $(11,2)$, with a $5$-point FD.

In Figure~\ref{fig:SG_vs_SDD}, as in Figure~\ref{fig:mag-response}, we also show the relative magnitude responses~\eqref{eq:magnitudeR} for smoothing in (b) and second-order differentiation approximation in (c). Both SG and SDD attenuate high frequencies while staying close to ideal differentiation near zero frequency. Notice that the magnitude response of SG transitions sharply from passband to stopband but shows oscillatory sidelobes, a consequence of the hard truncation $\Pi$ in~\eqref{eq:local_polynomial}; SDD, which instead uses a smooth Gaussian weight, shows the opposite behavior: a wider transition from passband to stopband with suppressed sidelobes. This trade-off between transition sharpness and stopband ripples reflects the well-known uncertainty principle.

\subsection{Frequency domain approach: MaxPol differentiation}\label{sec:maxflat}

Unlike the aforementioned approaches, which design the weight function in the physical domain, the maximally flat principle~\cite{hosseini2017finite} requires that the frequency response of $\mu_s$ satisfies the derivative conditions~\eqref{eq_frequency_condition} at frequency $\zeta=0$ for consistency, and that its derivatives up to a certain order at the Nyquist frequency remain $0$ for high-frequency suppression.

One particular design is the MaxPol filter~\cite{hosseini2017finite}. For some integers $A\geq 0$, $B\geq 0$, and $s$ satisfying $0\leq s\leq A$, the frequency response $H_{s}^h$ of the desired weight $\mu_s$ approximating the $s$-th order derivative needs to satisfy:
\begin{equation}
\begin{cases}
\frac{d^a }{d\zeta^a}\bigg|_{\zeta=0}H_s^h(\zeta) = 0&\quad a \in\{0,1,\dots,A\}\setminus\{s\}\\
\frac{d^a }{d\zeta^a}\bigg|_{\zeta=0}H_s^h(\zeta) =s!\cdot (2\pi i)^{s}&\quad a=s\\
\frac{d^b }{d\zeta^b}\bigg|_{\zeta=(2\Delta x)^{-1}}H_s^{h}(\zeta) = 0\;,&b=0,1,\dots,B\\
\end{cases},\label{eq_system}
\end{equation}
where $\Delta x>0$ is the sampling interval. Since~\eqref{eq_system} is a linear system with $A+B+2$ equations and $N$ unknowns, the number of terms in the weighted sum~\eqref{eq_linear_filter} needs to be $N=A+B+2$. If $N=2L+1$ is odd for some $L\geq 1$, then \eqref{eq_system} can be expressed as
\begin{equation}\label{eq_system_matrix}
\begin{pmatrix}
\mathbf{V}_{A}\\
\mathbf{V}_{B}\mathbf{S}
\end{pmatrix}\bc_s=\bb_s,
\end{equation}
where, for any positive integer $M$,
\begin{equation*}
\bV_M = \begin{pmatrix}
1&1&\cdots&1\\
-L&-L+1&\cdots&L\\
(-L)^2&(-L+1)^2&\cdots&L^2\\
\vdots&\vdots&\ddots&\vdots\\
(-L)^M&(-L+1)^M&\cdots&L^M\\
\end{pmatrix}\in\mathbb{R}^{(M+1)\times N}\;,~\bc_s=\begin{pmatrix}
c_{s,-L}\\
c_{s,-L+1}\\
\vdots\\
c_{s,L}
\end{pmatrix},
\end{equation*}
$\bc_s$ collects the weight coefficients, $\mathbf{S}\in\mathbb{R}^{N\times N}$ is a diagonal matrix with $n$-th diagonal element $\mathbf{S}_{nn} =(-1)^{n}$ for $n=1,\dots,N$, and $\bb_s\in\mathbb{R}^{A+B+2}$ is a zero vector except for the $s$-th entry, which is $(\Delta x)^{-s}\cdot s!$. Once $\bc_s$ is computed from~\eqref{eq_system_matrix}, the MaxPol approximation for the $s$-th derivative is  
\begin{equation}\label{eq_weight_maximal_flat}
g^{(s)}(x_n)\approx  \sum_{\ell=-L}^{L}c_{s,\ell}\cdot y_{n+\ell}\;,~s=0,1,\dots,A.
\end{equation}
As shown in Figure~\ref{fig:mag-response} (a) and (b), the relative magnitude responses of MaxPol remain $0$ near the Nyquist frequency and do not present stopband ripples; this does not necessarily guarantee that the variance of the estimated derivatives will be small. In addition, like SG differentiation, the MaxPol method also allows control over the order of accuracy via $A$ in \eqref{eq_system}, although the computation of the weights $\bc_s$ in \eqref{eq_system_matrix} is more involved.

\subsection{Other methods
} \label{sec:other}

There are other methods for estimating derivatives, and we briefly mention a few important ones to highlight the richness of this field.
In~\cite{lanczos1956applied}, Lanczos suggested that differentiation can be computed via integration and defined
\begin{equation*}
D_hg(x) := \frac{3}{2h^3}\int_{-h}^hg(x+t)t\,dt
\end{equation*}
for small $h>0$. Groetsch proved~\cite{groetsch1998lanczos} that $\lim_{h\to 0^+}D_h g(x)=(D^+g(x)+D^-g(x))/2$ whenever the left and right derivatives $D^-g(x)$ and $D^+g(x)$ exist at $x$. For any $\varepsilon>0$, if $g$ is sufficiently smooth and $h\sim \varepsilon^{1/3}$, then $\|D_hf - g'\|_{\infty}=\mathcal{O}(\varepsilon^{1/3})$ for any $f\in L^2$ with $\|f-g\|_{\infty}\leq \varepsilon$. This mild smoothing of the Lanczos derivative was explained by Burch et al.~\cite{burch2005least}, who found that $\alpha^*=D_hg(x)$ and $\beta^*=\frac{1}{2h}\int_{-h}^hg(x+t)\,dt$ minimize the following error
\begin{equation*}
E(\alpha,\beta):=\int_{-h}^h (g(x+t)-(\alpha t+\beta))^2\,dt\;.
\end{equation*}
This naturally generalizes to higher-order Lanczos derivatives by replacing $\alpha t +\beta$ with higher-degree polynomials; see~\cite{rangarajan2005lanczos} for the details. 
Weighted local polynomial regression is a generalization of~\eqref{eq_SG_fitting} obtained by replacing the rectangular window function with general kernel functions~\cite{cleveland1979robust}. The asymptotic efficiency of these estimators in the mean squared error sense was established in~\cite{fan1993local}, and the quadratic Epanechnikov kernel~\cite{epanechnikov1969non} was found to achieve the greatest efficiency. With general kernels, the estimators can be expressed as a weighted sum, as in~\eqref{eq_linear_filter}, and the corresponding weights are deduced from equivalent kernels.
Designing finite impulse response filters is important in signal processing~\cite{parks1987digital}. Their construction mostly relies on specific constraints on filter parameters such as peak passband and stopband ripples as well as practical requirements such as causality. The Parks–McClellan algorithm~\cite{parks1972chebyshev} provides a powerful technique for their construction. 
Furthermore, a variety of methods exist, including those based on regularization~\cite{nayak2020new}, splines~\cite{unser2002b}, Lyapunov functions~\cite{polyakov2014homogeneous}, and wavelets~\cite{bozzini2003numerical}.

\subsection{Weak versus strong form in terms of convolution }\label{sec:WeakvsStrong}

The discussion above offers a new perspective that reveals similarities and differences
between PDE identification using the weak-form~\cite{messenger2021weak,tang2023weakident} and the differential-form~\cite{kang2021ident,he2022robust}.
When $\mu_s$ is absolutely continuous, the weighted integral form~\eqref{eq_general_filter} for approximating derivatives can be written as
\begin{equation}\label{eq:conv-form-denoise}
    \cT_{\mu_s}^h(g)(x) = \int_{-h}^{h} g(x-y)\, \varphi^h_{\mu_s}(y)\, dy,
\end{equation}
where $\varphi^h_{\mu_s}(y) = h^{-1} \mu_s'(-h^{-1} y)$. If $\varphi_{\mu_s}^h$ vanishes outside $(-h, h)$, the denoised differentiation is formally accomplished by convolving the given data with a kernel function, which may vary depending on the targeted derivative order, the desired accuracy, and the desired regularity.
In the weak formulation~\cite{messenger2021weak,tang2023weakident}, if we take a family of test functions induced by shifting a fixed kernel $\psi_s^h$ compactly supported in $(-h, h)$, then the weak-form features are also of convolution type:
\begin{equation}\label{eq:weak-feature}
    \cW_{\psi^h_s}(g)(x) := \int_{-h}^{h} g(x-y)\, \psi_s^h(y)\, dy.
\end{equation}
Both the denoised differentiation and the weak-form feature construction can be regarded as kernel convolutions with the given data. Thanks to this, both approaches benefit from the Fast Fourier Transform (FFT) for computational efficiency.
In the general weighted integral form~\eqref{eq:conv-form-denoise}, the convolution kernel must satisfy the moment vanishing condition~\eqref{eq_moment_vanishing} to guarantee that the truncation error has the desired order of accuracy, i.e., the condition~\eqref{eq_moment_vanishing} certifies~\eqref{eq:conv-form-denoise} as a valid \textit{approximation of the derivative}. 
In contrast, the weak-form convolution kernel~\eqref{eq:weak-feature} does not need to satisfy this condition in general. 
For example, to construct a weak-form feature corresponding to $g^{(s)}$, one takes
\begin{equation*}
    \psi_s^{h} = (-1)^{s} (\Psi^h)^{(s)},
\end{equation*}
as a test function, where $\Psi^h$ is a kernel compactly supported in $(-h, h)$ that acts as a low-pass filter. When the underlying function lacks sufficient regularity, this construction is related to the weak derivative, and the weak form is indeed indispensable when identifying PDEs from discontinuous solutions~\cite{messenger2021weak,tang2023weakident}. With sufficient regularity, this yields the derivative of a smoothed function, which does not necessarily approximate the derivative of the original function. 
Yet, for Type-W terms, thanks to their special structure, all spatial derivatives can be transferred to the convolution kernel by integration by parts, which converts the differential-form features to the weak-form features~\eqref{eq:weak-feature}. 
The sparsity pattern, i.e., which terms are active or inactive, remains unchanged. 

For differential-form identification such as IDENT, Robust-IDENT and the proposed S-IDENT, once the model is identified, the denoised data approximately satisfies the PDE. 
For the weak form, the observed data approximately satisfies a family of transformed equations associated with the chosen convolution kernels. This distinction highlights the significance of identifying PDEs in their differential form, as the resulting models are not strongly coupled to the choice of test functions.

\section{Numerical Experiments}\label{sec_numerical}

We first present numerical experiments to validate the proposed S-IDENT and compare it with other methods.  In Subsection~\ref{sec:DPspecific}, we present the performance of S-IDENT with the Type-S dictionary for identifying general equations, and in Subsection~\ref{sec:comparison_Stype}, we compare it with methods using Type-S dictionaries, namely SINDy-PDE~\cite{rudy2017data} and Robust-IDENT~\cite{he2022robust}. We compare S-IDENT using Type-S and Type-W dictionaries with the state-of-the-art Weak-IDENT~\cite{tang2023weakident} and Weak-SINDy~\cite{messenger2021weak} in Subsection~\ref{sec:comparison_typeW}. We further investigate the polynomial order for SURE-SG in relation to accuracy in Subsection~\ref{sec:order_accuracy}. In Subsection~\ref{sec:ablation_feature_approx}, we present comparisons of various differentiation approximation methods and explore different strategies for higher-order derivative approximation.

\paragraph{Details of S-IDENT.} We present the pseudo-code of S-IDENT in Appendix~\ref{Asec:Algorithm}. For the hyper-parameters, we set the default values as follows: the SG polynomial degree is fixed at $7$ for the spatial dimension and $5$ for the temporal dimension (see Subsection~\ref{sec:order_accuracy}), $M=20$ for the maximal number of candidate models in~\eqref{eq_sp_sparse}, $\tau = 0.1$ for the trimming~\eqref{eq:trimming}, and $\rho=0.005$ and $L=3$ for the RR model selection~\eqref{eq:cand_score}--\eqref{eq_thresh}.

\paragraph{Noise model.}  
On each simulated trajectory, we add independent centered Gaussian noise as in~\cite{he2022robust, messenger2021weak} with standard deviation $\sigma$ proportional to the root mean square (RMS) of the clean data $\{U_n: n=1,\dots,N\}$:
\begin{equation}\label{eq_nsr_energy}
\sigma  = \NSR \times \sqrt{\frac{1}{N}\sum_{n=1}^N |U_n|^2}\;,
\end{equation}
where $\NSR\geq 0$ denotes the noise-to-signal ratio (NSR). 
Different noise models are used in the related literature, and we discuss and clarify them in Appendix~\ref{noise discussion}. 

\paragraph{Evaluation metrics.}  We primarily evaluate the identification performance in terms of support recovery and coefficient value recovery:
\begin{enumerate}
\item \textit{Support recovery}: Let $S$ denote the identified support and $S^*$ the true support.
The true positive rate and positive predictive value are computed as 
\begin{align}\label{eq:tprppv}
    \mathrm{TPR} := \frac{|S \cap S^*|}{|S^*|}, \quad \text{ and } \quad
    \mathrm{PPV} := \frac{|S \cap S^*|}{|S|},
\end{align}
respectively, where $|A|$ denotes the number of elements in any set $A$.
These values are between $0$ and $1$, and the higher, the better; moreover, when $\mathrm{TPR}=\mathrm{PPV}=1$, we have exact recovery, that is, $S=S^*$. To quantify the frequency of exact recovery from $M$ independent experiments when studying the robustness to noise, we define the exact recovery rate:
\begin{equation}\label{eq:ER}
\mathrm{E.R.} := \frac{|\{S_m: S_m=S^*, m=1,\dots,M\}|}{M}\times 100\%,
\end{equation}
where $S_m$ is the identified support of the $m$-th run.
\item \textit{Coefficient value recovery}: Let $\bc=(c_1,\dots, c_K)$ be the estimated coefficient vector, $\bc^*=(c_1^*,\dots, c^*_K)$ be the true coefficient vector, and $S^*=\{k: c_k^*\neq 0, k=1,\dots,K\}$. We evaluate the coefficient value errors using the following metrics:
\begin{itemize}
\item  Relative in-coefficient error: 
\begin{equation}\label{eq:Ein}
E_{\text{in}}(\bc, \bc^*)=\frac{1}{|S^*|}\sum_{i\in S^*}\frac{|c_i-c^*_i|}{|c_i^*|}\times 100\%.
\end{equation}
It measures the relative error of the coefficients reconstructed for the true features. The smaller, the better, as it indicates accurate coefficient reconstruction on the true support.
\item Relative out-coefficient energy:
\begin{equation}\label{eq:Eout}
E_{\text{out}}(\bc, \bc^*)=\frac{\sum_{i\not\in S^*}|c_i|}{\|\bc\|_1}\times 100\%.
\end{equation}
It quantifies the scale of the coefficients of the wrongly identified features. The smaller, the better, since it measures the relative energy of the wrongly identified coefficients.
\end{itemize}
In addition, we consider the coefficient of determination:
\begin{equation}\label{eq:R2}
  R^{2}(\bc)
  :=
  1-\frac{\bigl\lVert \bm b-\bF\,\bm c\bigr\rVert_2^{2}}
         {\bigl\lVert \bm b-\overline{b}\,\bm 1\bigr\rVert_2^{2}},
  \qquad
  \overline{b}=\frac{1}{N}\sum_{i=1}^{N} b_i,
\end{equation}
which measures the proportion of the variance of the dependent variable that can be explained by the chosen features.

\end{enumerate}

\subsection{S-IDENT for identification of  PDEs}\label{sec:DPspecific}

We present the results of S-IDENT for identifying general PDEs. Using the Type-S dictionary allows us to find a wider range of terms, as illustrated in Table~\ref{tab_dictionary_sizes}.
We consider the following PDEs with periodic boundary conditions:
\begin{itemize}
\item[(A)] Harry Dym equation~\cite{kruskal2005nonlinear}:
\begin{equation}\label{eq:harry-dym}
u_t = u^3\,u_{xxx}, \qquad x\in [0,2\pi),\, t\in [0,0.25],
\end{equation}
with initial condition: a Gaussian bump with a positive background,
        \[
          u(x,0)=0.3+0.7\,\exp\!\Big(-\frac{1}{2}\left(\frac{x-\pi}{\pi/4}\right)^2\Big).
        \]
Data are collected over a regular grid with $256$ points in space and $500$ points in time. Note that in~\eqref{eq:harry-dym}, the dispersion and nonlinearity are coupled. The equation was originally derived from a classical string problem with a varying elastic constant~\cite{hereman1989derivation}.
\item[(B)] Thin film equation~\cite{bernis1990higher}:
\begin{equation}\label{eq:thin-film}
u_t = -\partial_x(u^2\,u_{xxx}) = -2uu_xu_{xxx} -u^2u_{xxxx}, \qquad x\in [0,2\pi),\, t\in [0,0.08],
\end{equation}
with initial condition: two Gaussian droplets,
\[
    u(x,0)=0.05
        +\exp\!\Big(-\frac{1}{2}\left(\frac{x-2\pi/3}{\pi/5}\right)^2\Big)
        +0.6\,\exp\!\Big(-\frac{1}{2}\left(\frac{x-4\pi/3}{\pi/5}\right)^2\Big).
\]
Data are collected over a regular grid with $256$ points in space and $500$ points in time. This is a well-known nonlinear high-order PDE modeling the time evolution of the thickness of a liquid film resting on a surface~\cite{oron1997long}. As the exponent of $u$ is $2$, it describes an intermediate slip condition at the liquid--solid interface. 
\item[(C)] Viscous Hamilton--Jacobi equation~\cite{ben1992global}:
\begin{equation}\label{eq:viscousHJ}
u_t = -\frac{1}{2}\left(u_x^2 + u_y^2\right) + 0.1\Delta u, \qquad (x,y)\in [0,2\pi)\times[0,2\pi),\, t\in [0,0.5],
\end{equation}
with initial condition: random low-frequency Fourier modes,
\[
          u(x,y,0)=\frac{1}{\max|w|}\,w(x,y),
          \qquad
          w(x,y)=\!\!\sum_{\substack{k,l\in\mathbb Z\\ 0<k^2+l^2\le k_{\max}^2}}\!\!
                 \frac{a_{kl}}{|k|+|l|+1}\,
                 \sin\!\Big(kx+ly+\phi_{kl}\Big),
        \]
        with $k_{\max}=5$, $a_{kl}\sim\mathcal U(-1,1)$,
        $\phi_{kl}\sim\mathcal U(0,2\pi)$.
Data are collected over a regular grid with $64\times 64$ points in space and $50$ points in time. This model describes front propagation and is also known as the deterministic Kardar–Parisi–Zhang (KPZ) equation~\cite{kardar1986dynamic}. 
\item[(D)] 2D nonlinear advection (linear gradient):
\begin{equation}\label{eq:2dnladv_v1}
u_t=\frac{1}{2} u\,u_{xx}+\frac{1}{2} u\,u_{yy}+3u_x+2u_y, \qquad (x,y)\in [-\pi,\pi)^2,\, t\in [0,0.2].
\end{equation}
\item[(E)] 2D nonlinear advection (squared gradient):
\begin{equation}\label{eq:2dnladv_v2}
u_t=\frac{1}{2} u\,u_{xx}+\frac{1}{2} u\,u_{yy}+3u_x^2+2u_y^2, \qquad (x,y)\in [-\pi,\pi)^2,\, t\in [0,0.2].
\end{equation}
Both~\eqref{eq:2dnladv_v1} and~\eqref{eq:2dnladv_v2} are simulated with the initial condition:
 \[
          u(x,y,0)=\exp\!\big(\cos x+\frac{1}{2}\cos 3y\big)
                  +\exp\!\big(\frac{1}{2}\cos 2x+\cos y\big).
        \]
\item[(F)] Drinfeld–Sokolov–Wilson (DSW) equation~\cite{drinfeld1981equations}:
\begin{equation}\label{eq:dsw}
\begin{cases}
u_t = 3vv_x\\
v_t = 2v_{xxx} + u_x\,v + 2u\,v_x
\end{cases}, \qquad x\in [0,2\pi),\, t\in [0,0.5],
\end{equation}
with initial condition: sinusoids with two modes
        \[
          u(x,0)=\sin(x)+\frac{1}{2}\sin(2x),
          \qquad
          v(x,0)=\frac{1}{2}\cos(x)+\frac{3}{10}\cos(2x).
        \]
Data are collected over a regular grid with $256$ points in space and $200$ points in time. This model describes solitonic wave interactions in nonlinear dispersive systems and finds applications in fluid dynamics and plasma physics~\cite{patel2025analytical}. 
\end{itemize}

These equations are exclusively representable by Type-S features. For example, in the thin film equation~\eqref{eq:thin-film}, although the outermost differential operator can be transferred to a smooth test function via integration by parts, as done in the weak form, it is not straightforward to address the innermost third-order differentiation in that paradigm.
For the DSW equation~\eqref{eq:dsw}, we remark that although $v v_x = \tfrac{1}{2}(v^2)_x$ and $(uv)_x = u_x v + u v_x$, it cannot be represented by Type-W features, because the coefficient of $u v_x$ differs from that of $u_x v$. This example shows that Type-W dictionaries require specific linear combinations of Type-S features.

\begin{table}
\centering
\small
\begin{tabular}{ccccccc}
\toprule
NSR & TPR~\eqref{eq:tprppv} & PPV~\eqref{eq:tprppv} & $E_{\text{in}}$~\eqref{eq:Ein} & $E_{\text{out}}$~\eqref{eq:Eout} & $R^2$~\eqref{eq:R2} & E.R. (\%)~\eqref{eq:ER} \\
\midrule
\multicolumn{7}{c}{(A) \textit{Harry Dym Equation}~\eqref{eq:harry-dym} ($\Nd=330$)} \\
\midrule
0\% & 1.000 & 1.000 & 0.041 & 0.000 & 1.000 & 100.0 \\
1\% & 1.000 $\pm$ 0.000 & 0.977 $\pm$ 0.116 & 2.570 $\pm$ 0.618 & 0.374 $\pm$ 1.830 & 0.973 $\pm$ 0.004 & 96.0 $\pm$ 19.6 \\
3\% & 0.900 $\pm$ 0.300 & 0.885 $\pm$ 0.313 & 14.592 $\pm$ 28.474 & 11.460 $\pm$ 31.230 & 0.937 $\pm$ 0.009 & 88.0 $\pm$ 32.5 \\
5\% & 0.660 $\pm$ 0.474 & 0.417 $\pm$ 0.352 & 50.637 $\pm$ 36.531 & 42.566 $\pm$ 42.107 & 0.862 $\pm$ 0.015 & 18.0 $\pm$ 38.4 \\
10\% & 0.280 $\pm$ 0.449 & 0.208 $\pm$ 0.371 & 86.938 $\pm$ 26.068 & 80.404 $\pm$ 36.677 & 0.642 $\pm$ 0.061 & 16.0 $\pm$ 36.7 \\
\midrule
\multicolumn{7}{c}{(B) \textit{Thin Film Equation}~\eqref{eq:thin-film} ($\Nd=330$)} \\
\midrule
0\%& 1.000 & 1.000 & 1.188 & 0.000 & 0.998 & 100.0 \\
1\% & 1.000 $\pm$ 0.000 & 0.873 $\pm$ 0.162 & 9.047 $\pm$ 5.197 & 18.004 $\pm$ 23.002 & 0.980 $\pm$ 0.002 & 62.0 $\pm$ 48.5 \\
3\% & 0.370 $\pm$ 0.467 & 0.247 $\pm$ 0.319 & 73.336 $\pm$ 33.526 & 77.318 $\pm$ 31.483 & 0.894 $\pm$ 0.110 & 2.0 $\pm$ 14.0 \\
\midrule
\multicolumn{7}{c}{(C) \textit{Viscous Hamilton-Jacobi Equation}~\eqref{eq:viscousHJ} ($\Nd=816$)} \\
\midrule
0\% & 1.000 & 1.000 & 3.123 & 0.000 & 0.990 & 100.0 \\
1\% & 1.000 $\pm$ 0.000 & 1.000 $\pm$ 0.000 & 3.108 $\pm$ 0.031 & 0.000 $\pm$ 0.000 & 0.975 $\pm$ 0.000 & 100.0 $\pm$ 0.0 \\
3\% & 1.000 $\pm$ 0.000 & 1.000 $\pm$ 0.000 & 3.010 $\pm$ 0.092 & 0.000 $\pm$ 0.000 & 0.874 $\pm$ 0.001 & 100.0 $\pm$ 0.0 \\
5\% & 1.000 $\pm$ 0.000 & 1.000 $\pm$ 0.000 & 2.822 $\pm$ 0.153 & 0.000 $\pm$ 0.000 & 0.722 $\pm$ 0.002 & 100.0 $\pm$ 0.0 \\
10\% & 1.000 $\pm$ 0.000 & 0.981 $\pm$ 0.065 & 2.238 $\pm$ 0.299 & 0.323 $\pm$ 1.131 & 0.394 $\pm$ 0.003 & 92.0 $\pm$ 27.1 \\
20\% & 1.000 $\pm$ 0.000 & 0.961 $\pm$ 0.084 & 7.745 $\pm$ 0.764 & 1.095 $\pm$ 2.843 & 0.132 $\pm$ 0.002 & 82.0 $\pm$ 38.4 \\
\midrule
\multicolumn{7}{c}{(D) \textit{2D Nonlinear Advection (linear gradient)}~\eqref{eq:2dnladv_v1} ($\Nd=816$)} \\
\midrule
0\% & 1.000 & 1.000 & 0.013 & 0.000 & 1.000 & 100.0 \\
1\% & 1.000 $\pm$ 0.000 & 1.000 $\pm$ 0.000 & 0.220 $\pm$ 0.016 & 0.000 $\pm$ 0.000 & 0.995 $\pm$ 0.000 & 100.0 $\pm$ 0.0 \\
3\% & 1.000 $\pm$ 0.000 & 1.000 $\pm$ 0.000 & 1.983 $\pm$ 0.059 & 0.000 $\pm$ 0.000 & 0.953 $\pm$ 0.001 & 100.0 $\pm$ 0.0 \\
5\% & 1.000 $\pm$ 0.000 & 1.000 $\pm$ 0.000 & 5.202 $\pm$ 0.110 & 0.000 $\pm$ 0.000 & 0.878 $\pm$ 0.001 & 100.0 $\pm$ 0.0 \\
10\% & 0.995 $\pm$ 0.035 & 0.992 $\pm$ 0.039 & 16.573 $\pm$ 2.071 & 0.432 $\pm$ 2.116 & 0.630 $\pm$ 0.017 & 94.0 $\pm$ 23.7 \\
20\% & 0.490 $\pm$ 0.049 & 0.884 $\pm$ 0.198 & 56.545 $\pm$ 5.666 & 4.560 $\pm$ 9.868 & 0.196 $\pm$ 0.033 & 0.0 $\pm$ 0.0 \\
\midrule
\multicolumn{7}{c}{(E) \textit{2D Nonlinear Advection (squared gradient)~\eqref{eq:2dnladv_v2}} ($\Nd=816$)} \\
\midrule
0\%& 1.000 & 1.000 & 0.649 & 0.000 & 0.997 & 100.0 \\
1\% & 1.000 $\pm$ 0.000 & 1.000 $\pm$ 0.000 & 0.413 $\pm$ 0.024 & 0.000 $\pm$ 0.000 & 0.989 $\pm$ 0.000 & 100.0 $\pm$ 0.0 \\
3\% & 1.000 $\pm$ 0.000 & 1.000 $\pm$ 0.000 & 3.509 $\pm$ 0.090 & 0.000 $\pm$ 0.000 & 0.933 $\pm$ 0.001 & 100.0 $\pm$ 0.0 \\
5\% & 0.860 $\pm$ 0.124 & 0.976 $\pm$ 0.065 & 22.217 $\pm$ 10.751 & 0.615 $\pm$ 1.665 & 0.774 $\pm$ 0.056 & 32.0 $\pm$ 46.6 \\
10\% & 0.745 $\pm$ 0.169 & 0.898 $\pm$ 0.183 & 42.710 $\pm$ 18.051 & 6.640 $\pm$ 14.346 & 0.511 $\pm$ 0.039 & 20.0 $\pm$ 40.0 \\
20\%& 0.410 $\pm$ 0.120 & 0.650 $\pm$ 0.198 & 79.509 $\pm$ 15.160 & 20.568 $\pm$ 25.177 & 0.182 $\pm$ 0.025 & 0.0 $\pm$ 0.0 \\
\midrule
\multicolumn{7}{c}{(F) \textit{Drinfeld\textendash{}Sokolov\textendash{}Wilson Equation}~\eqref{eq:dsw} ($\Nd=3{,}060$)} \\
\midrule
0\% & 1.000 & 1.000 & 0.061 & 0.000 & 1.000 & 100.0 \\
1\% & 0.983 $\pm$ 0.050 & 1.000 $\pm$ 0.000 & 2.259 $\pm$ 5.273 & 0.000 $\pm$ 0.000 & 0.989 $\pm$ 0.002 & 90.0 $\pm$ 30.0 \\
3\% & 0.997 $\pm$ 0.023 & 0.931 $\pm$ 0.108 & 4.324 $\pm$ 2.954 & 5.433 $\pm$ 9.487 & 0.956 $\pm$ 0.003 & 70.0 $\pm$ 45.8 \\
5\% & 0.980 $\pm$ 0.054 & 0.915 $\pm$ 0.126 & 9.653 $\pm$ 5.085 & 7.605 $\pm$ 12.843 & 0.890 $\pm$ 0.008 & 64.0 $\pm$ 48.0 \\
10\% & 0.520 $\pm$ 0.064 & 0.518 $\pm$ 0.059 & 49.767 $\pm$ 2.729 & 49.404 $\pm$ 2.213 & 0.754 $\pm$ 0.044 & 0.0 $\pm$ 0.0 \\
\bottomrule
\end{tabular}
\caption{Performance of S-IDENT for identifying general PDEs (Type-S). 
For each PDE and each NSR, the mean $\pm$ standard deviation over $50$ independent trials is shown. (A) and (B) are 1D PDEs with a dictionary of size $330$, (C)--(E) are 2D PDEs with a dictionary of size $816$, and (F) is a PDE system with a dictionary of size $3{,}060$. In most cases, even with up to $3\%$ noise, the identification results are satisfactory, with TPR and PPV near $1$ and low coefficient errors.}
\label{tab:group1_all}
\end{table}

For the single-state 1D-space PDEs~\eqref{eq:harry-dym} and~\eqref{eq:thin-film}, we use Type-S $(6,4)$ dictionaries, giving $\mathbf{330}$ features; for the single-state 2D-space PDEs~\eqref{eq:viscousHJ},~\eqref{eq:2dnladv_v1}, and~\eqref{eq:2dnladv_v2}, we use Type-S $(4,3)$ dictionaries, giving $\mathbf{816}$ features; and for the PDE system~\eqref{eq:dsw}, we use the Type-S $(6,4)$ dictionary, giving $\mathbf{3{,}060}$ features. Figure~\ref{fig:solution_group1} in Appendix~\ref{Asec:StypePDE} shows heatmaps of clean trajectories.

We summarize the performance (mean $\pm$ standard deviation) in Table~\ref{tab:group1_all} for different noise levels (NSR), showing TPR~\eqref{eq:tprppv}, PPV~\eqref{eq:tprppv}, $E_{\text{in}}$~\eqref{eq:Ein}, $E_{\text{out}}$~\eqref{eq:Eout}, $R^2$~\eqref{eq:R2}, and E.R. (\%)~\eqref{eq:ER}. 
For each PDE, we test the identification on data with additive Gaussian noise as the NSR~\eqref{eq_nsr_energy} grows, and for each noise level, we conduct 50 independent trials. 
In all cases, we observe that the identifications are successful when there is no noise ($\NSR = 0\%$). 
This is particularly surprising for (F) the DSW equation~\eqref{eq:dsw}, where the single true feature $vv_x$ for $u_t$ can be identified from $3{,}060$ candidates, and only $3$ true features ($v_{xxx}$, $u_x v$, $uv_x$) are selected for $v_t$. 
For (A) the Harry Dym equation~\eqref{eq:harry-dym}, the identification naturally becomes less stable as the NSR increases, yet the mean exact recovery rate remains above $50\%$ until the NSR exceeds $3\%$. 
(B) The thin film equation~\eqref{eq:thin-film}, with its highly nonlinear and high-order features, exhibits sensitivity to noise: the exact recovery rate drops to around $60\%$ at an NSR of $1\%$. 
We observe that the identification for (C) the viscous HJ equation~\eqref{eq:viscousHJ}, which consists of up to second-order derivatives and two product terms, remains the most stable, and the coefficient reconstruction remains accurate. 
The nonlinear advection equations~\eqref{eq:2dnladv_v1} and~\eqref{eq:2dnladv_v2} in (D) and (E) are successfully identified up to $3\%$ noise, even with large dictionary sizes. With $3{,}060$ candidate features, it is notable that the recovery rate of (F) remains above $60\%$ until the NSR grows to $5\%$.

The phenomenon that high-order and strongly nonlinear PDEs are more sensitive to noise  is general. In IDENT~\cite{kang2021ident}, an error analysis explaining this aspect was presented. 
This can also be explained from a frequency perspective. Nonlinear features generate high-frequency components, as multiplication in the physical domain translates to convolution in the frequency domain. Meanwhile, high-order derivatives are most detectable when the signal contains high-frequency components. Since Gaussian noise is white, i.e., uniform across all frequencies up to the Nyquist frequency, and differentiation amplifies high-frequency components, the useful information for accurately identifying the high-order derivatives is largely destroyed. This is a fundamental challenge of identifying PDEs from noisy data.

\subsection{Comparison with methods using general dictionaries (Type-S) 
}\label{sec:comparison_Stype}

We compare S-IDENT with differential-form (strong-form) identification methods: SINDy-PDE~\cite{rudy2017data} and Robust-IDENT~\cite{he2022robust}. These methods use general dictionaries, and we compare the performance on the following classical differential equations using the Type-S dictionary. We note that the following equations (a)--(c) can also be expressed using a Type-W dictionary, and we consistently use these equations for the comparison in the rest of the section. 
\begin{itemize}
\item[(a)] Viscous Burgers equation
\begin{equation}\label{eq:vb}
u_t = auu_x + \nu u_{xx} = \frac{a}{2}(u^2)_x + \nu u_{xx}, \qquad x\in [0,2\pi),\, t\in [0,1.0],
\end{equation}
where $a = -1$ and $\nu = 0.1$.
\item[(b)] KdV equation
\begin{equation}\label{eq:kdv}
u_t = auu_x + b u_{xxx} = \frac{a}{2}(u^2)_x + b u_{xxx}, \qquad x\in [0,2\pi),\, t\in [0,0.5],
\end{equation}
where $a = -1.0$ and $b = -0.1$.
For both~\eqref{eq:vb} and~\eqref{eq:kdv}, we use a random initial condition: 
\begin{equation}\label{eq:initial_cond_random}
u(x,0)=w(x), \text{ where } \;\; w(x)=\sum_{k=1}^{8}\frac{a_k}{k}\,\sin\!\big(kx+\phi_k\big)
\end{equation}
where the magnitude and phase are sampled from a uniform distribution: $a_k\sim\mathcal U(-1,1)$ and $\phi_k\sim\mathcal U(0,2\pi)$. For each experimental run, these coefficients are sampled independently.
\item[(c)] Allen--Cahn equation
\begin{equation}\label{eq:allencahn}
u_t = a  u_{xx} +  bu(1-u^2), \qquad x\in [0,2\pi),\, t\in [0,5.0],
\end{equation}
where $a =0.5$ and $b = 1.0$, with initial condition:
 \[
          u(x,0)=0.5\times\frac{w(x)}{\max_x|w(x)|},
        \]
where $w$ is as given in~\eqref{eq:initial_cond_random}.
After generating the data, we add Gaussian noise with NSR ranging from $0\%$ (no noise) to $50\%$ and repeat the identification for 50 independent random noise samples.
\end{itemize}

\begin{figure}
    \centering
\begin{tabular}{c}
(a) Viscous Burgers equation~\eqref{eq:vb} \\
 \includegraphics[width=\textwidth]{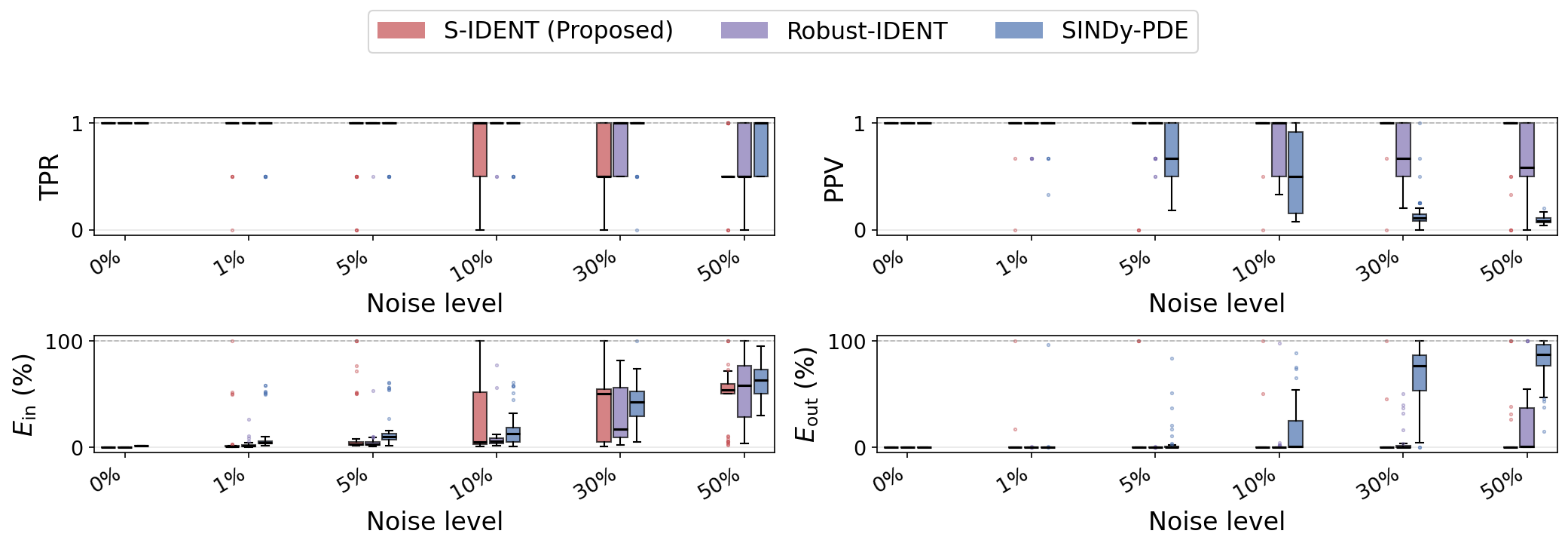}\\
 (b) KdV equation~\eqref{eq:kdv}\\
  \includegraphics[width=\textwidth]{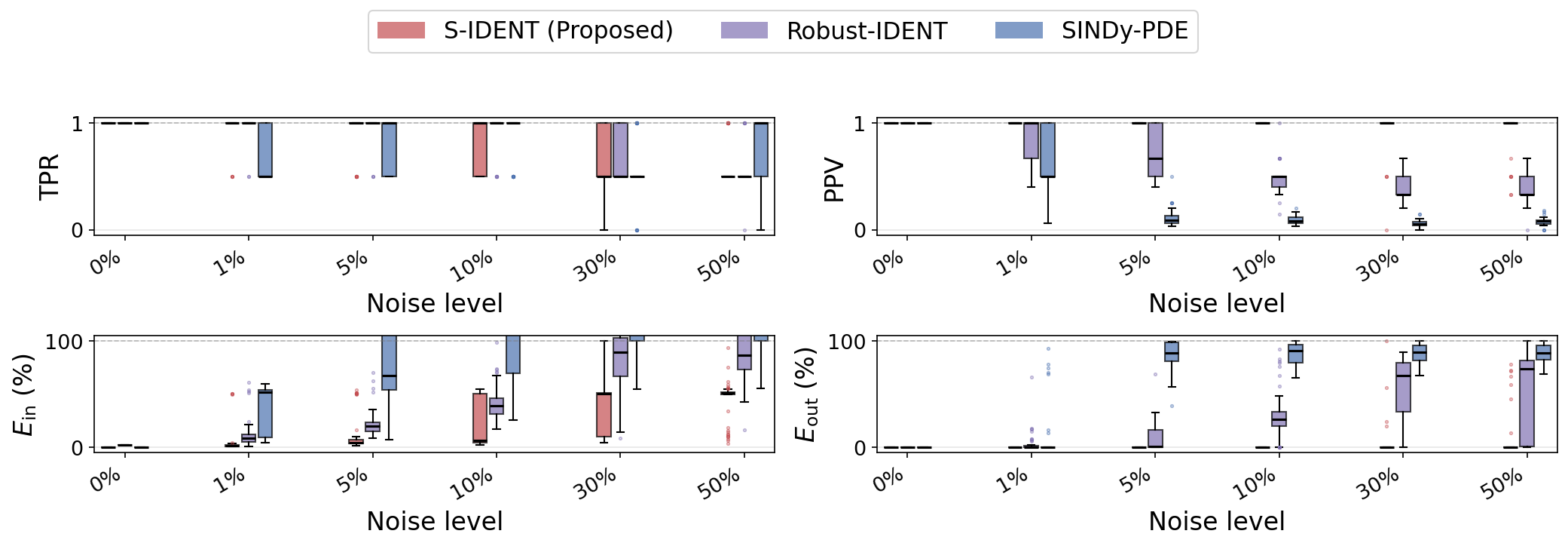}\\
(c) Allen--Cahn equation~\eqref{eq:allencahn}\\
 \includegraphics[width=\textwidth]{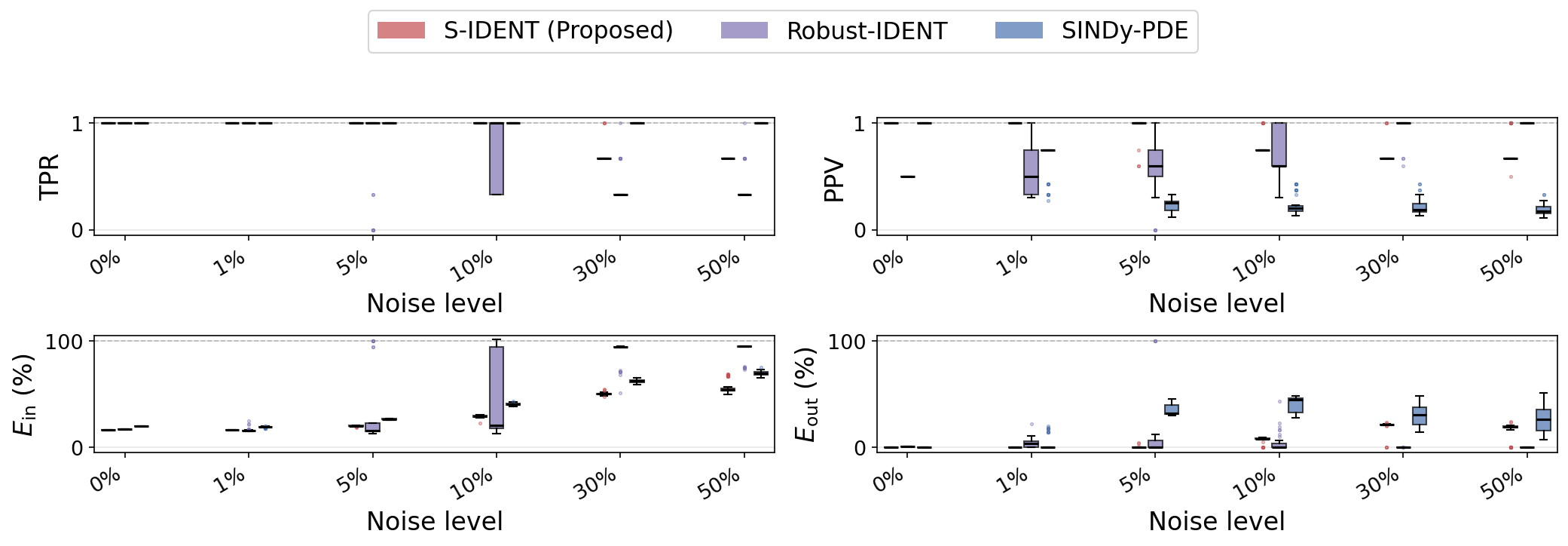}
 \end{tabular}
\caption{Comparison: S-IDENT, Robust-IDENT~\cite{he2022robust}, and SINDy-PDE~\cite{rudy2017data}. 
For each noise level, TPR and PPV are presented in the first row for support identification, and $E_{\text{in}}$ and $E_{\text{out}}$ in the second row for coefficient recovery. Each panel shows box plots from $50$ independent experiments. S-IDENT shows better support and coefficient recovery.}
\label{fig:compare-typeS}
\end{figure}

In Figure~\ref{fig:compare-typeS}, we present comparisons of the proposed S-IDENT, Robust-IDENT~\cite{he2022robust}, and SINDy-PDE~\cite{rudy2017data}.
For all methods, we test with the Type-S $(6,4)$ dictionary containing $330$ features. For both Robust-IDENT and SINDy-PDE, while we mostly use their default parameters, a few parameters need to be tuned for them to achieve exact recovery on clean data in our examples. 
For Robust-IDENT, we also remove  coefficients whose absolute values are very small.  
We present these details of parameter choices for compared methods in Appendix~\ref{sec:compare_details}. 

In Figure~\ref{fig:compare-typeS}, for each noise level, we show box plots for S-IDENT (red), Robust-IDENT (purple), and SINDy-PDE (blue) from $50$ independent experiments. 
Compared to SINDy-PDE and Robust-IDENT, S-IDENT shows significantly improved results. In all these examples, S-IDENT maintains a high TPR across different noise levels and remarkably avoids including incorrect terms when the noise becomes very high, yielding the highest PPV; consequently, the $E_{\mathrm{out}}$ of S-IDENT is the lowest, remaining comparable to that of Robust-IDENT. As the noise level increases, the coefficient reconstruction of the true features by S-IDENT, i.e., $E_{\mathrm{in}}$, remains satisfactory and achieves the lowest level in most cases. This validates the benefit of the noise-adaptive feature approximation enabled by SURE-SG.

\begin{figure}
    \centering
        \centering
 \includegraphics[width=\textwidth]{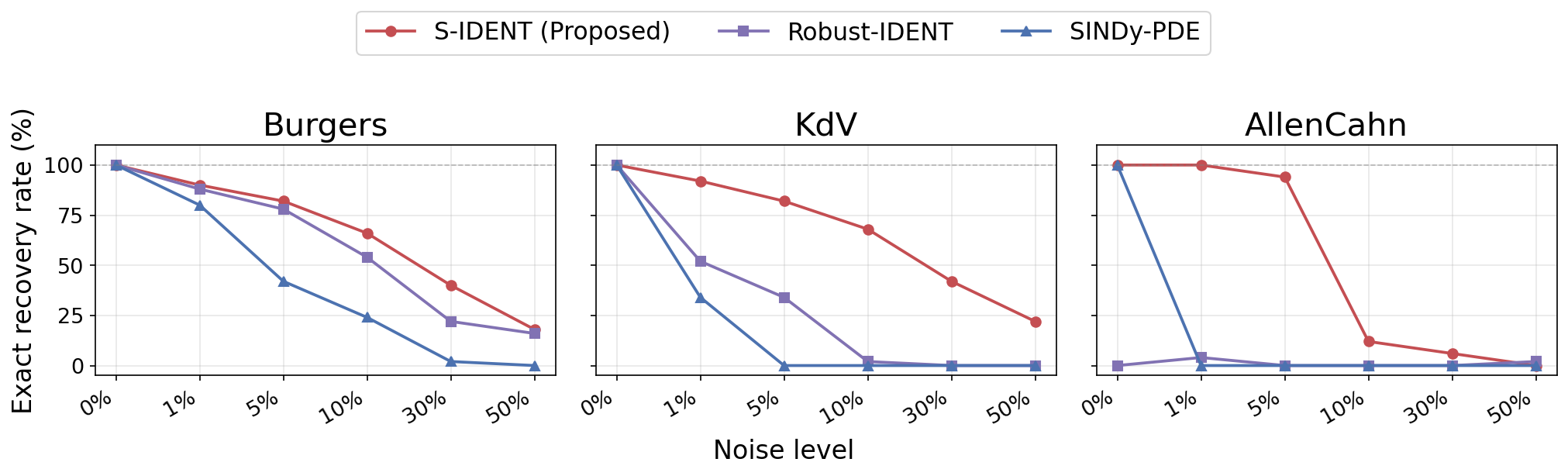}
\caption{Comparison: S-IDENT, Robust-IDENT~\cite{he2022robust}, and SINDy-PDE~\cite{rudy2017data}. Mean exact recovery rate (E.R.)~\eqref{eq:ER} from $50$ independent trials is presented for viscous Burgers~\eqref{eq:vb}, KdV~\eqref{eq:kdv}, and Allen--Cahn~\eqref{eq:allencahn}. S-IDENT shows the best results.}
\label{fig:compare-typeS-ER}
\end{figure}

In Figure~\ref{fig:compare-typeS-ER}, we present the mean E.R.~\eqref{eq:ER}. This confirms the effectiveness of S-IDENT: while the identification by all methods gradually deteriorates as the noise increases, S-IDENT offers a noticeably better chance of exactly recovering the features of the underlying PDEs across most noise levels.

\subsection{Comparison with methods using weak-form dictionaries (Type-W)}
\label{sec:comparison_typeW}

We compare with methods that use the weak-form approach: Weak-SINDy~\cite{messenger2021weak} and Weak-IDENT~\cite{tang2023weakident}, where partial derivatives are transferred to the test functions via integration by parts. In contrast, S-IDENT directly approximates the derivatives from data via SURE-SG. For the weak-form approaches---S-IDENT (W), Weak-SINDy, and Weak-IDENT---we use the Type-W $(6,4)$ dictionary containing $\mathbf{29}$ features. We also add a comparison with S-IDENT (S) using the Type-S $(6,4)$ dictionary to reveal the effectiveness of S-IDENT even when the dictionary size grows from $\mathbf{29}$ to $\mathbf{330}$. For Weak-SINDy and Weak-IDENT, we use the default parameters; see Appendix~\ref{sec:compare_details}.

\begin{figure}
    \centering
\begin{tabular}{c}
(a) Viscous Burgers equation~\eqref{eq:vb} \\
 \includegraphics[width=\textwidth]{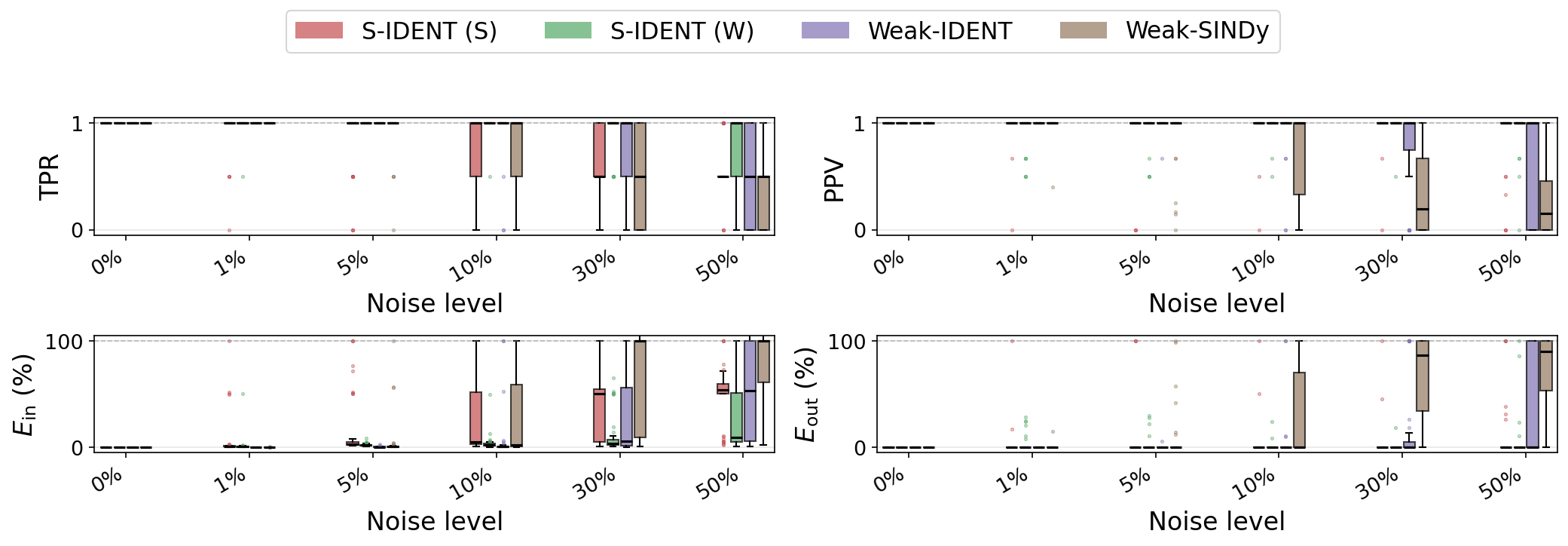}\\
 (b) KdV equation~\eqref{eq:kdv}\\
 \includegraphics[width=\textwidth]{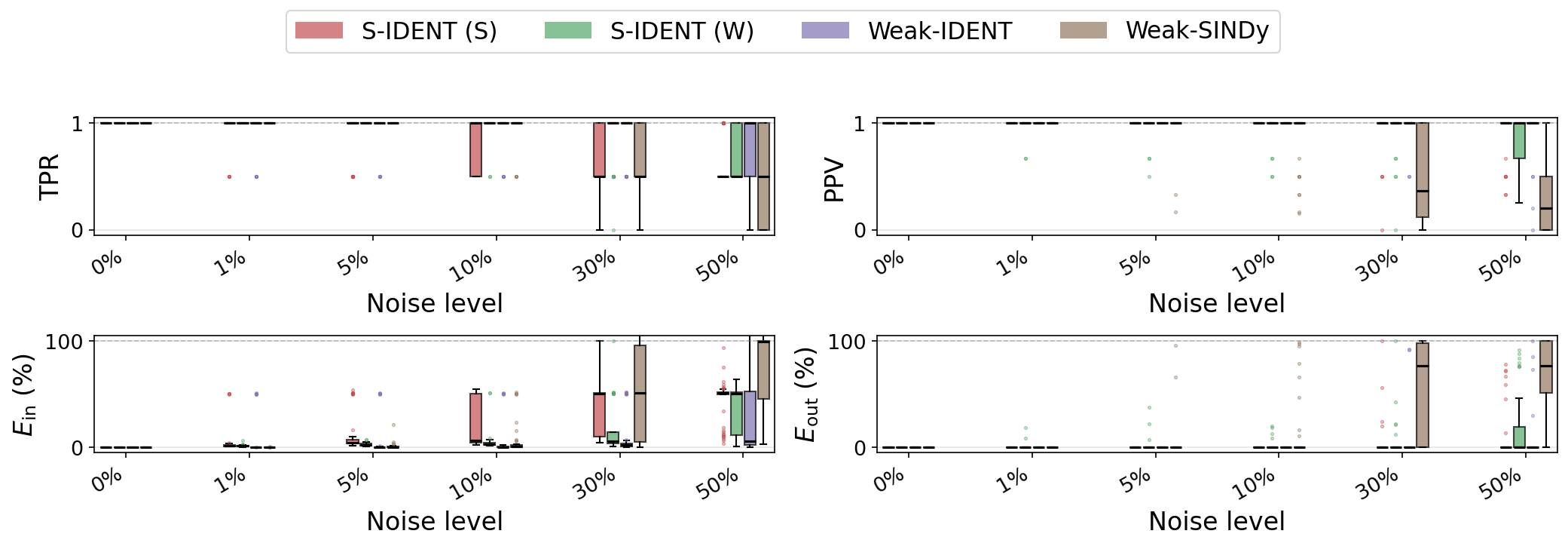}\\
(c) Allen--Cahn equation~\eqref{eq:allencahn}\\
 \includegraphics[width=\textwidth]{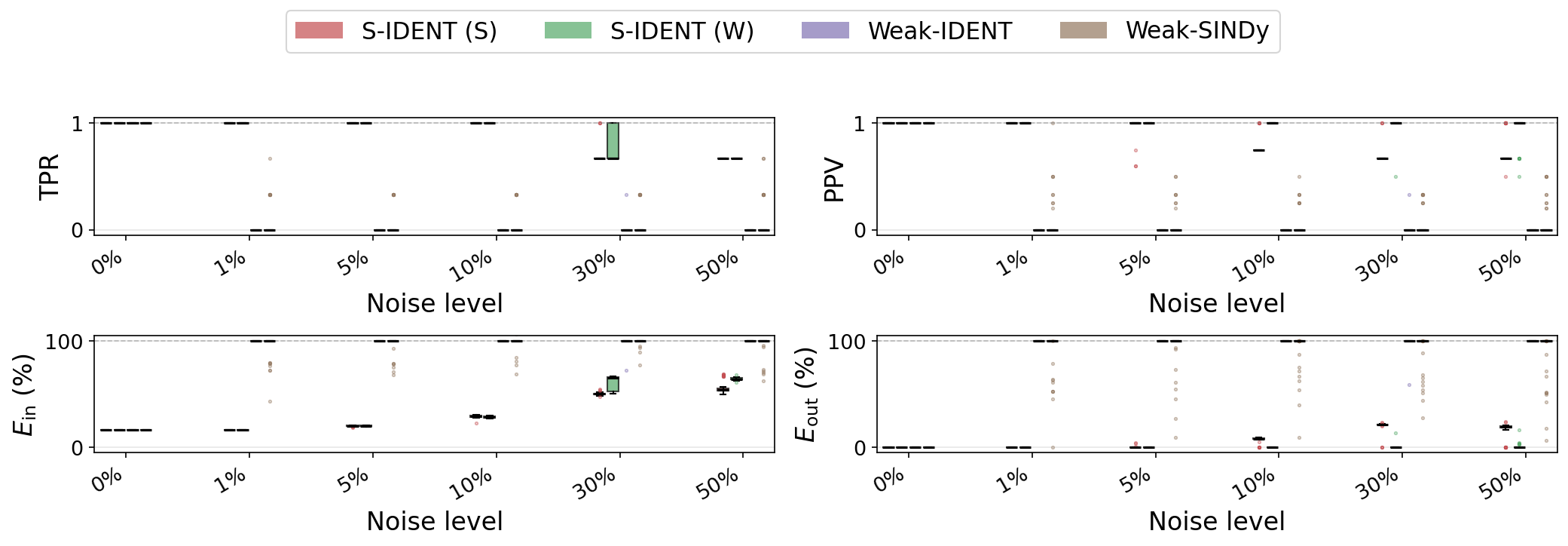}
 \end{tabular}
\caption{Comparison: S-IDENT with Type-S (S-IDENT (S)), S-IDENT with Type-W (S-IDENT (W)), Weak-IDENT~\cite{tang2023weakident}, and Weak-SINDy~\cite{messenger2021weak}. 
For each noise level, TPR and PPV are presented in the first row for support identification, and $E_{\text{in}}$ and $E_{\text{out}}$ in the second row for coefficient recovery. Each panel shows box plots from $50$ independent experiments. S-IDENT (W) shows support and coefficient recovery comparable to or even better than that of Weak-IDENT.}
\label{fig:compare-typeW}
\end{figure}

\begin{figure}
    \centering
        \centering
 \includegraphics[width=\textwidth]{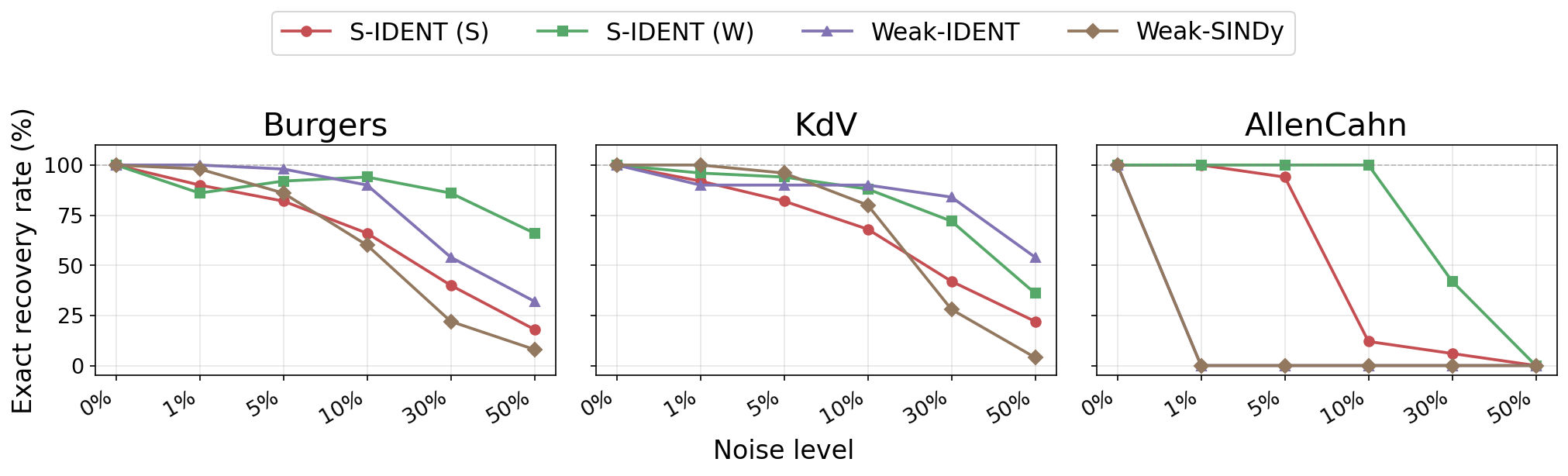}
\caption{Comparison: S-IDENT with Type-S (S-IDENT (S)), S-IDENT with Type-W (S-IDENT (W)), Weak-IDENT~\cite{tang2023weakident}, and Weak-SINDy~\cite{messenger2021weak}.
Mean exact recovery rate (E.R.)~\eqref{eq:ER} from $50$ independent trials is presented for viscous Burgers~\eqref{eq:vb}, KdV~\eqref{eq:kdv}, and Allen--Cahn~\eqref{eq:allencahn}. 
S-IDENT (S) shows results comparable to those of Weak-SINDy, while S-IDENT (W) offers better performance.}
\label{fig:compare-typeW-ER}
\end{figure}

In Figure~\ref{fig:compare-typeW}, we report and compare the performance of the aforementioned methods across $50$ independent experiments. Figure~\ref{fig:compare-typeW-ER} shows the average E.R.~\eqref{eq:ER} of the compared methods. Notably, the results confirm that S-IDENT (W) gives identification accuracy comparable to that of Weak-IDENT, and in the case of Allen--Cahn~\eqref{eq:allencahn}, both S-IDENT (W) and S-IDENT (S) perform better than Weak-IDENT and Weak-SINDy. 
We note that S-IDENT (S) can achieve a level of accuracy similar to that of S-IDENT (W) which uses only $29$ features. For instance, in the case of KdV shown in Figure~\ref{fig:compare-typeW} (b), S-IDENT (S) yields a higher PPV than S-IDENT (W). Figure~\ref{fig:compare-typeW-ER} shows that in the case of KdV, S-IDENT (W) achieves an E.R. similar to that of Weak-IDENT, even though the derivatives are approximated directly using SURE-SG instead of relying on a smooth test function. In the case of the viscous Burgers equation, we observe that Weak-IDENT remains highly competitive until the noise grows to $10\%$, after which S-IDENT (W) maintains an E.R. of approximately $75\%$, the highest among the methods. The most significant results are observed for Allen--Cahn, where S-IDENT with both the Type-S and Type-W dictionaries maintains close to $100\%$ E.R. up to $5\%$, while S-IDENT (W) keeps this perfect recovery up to $10\%$; in contrast, both Weak-IDENT and Weak-SINDy fail to exactly recover the equation once the noise reaches $1\%$ or higher.  

Overall, S-IDENT gives results comparable to those of state-of-the-art weak-form methods. This is noteworthy since S-IDENT does not rely on an integral formulation; instead, it uses approximate differentiation to estimate the differential features.

\subsection{The choice of the SG polynomial degree}
\label{sec:order_accuracy}

We propose using SURE to automatically select the window-length parameter for SG when approximating derivatives, adapting to the noise in the data. The SG method also involves a choice of the polynomial degree, which governs the order of accuracy of the estimated derivatives.
For the polynomial degree, we fix $d=7$ as the default, and we provide empirical evidence to support this choice.

\begin{figure}
\centering
\setlength{\tabcolsep}{2pt}
\renewcommand{\arraystretch}{1.2}
\begin{tabular}{cccc}
& \textbf{Viscous Burgers} & \textbf{KdV} & \textbf{Allen-Cahn} \\
\adjustbox{valign=c}{\rotatebox{90}{\textbf{TPR}}} &
\adjustbox{valign=c}{\includegraphics[width=0.28\textwidth]{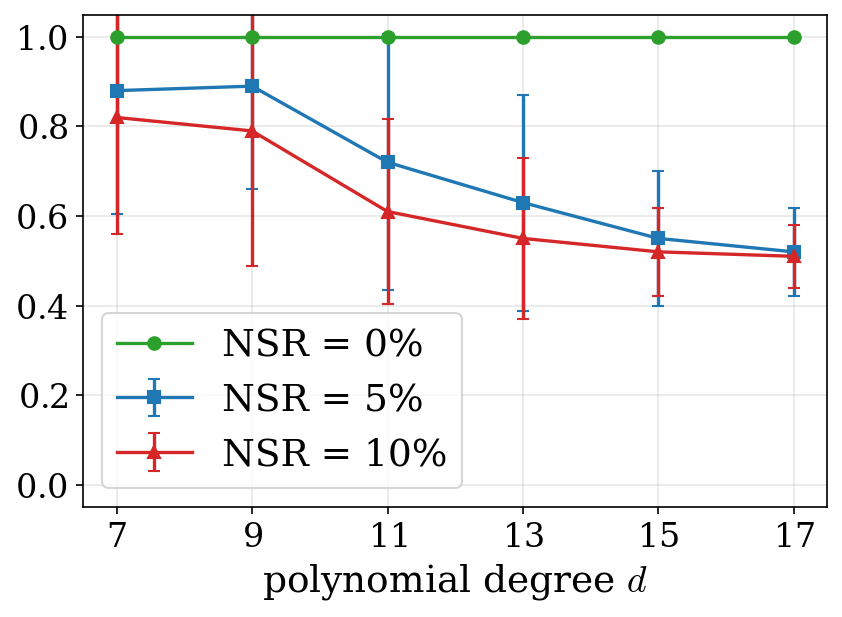}} &
\adjustbox{valign=c}{\includegraphics[width=0.28\textwidth]{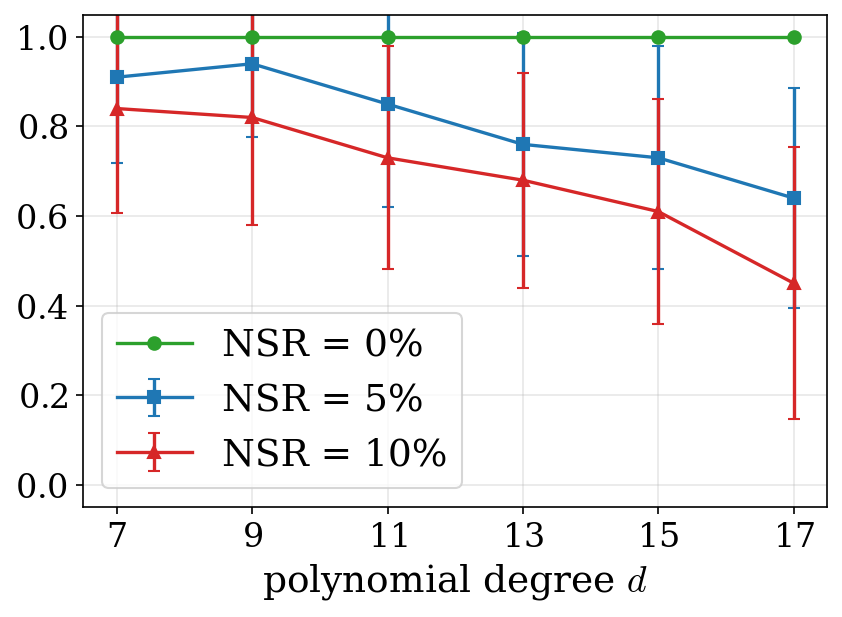}} &
\adjustbox{valign=c}{\includegraphics[width=0.28\textwidth]{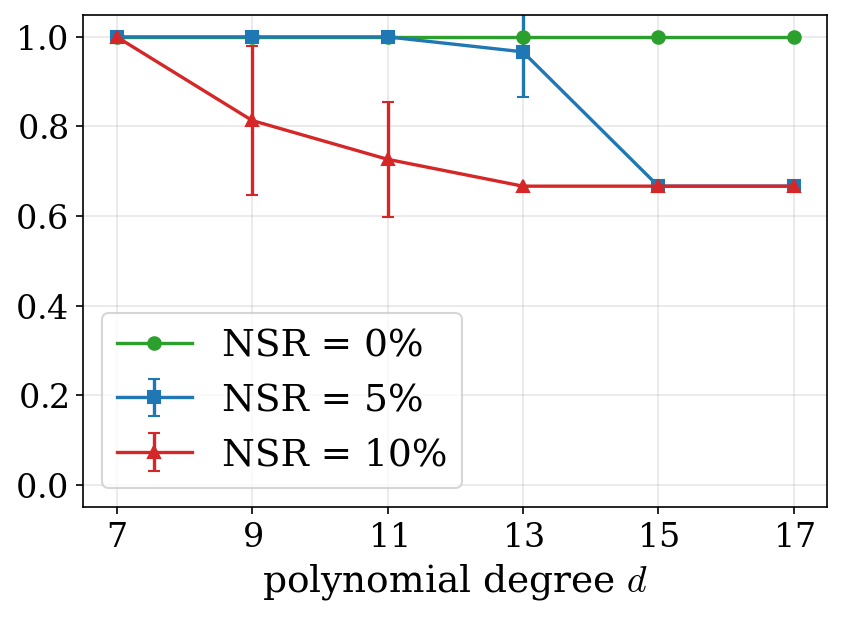}} \\
\adjustbox{valign=c}{\rotatebox{90}{\textbf{PPV}}} &
\adjustbox{valign=c}{\includegraphics[width=0.28\textwidth]{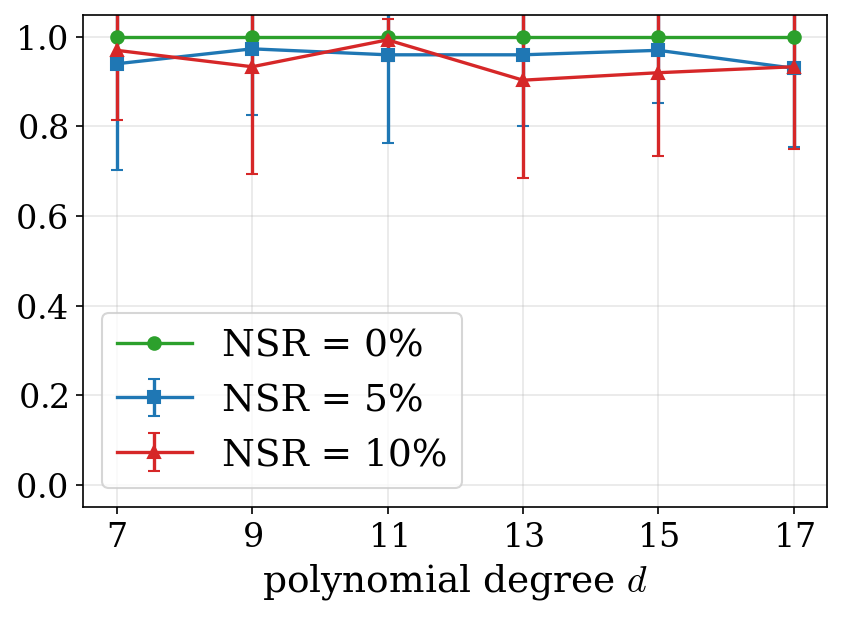}} &
\adjustbox{valign=c}{\includegraphics[width=0.28\textwidth]{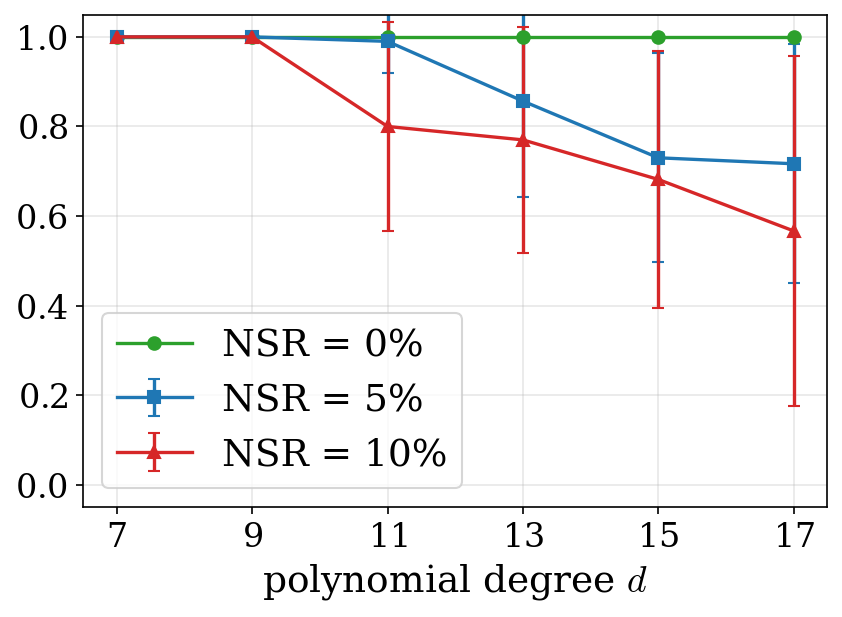}} &
\adjustbox{valign=c}{\includegraphics[width=0.28\textwidth]{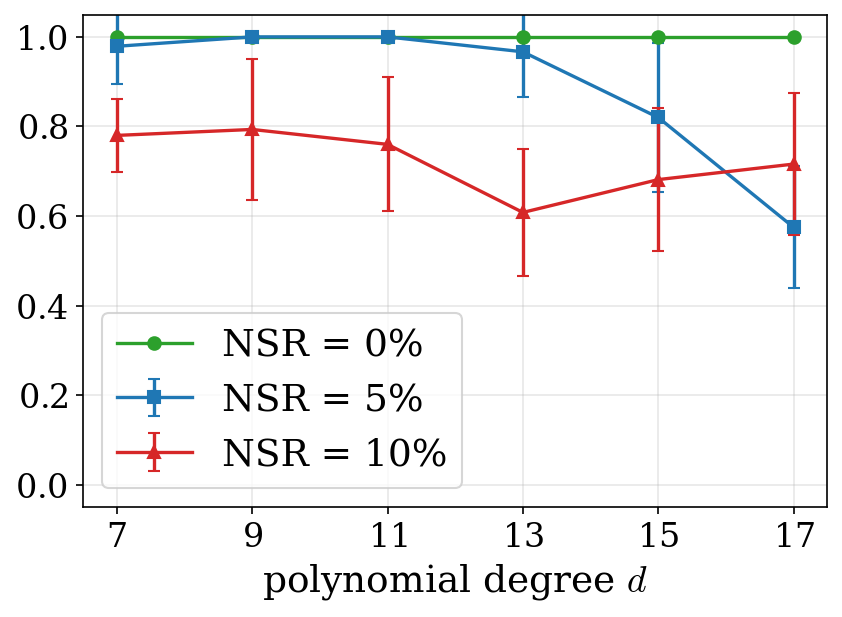}} \\
\adjustbox{valign=c}{\rotatebox{90}{$\boldsymbol{E_{\mathrm{in}}}$}} &
\adjustbox{valign=c}{\includegraphics[width=0.28\textwidth]{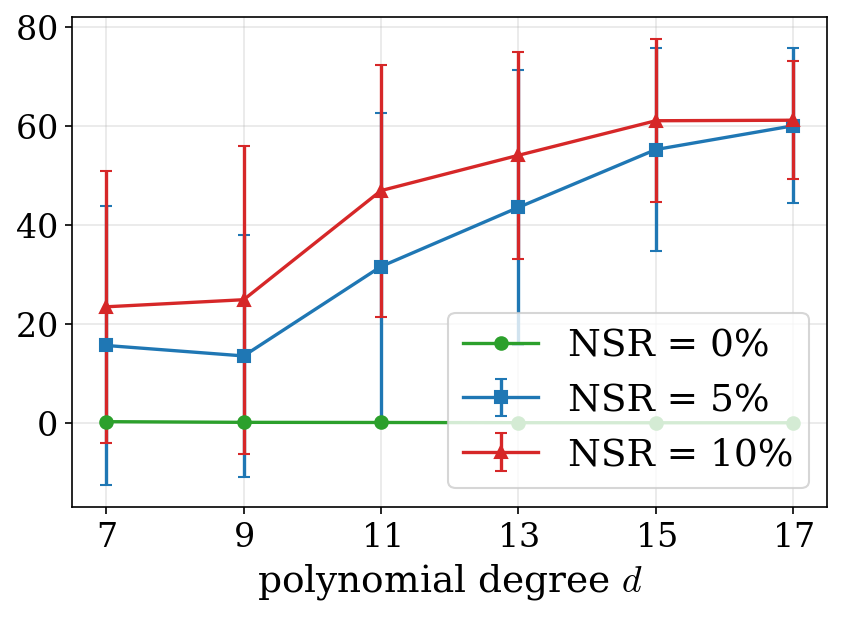}} &
\adjustbox{valign=c}{\includegraphics[width=0.28\textwidth]{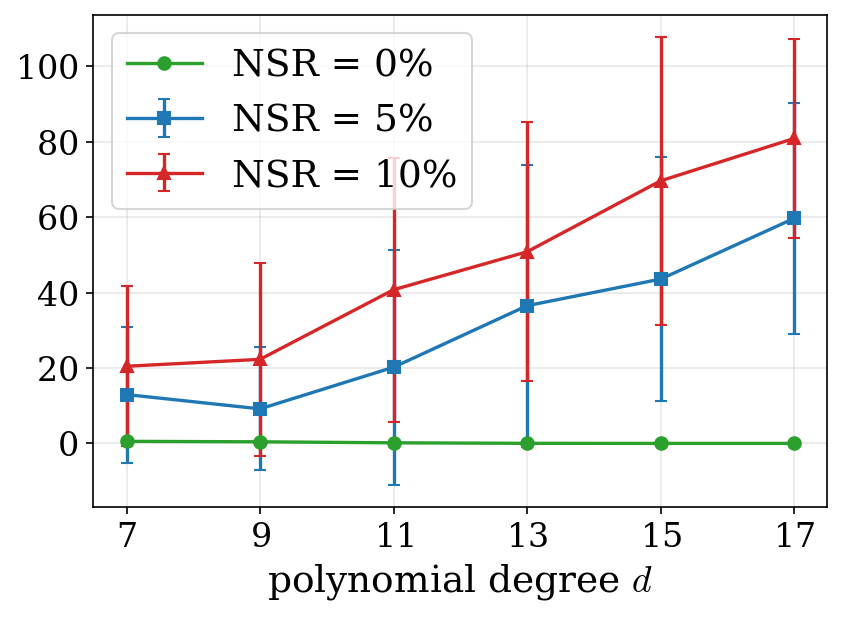}} &
\adjustbox{valign=c}{\includegraphics[width=0.28\textwidth]{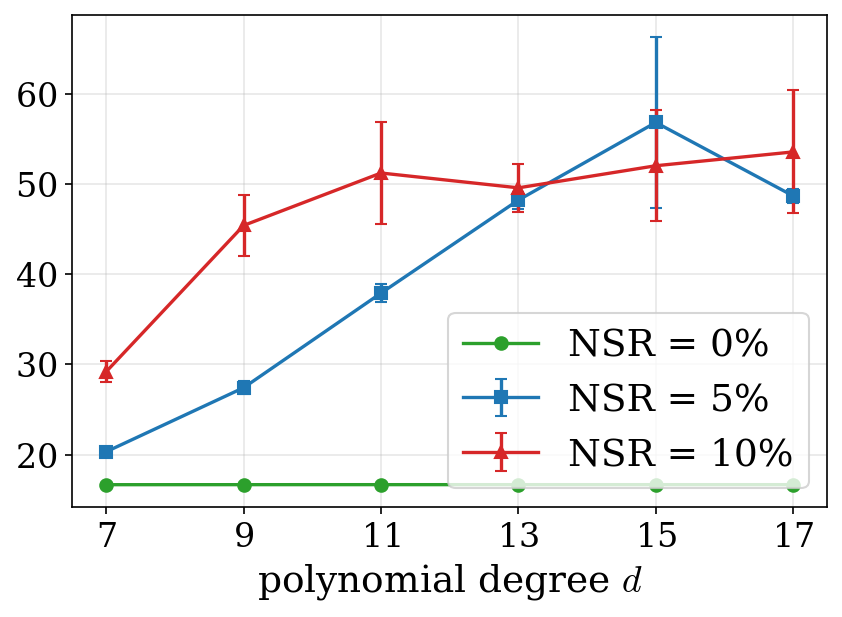}} \\
\adjustbox{valign=c}{\rotatebox{90}{$\boldsymbol{E_{\mathrm{out}}}$}} &
\adjustbox{valign=c}{\includegraphics[width=0.28\textwidth]{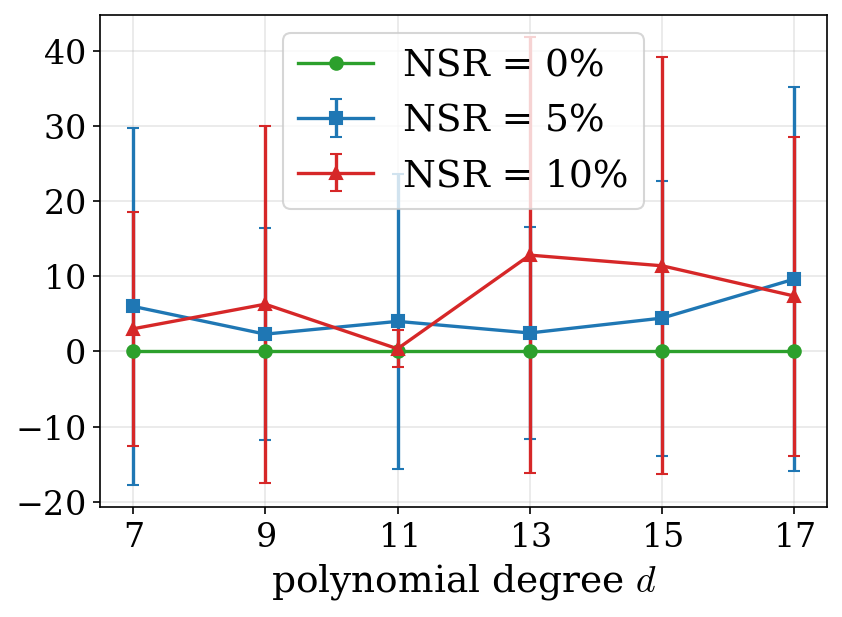}} &
\adjustbox{valign=c}{\includegraphics[width=0.28\textwidth]{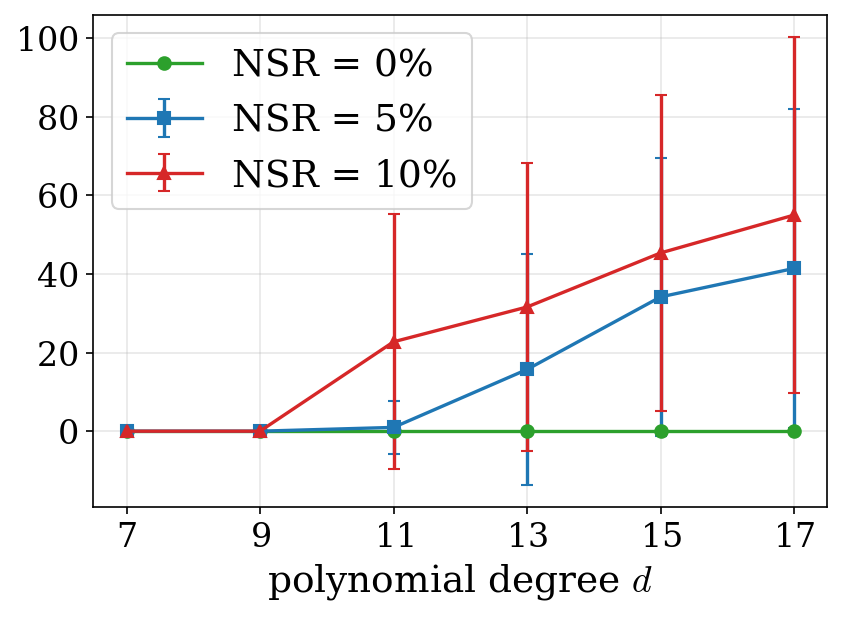}} &
\adjustbox{valign=c}{\includegraphics[width=0.28\textwidth]{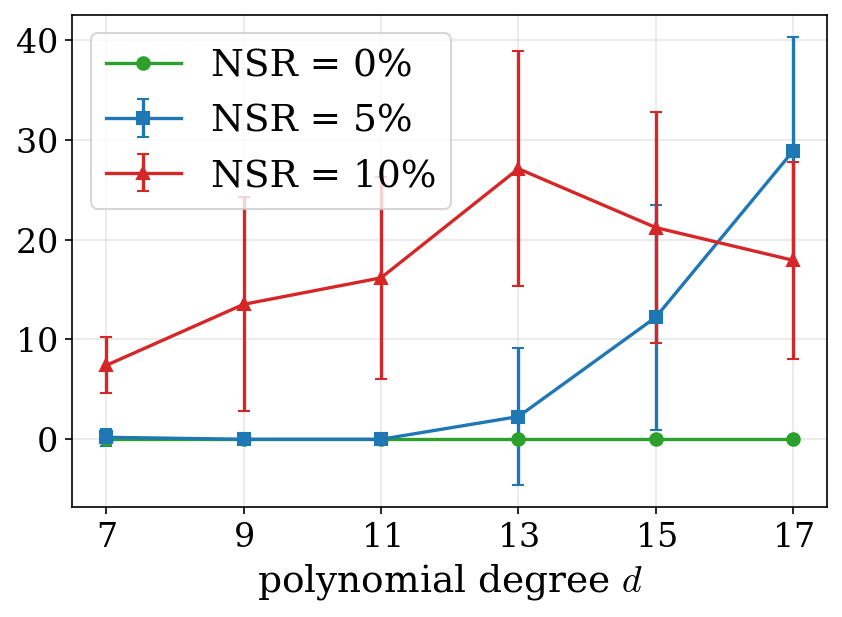}} \\
\adjustbox{valign=c}{\rotatebox{90}{\textbf{Recovery}}} &
\adjustbox{valign=c}{\includegraphics[width=0.28\textwidth]{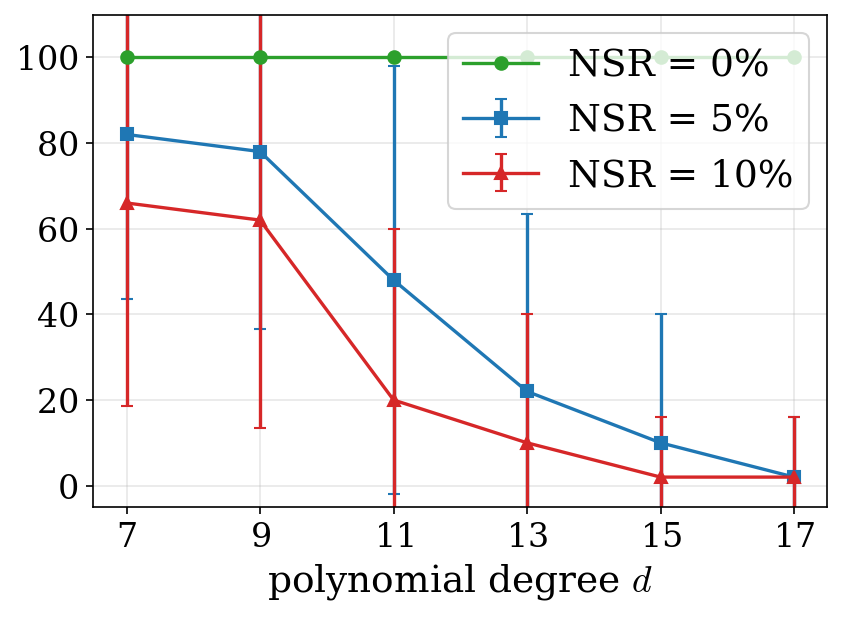}} &
\adjustbox{valign=c}{\includegraphics[width=0.28\textwidth]{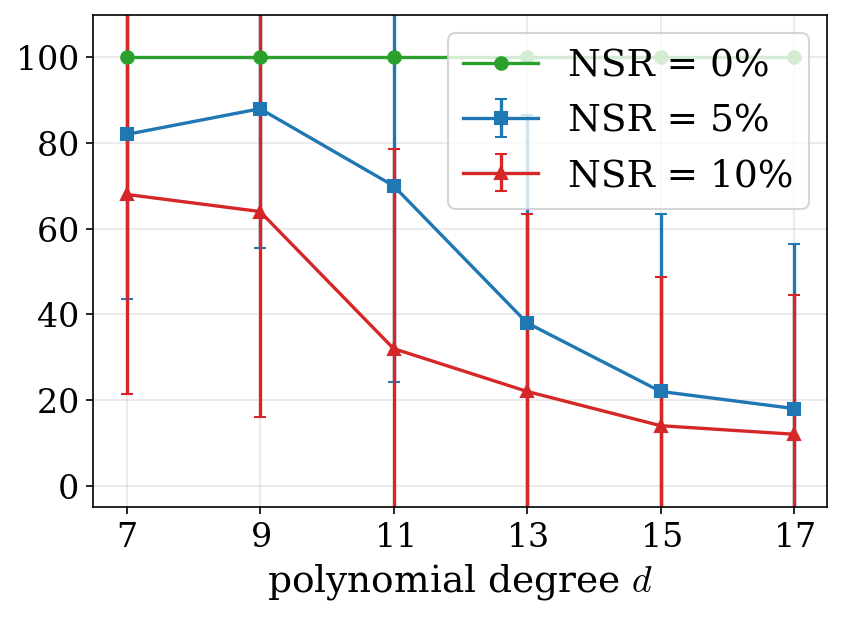}} &
\adjustbox{valign=c}{\includegraphics[width=0.28\textwidth]{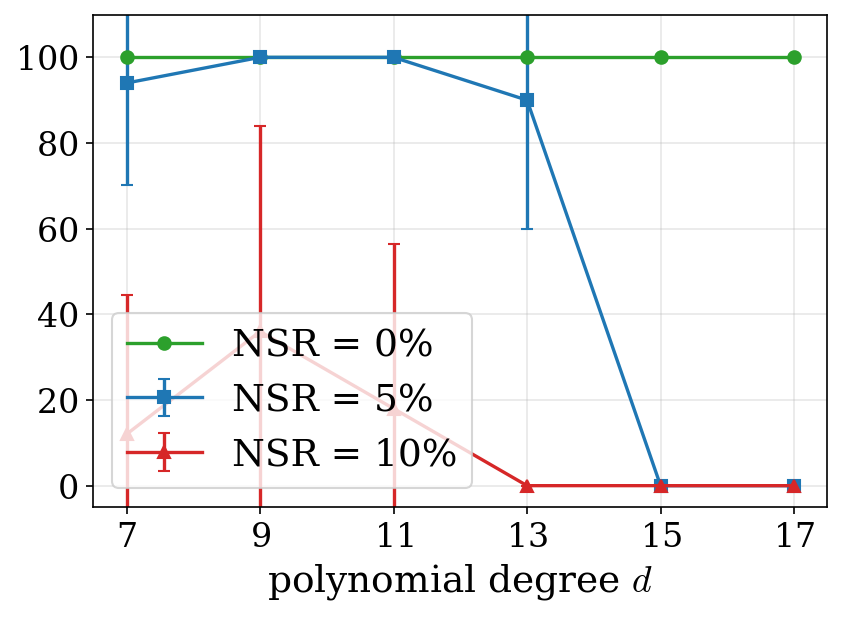}} \\
\end{tabular}
\caption{SG polynomial degree vs S-IDENT results: Identification performance with the Type-S dictionary for the viscous Burgers equation~\eqref{eq:vb}, the KdV equation~\eqref{eq:kdv}, and the Allen--Cahn equation~\eqref{eq:allencahn}. Results for polynomial degrees 7 to 11 seem reasonable for Allen--Cahn, while degree 7 gives the best results for the viscous Burgers and KdV equations. We fix the SG polynomial degree to 7 in this paper.}
\label{fig:poly-order-DP}
\end{figure}

In Figure~\ref{fig:poly-order-DP}, we show S-IDENT results while increasing the polynomial degree used for the approximation. On both clean and noisy data ($5\%$ and $10\%$ NSR), we test the identification with the Type-S $(6,4)$ dictionary for the viscous Burgers~\eqref{eq:vb}, KdV~\eqref{eq:kdv}, and Allen--Cahn~\eqref{eq:allencahn} equations. For each choice of polynomial degree, we run $50$ independent experiments, as before. 
For clean data, the model identification precision and the coefficient reconstruction accuracy remain steady. In contrast, on noisy data, the performance generally deteriorates as higher-degree polynomials are used. Compared to $5\%$ noise, when the data contains higher noise ($10\%$), the coefficient reconstructions are significantly worse for higher degrees. This is due to the amplification of the high-frequency regime as the order of accuracy increases; see the analysis in Subsection~\ref{sec:bias} and Figure~\ref{fig:mag-response}. Notably, the feature identification metrics may be less sensitive to the polynomial degree. For instance, the TPR and PPV for the viscous Burgers and Allen--Cahn equations remain relatively stable as the degree changes; however, such behavior depends on the underlying PDE, as shown by the example of KdV.

In general, we see that the polynomial degree $d=7$ in our SURE-SG method for feature approximation is a valid default value for the experiments in this paper. This numerical evidence also shows the importance of the hyper-parameters in the differentiation methods. Thanks to SURE, our method is easy to use, and choosing the polynomial degree just high enough to cover the highest order of derivative in the dictionary is a robust option. The results obtained using the Type-W dictionaries are similar, and we collect them in Appendix~\ref{sec:w-type-more}, Figure~\ref{fig:poly-order-PD}.

\subsection{Comparison of variants for feature approximation} 
\label{sec:ablation_feature_approx}
\begin{figure}
    \centering
        \centering
 \includegraphics[width=\textwidth]{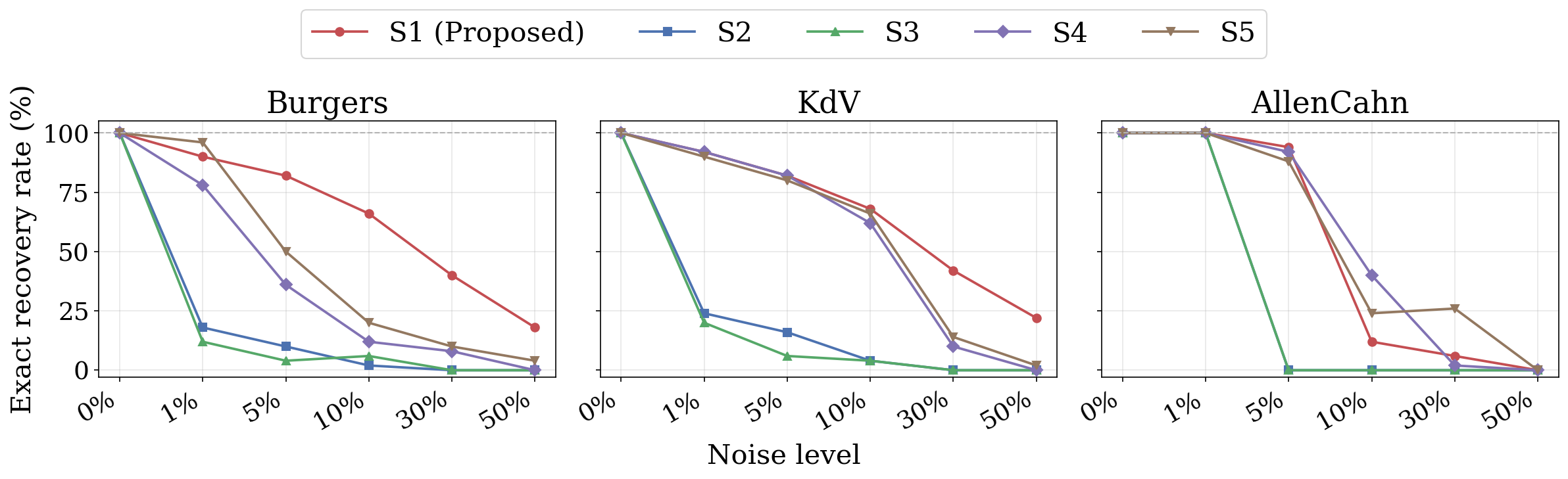}
\caption{Variant study of S1--S5 in Subsection~\ref{sec:ablation_feature_approx}: Exact recovery rate (E.R.)~\eqref{eq:ER} from $50$ independent trials for identifying the viscous Burgers~\eqref{eq:vb}, KdV~\eqref{eq:kdv}, and Allen--Cahn~\eqref{eq:allencahn} equations using the feature approximation methods S1--S5. The proposed S-IDENT shows the best or comparable results.}
\label{fig:compare_diff_methods}
\end{figure}

We compare S-IDENT with several variants in which the feature estimation method SURE-SG is substituted by MaxPol (Section~\ref{sec:maxflat}) in space, MaxPol in both space and time, and SDD (Section~\ref{sec:sg}), where differentiation uses finite differences instead of ENO as in~\cite{he2022robust}. Specifically, we compare the following feature approximation methods:
\begin{itemize}
\item[\textbf{S1.}] \textbf{SURE-SG (proposed):} feature derivatives are estimated by SURE-SG using polynomial degree $7$ in space and $5$ in time, the same as the default.
\item[\textbf{S2.}] \textbf{MaxPol (space):} the spatial feature derivatives are estimated by MaxPol$(B,A)$ with $A=7$ and $B$ chosen so that the stencil length equals the spatial SURE-SG window length; the temporal derivative uses SURE-SG.
\item[\textbf{S3.}] \textbf{MaxPol (space and time):} MaxPol$(B,A)$ replaces SURE-SG on both axes, matched as above, with $A=7$ in space and $A=5$ in time, and $B$ chosen to match the SURE-SG window length in each dimension.
\item[\textbf{S4.}] \textbf{SDD (window-matched):} the $(\mathcal{S}\mathcal{D})^{k}\mathcal{S}$ scheme of~\cite{he2022robust} on both axes for approximating the $k$-th order derivative; here $\mathcal{S}$ is the MLS smoother~\eqref{eq_SDD_mls} with polynomial degree $7$ in space and $5$ in time, and $\mathcal{D}$ a finite difference, both using a stencil length equal to the matched SURE-SG window length.
\item[\textbf{S5.}] \textbf{SDD (accuracy-matched):} the same $(\mathcal{S}\mathcal{D})^{k}\mathcal{S}$ scheme, but with $\mathcal{D}$ a finite difference using $9$ points in space and $7$ in time.
\end{itemize}
The parameters in MaxPol are chosen so that the differentiation attains the same order of accuracy and the resulting window lengths match those of SURE-SG. As for SDD, the window length of the smoother $\cS$ matches the window length of SURE-SG. Since the stencil size of a finite difference $\cD$ is directly related to the order of accuracy, we compare SDD (window-matched), which matches the window size (and thus uses a higher order of accuracy), and SDD (accuracy-matched), which matches the order of accuracy with a smaller stencil. Figure~\ref{fig:compare_diff_methods} shows the comparison of mean E.R.~\eqref{eq:ER} when identifying the viscous Burgers~\eqref{eq:vb}, KdV~\eqref{eq:kdv}, and Allen--Cahn~\eqref{eq:allencahn} equations with the above strategies. We find that using MaxPol (both S2 and S3) deteriorates the identification precision, mostly due to the variability of the estimated features. In the case of Allen--Cahn, we see that the SDD variants (both S4 and S5) yield results comparable to those of S-IDENT using SURE-SG, while in the case of viscous Burgers, the proposed method remains the best across most noise levels. These results justify the feature approximation method of S-IDENT.

Next, we explore and compare different differentiation strategies for estimating higher-order derivatives in space. For SURE-SG in this paper, we fix the polynomial degree to 7 in space and approximate derivatives of different orders from this polynomial, assuming the underlying differential equation does not have any higher-order terms. Alternatively, one can consider the SDD strategy, where first-order SG differentiation is applied repeatedly. Another way is to match the order of approximation for each feature term by using a different polynomial degree for each derivative order. 
Specifically, let ${\cP}_{d,w}^{(s)}$ denote the SG differentiation (Section~\ref{sec:sg}) that approximates the $s$-th order derivative using polynomial degree $d$ and window size $w$. We consider the following strategies to estimate the $s$-th order derivative:
\begin{itemize}
\item \textbf{Direct$(w,d)$ (proposed)}: The $s$-th order derivative of the local fitting polynomial of degree $d>s$ is used to approximate the data derivative, i.e., we directly apply $\cP_{d,w}^{(s)}$, where the window length $w$ is determined by SURE minimization (Section~\ref{sec:sure_sg}) with $\cP_{d,w}^{(0)}$.
\item \textbf{Repeated$(w,d)$}: Apply the SG differentiation of polynomial degree $d$ to the data $s$ times, i.e., $\left(\cP_{d,w}^{1}\right)^{s}$. This is analogous to the SDD~\cite{he2022robust}, and $w$ is determined by SURE minimization with $\cP_{d,w}^{(0)}$, as in the direct approach.
\item \textbf{Adaptive$(w,d)$}: Apply the SG differentiation with local polynomials of degree $s+k$ and estimate the $s$-th order derivative from the fitted polynomial, so that the excess degree $k$, and hence the order of approximation, is the same for every derivative order; i.e., we apply $\cP_{k+s,w}^{(s)}$, where $w$ is determined by SURE minimization with $\cP_{k+s,w}^{(0)}$. Note that the window length $w$ varies with $s$.
\end{itemize}

\begin{figure}
\centering
\begin{tabular}{c@{\vspace{2pt}}c@{\vspace{2pt}}c}
(a)&(b)&(c)\\
\includegraphics[width=0.33\textwidth,height=1.7in]{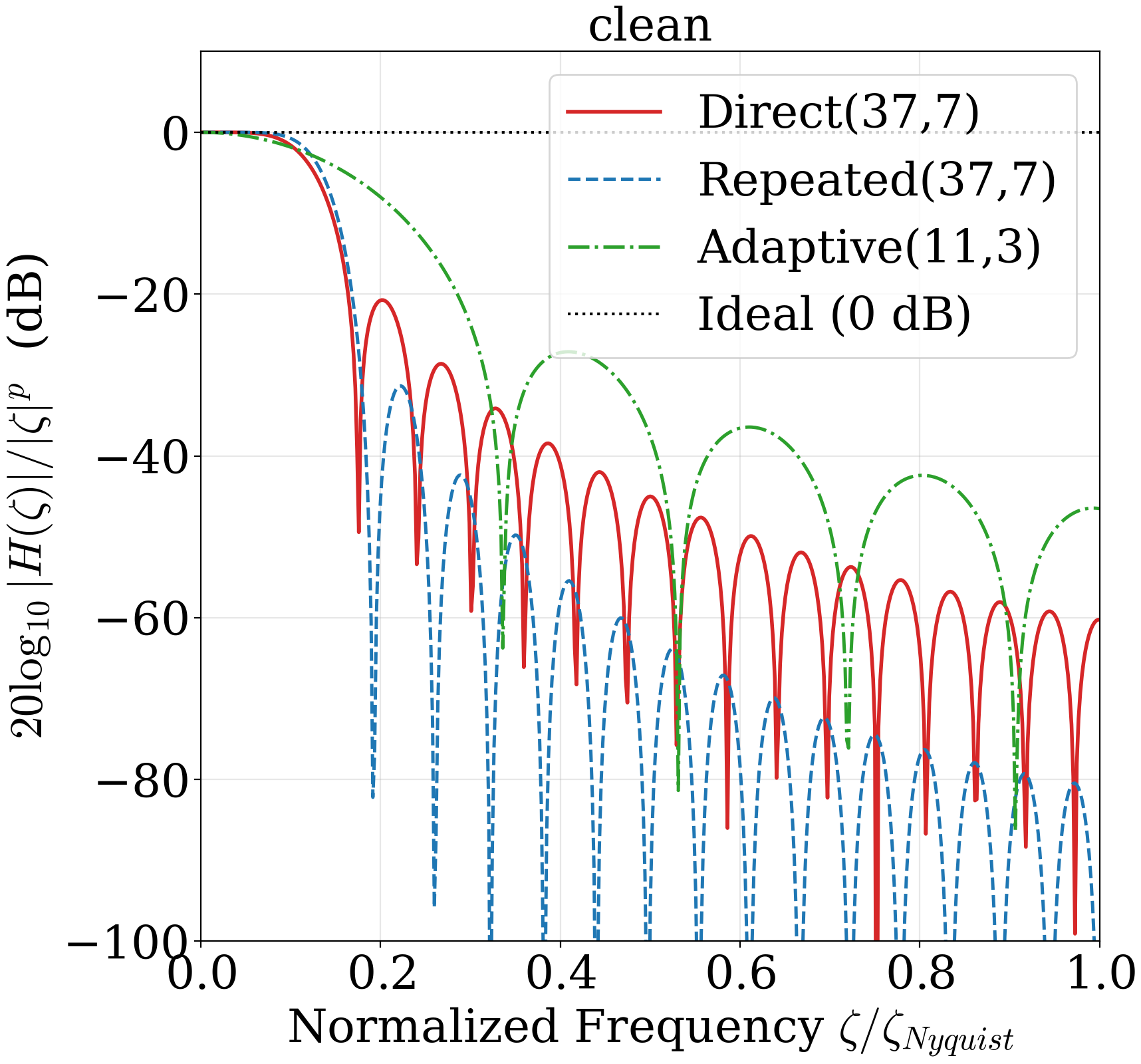}&
\includegraphics[width=0.33\textwidth,height=1.7in]{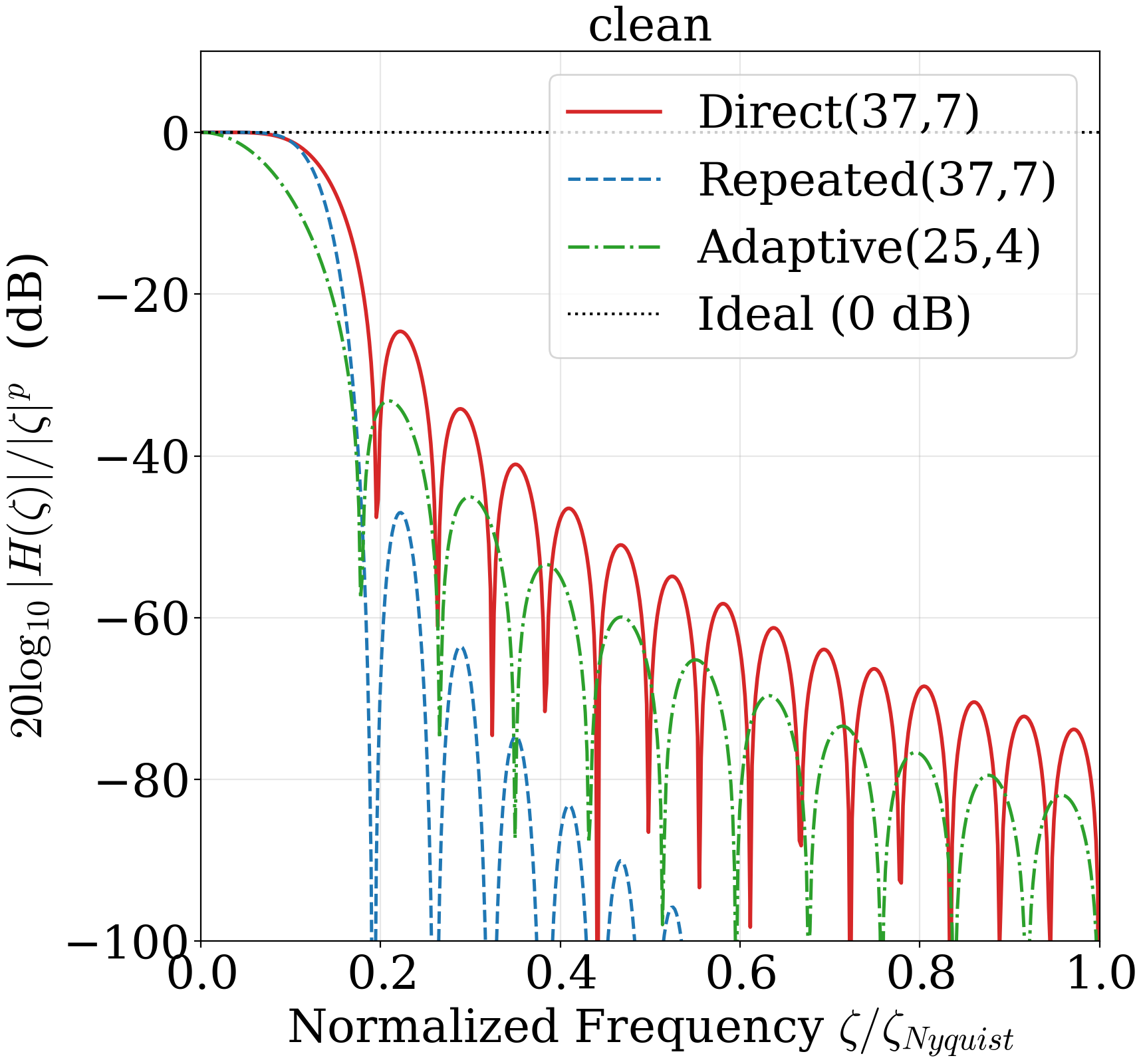}&
\includegraphics[width=0.33\textwidth,height=1.7in]{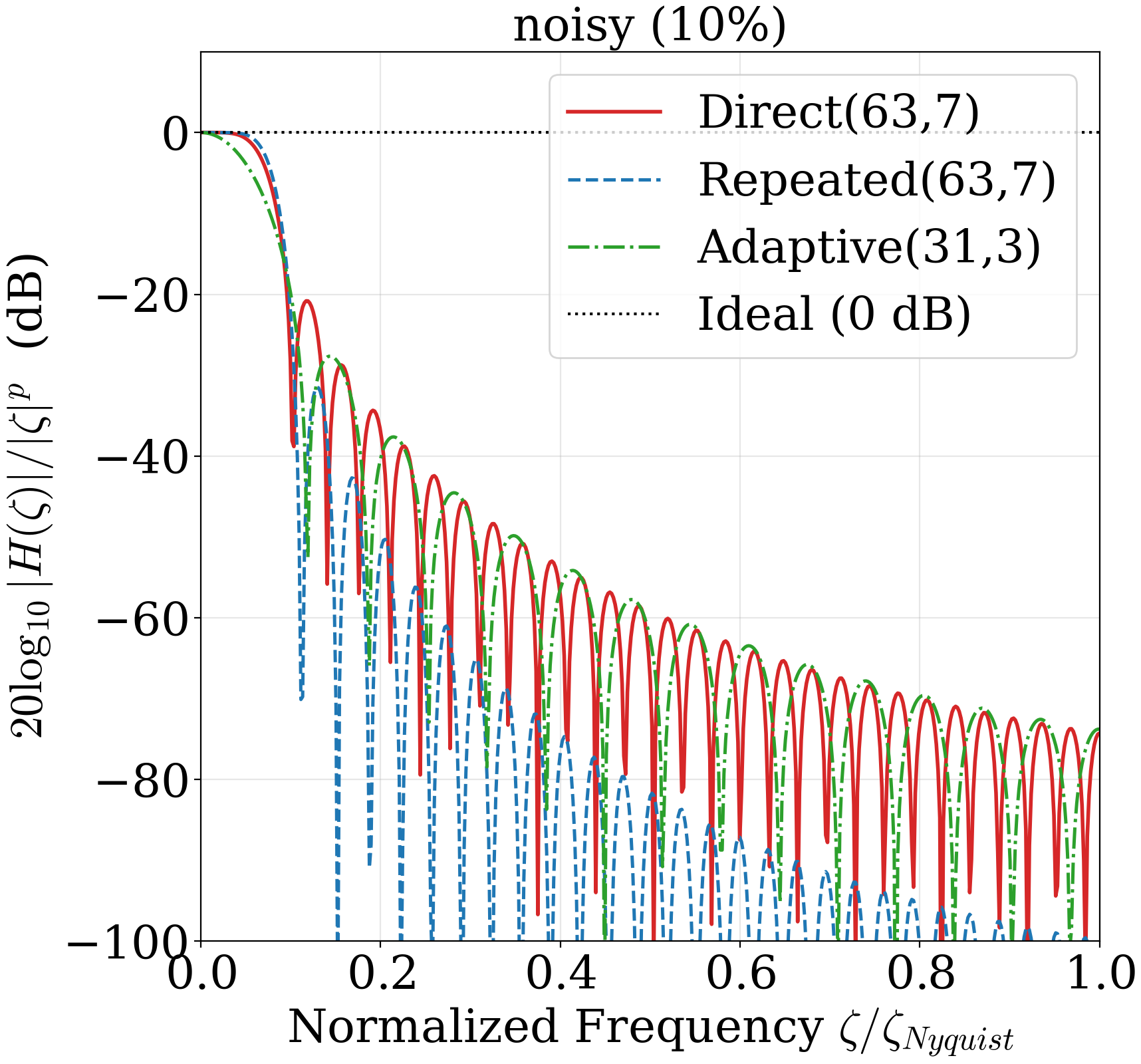}
\end{tabular}
\caption{ Variant study of Direct (Proposed), Repeated, and Adaptive in Subsection \ref{sec:ablation_feature_approx}:
Relative magnitude response~\eqref{eq:magnitudeR} is shows to approximate (a) $2$-nd order derivative and (b) $3$-rd order derivative from clean data, and (c) $2$-nd order derivative from data with $10\%$ noise.} \label{fig:higher-order}
\end{figure}

We compare these strategies by choosing $d=7$ and $k=1$ so that, for the features in the Type-S $(6,4)$ dictionary, they have matching orders of accuracy when estimating the highest-order derivative feature. 
Figure~\ref{fig:higher-order} shows the relative magnitude responses~\eqref{eq:magnitudeR} of these strategies when approximating (a) the second-order and (b) the third-order derivatives from clean data, and (c) the second-order derivative from data with $10\%$ noise. In general, we observe that all strategies exhibit low-pass filtering behavior while maintaining a tight approximation near zero frequency, indicating that condition~\eqref{eq_frequency_condition} is satisfied. We notice that, compared to the Direct strategy (used in S-IDENT), the relative magnitude response~\eqref{eq:magnitudeR} of the Repeated strategy decays faster, more strongly suppressing the higher-frequency content, which can lead to more bias. Compared to the Direct strategy, the Adaptive strategy has a lower but consistent order of accuracy for each derivative feature.

Table~\ref{tab:diff_compare_dp} shows the performance of S-IDENT using the above three feature approximation strategies. We test with the viscous Burgers~\eqref{eq:vb}, KdV~\eqref{eq:kdv}, and Allen--Cahn~\eqref{eq:allencahn} equations. We observe that the proposed S-IDENT with the Direct strategy outperforms the other strategies in most cases in terms of the exact recovery rate. The Repeated strategy maintains a high PPV, i.e., avoiding the wrong features, and the Adaptive strategy does not show decisive benefits over the others. The results and conclusions are similar for the identification using the Type-W $(6,4)$ dictionary, with overall better accuracy, and we show them in Appendix~\ref{sec:w-type-more}, Table~\ref{tab:diff_compare_pd}.

\begin{table}
\centering
\small
\begin{tabular}{llccccc}
\toprule
NSR & Method & TPR & PPV & $E_{\text{in}}$ & $E_{\text{out}}$ & E.R. (\%) \\
\midrule
\multicolumn{7}{c}{(a) \textit{Viscous Burgers} \eqref{eq:vb}} \\
\midrule
 & Direct & \textbf{0.950 $\pm$ 0.180} & 0.973 $\pm$ 0.147 & \textbf{6.090 $\pm$ 17.867} & 2.341 $\pm$ 14.154 & \textbf{90.0 $\pm$ 30.0} \\
1\% & Repeated & 0.860 $\pm$ 0.246 & \textbf{0.980 $\pm$ 0.140} & 15.686 $\pm$ 26.836 & \textbf{2.000 $\pm$ 14.000} & 74.0 $\pm$ 43.9 \\
 & Adaptive & 0.900 $\pm$ 0.224 & \textbf{0.980 $\pm$ 0.140} & 13.637 $\pm$ 22.613 & \textbf{2.000 $\pm$ 14.000} & 82.0 $\pm$ 38.4 \\
\cmidrule(lr){2-7}
 & Direct & 0.820 $\pm$ 0.260 & \textbf{0.970 $\pm$ 0.155} & 23.445 $\pm$ 27.476 & \textbf{3.006 $\pm$ 15.542} & \textbf{66.0 $\pm$ 47.4} \\
10\% & Repeated & 0.770 $\pm$ 0.287 & 0.940 $\pm$ 0.208 & 27.045 $\pm$ 30.482 & 4.661 $\pm$ 19.647 & 52.0 $\pm$ 50.0 \\
 & Adaptive & \textbf{0.840 $\pm$ 0.273} & 0.917 $\pm$ 0.222 & \textbf{21.017 $\pm$ 25.736} & 5.860 $\pm$ 20.077 & 62.0 $\pm$ 48.5 \\
\midrule
\multicolumn{7}{c}{(b) \textit{KdV} \eqref{eq:kdv}} \\
\midrule
 & Direct & 0.960 $\pm$ 0.136 & \textbf{1.000 $\pm$ 0.000} & 5.717 $\pm$ 13.146 & \textbf{0.000 $\pm$ 0.000} & \textbf{92.0 $\pm$ 27.1} \\
1\% & Repeated & \textbf{0.970 $\pm$ 0.119} & 0.993 $\pm$ 0.047 & \textbf{5.010 $\pm$ 16.732} & 0.328 $\pm$ 2.295 & \textbf{92.0 $\pm$ 27.1} \\
 & Adaptive & 0.570 $\pm$ 0.173 & 0.960 $\pm$ 0.136 & 55.646 $\pm$ 16.144 & 1.865 $\pm$ 6.818 & 14.0 $\pm$ 34.7 \\
\cmidrule(lr){2-7}
 & Direct & 0.840 $\pm$ 0.233 & \textbf{1.000 $\pm$ 0.000} & 20.460 $\pm$ 21.243 & \textbf{0.000 $\pm$ 0.000} & 68.0 $\pm$ 46.6 \\
10\% & Repeated & \textbf{0.850 $\pm$ 0.229} & \textbf{1.000 $\pm$ 0.000} & \textbf{19.880 $\pm$ 25.320} & \textbf{0.000 $\pm$ 0.000} & \textbf{70.0 $\pm$ 45.8} \\
 & Adaptive & 0.500 $\pm$ 0.283 & 0.750 $\pm$ 0.377 & 67.211 $\pm$ 22.161 & 22.853 $\pm$ 39.330 & 16.0 $\pm$ 36.7 \\
\midrule
\multicolumn{7}{c}{(c) \textit{Allen-Cahn} \eqref{eq:allencahn}} \\
\midrule
 & Direct & \textbf{1.000 $\pm$ 0.000} & \textbf{1.000 $\pm$ 0.000} & 16.801 $\pm$ 0.044 & \textbf{0.000 $\pm$ 0.000} & \textbf{100.0 $\pm$ 0.0} \\
1\% & Repeated & \textbf{1.000 $\pm$ 0.000} & \textbf{1.000 $\pm$ 0.000} & 16.867 $\pm$ 0.071 & \textbf{0.000 $\pm$ 0.000} & \textbf{100.0 $\pm$ 0.0} \\
 & Adaptive & \textbf{1.000 $\pm$ 0.000} & \textbf{1.000 $\pm$ 0.000} & \textbf{16.231 $\pm$ 0.094} & \textbf{0.000 $\pm$ 0.000} & \textbf{100.0 $\pm$ 0.0} \\
\cmidrule(lr){2-7}
 & Direct & \textbf{1.000 $\pm$ 0.000} & 0.780 $\pm$ 0.081 & 29.210 $\pm$ 1.137 & 7.410 $\pm$ 2.818 & 12.0 $\pm$ 32.5 \\
10\% & Repeated & \textbf{1.000 $\pm$ 0.000} & 0.805 $\pm$ 0.104 & 30.760 $\pm$ 1.853 & 7.457 $\pm$ 4.109 & 22.0 $\pm$ 41.4 \\
 & Adaptive & \textbf{1.000 $\pm$ 0.000} & \textbf{0.842 $\pm$ 0.125} & \textbf{18.208 $\pm$ 2.534} & \textbf{3.628 $\pm$ 2.986} & \textbf{38.0 $\pm$ 48.5} \\\bottomrule
\end{tabular}
\caption{Variant study of Direct (proposed), Repeated, and Adaptive for S-IDENT (Type-S dictionary) in Subsection~\ref{sec:ablation_feature_approx}: 
For each (PDE, NSR) and each metric, the best of the three SURE-SG strategies is shown in \textbf{bold}. Each entry is the mean $\pm$ standard deviation over $50$ independent trials. The Direct strategy shows the best overall performance.}
\label{tab:diff_compare_dp}
\end{table}

\section{Conclusion}\label{sec:conclude}

We investigate the problem of PDE identification from noisy data. We find that by carefully designing differential estimators, we can achieve robustness comparable to that of weak-form methods when identifying PDEs with general nonlinear features. 
We develop and propose S-IDENT, which leverages Savitzky--Golay (SG) differentiation together with an adaptive approach based on Stein's Unbiased Risk Estimate (SURE) for selecting the window length. We incorporate an effective trimming technique and a model selection criterion based on Residual in Reduction (RR) for more robust identification of differential equations. The proposed method can not only identify more general nonlinear PDEs but also achieve identification performance comparable to that of the state-of-the-art methods, with or without restrictions on the dictionary. 
While S-IDENT is effective with large Type-S dictionaries, we highlight the intrinsic difficulty of identifying PDEs with more general dictionaries by numerically showing that the Type-S feature matrices are generally more ill-conditioned than Type-W ones, especially when the observed trajectory is smooth. See Appendix~\ref{sec:S-correlation}. These observations also point to the gap between the numerical effectiveness of S-IDENT and the theoretical understanding of identifiability, which we leave to future work.

\appendix

\bibliographystyle{plain}
\bibliography{sample}

\appendix

\section{S-IDENT Algorithm}
\label{Asec:Algorithm}

We present the pseudo-code of our proposed S-IDENT in Algorithm~\ref{alg:proposed}. Here the sub-procedure $\textsc{SubspacePursuit}(\bar\bF,\bar{\bm b},\,k)$ refers to the SP algorithm~\cite{dai2009subspace} applied to the normalized system $(\bar\bF,\bar{\bm b})$ using sparsity level $k$; $\textsc{Trimming}(\bar\bF,\bar{\bm b}, S_k,\tau)$ is the trimming technique~\eqref{eq:trimming} described in Subsection~\ref{sec:identification} with threshold $\tau$; and $\textsc{RRScore}\big(\bar\bF,\bar{\bm b}, S_k;\,\rho,L\big)$ is the Reduction in Residual criterion described in Subsection~\ref{sec:identification} using hyper-parameters $\rho$ and $L$; see~\eqref{eq:cand_score} and~\eqref{eq_thresh}. 
\begin{algorithm}[t]
\DontPrintSemicolon
\SetKwInOut{Input}{Input}
\SetKwInOut{Output}{Output}
\SetKwInOut{Param}{Params}
 
\Input{Noisy samples $\{\tilde u_i\}_{i=1}^{N}$ on a space--time grid;
       spatial step $\Delta x$, temporal step $\Delta t$.}
\Param{Candidate features $\{f_1,\dots,f_M\}$;
       highest sparsity $k_{\max}\leq M$; trimming threshold $\tau\in[0,1]$;
       RR-validation parameters $(\rho, L)$.}
\Output{Identified PDE $\;u_t = \sum_{\ell\in S^\star} c_\ell\,f_\ell(u)$.}
\BlankLine
 
\tcp{1. Denoise and differentiate (Savitzky--Golay, SURE-selected window)}
$\hat u \leftarrow \mathcal S_x[\tilde u]$
   \tcp*{spatial smoothing, degree $d_x$, window by SURE}
$\mathbf{b} \leftarrow \mathcal D_t\,\mathcal S_t[\hat u]$
   \tcp*{response $u_t$, temporal degree $d_t$}
\BlankLine
 
\tcp{2. Build the feature matrix}
$\bF \leftarrow \big[\,\widehat{f}_1\;\cdots\;\widehat{f}_M\,\big]\in\mathbb{R}^{N\times M}$ \tcp*{spatial differentiation by SURE-SG}
\BlankLine
 
\tcp{3. Column normalization}
$\bm n_\ell \leftarrow \lVert \bF_{:,\ell}\rVert_2$,\quad
$\bar\bF_{:,\ell} \leftarrow \bF_{:,\ell}/\bm n_\ell$,\quad
$\bar{\bm b} \leftarrow \bm b/\lVert \bm b\rVert_2$\;
\BlankLine
 
\tcp{4. Candidate generation by Subspace Pursuit}
\For{$k \leftarrow 1$ \KwTo $k_{\max}$}{
   $ S_k \leftarrow \textsc{SubspacePursuit}(\bar\bF,\bar{\bm b},\,k)$
       \tcp*{support of size $k$}
   $\bm c_k \leftarrow \arg\min_{\bm c}\;
       \lVert \bar{\bm b} - \bar\bF_{:, S_k}\,\bm c\rVert_2$
       \tcp*{least squares on the support}
   $( S_k,\bm c_k) \leftarrow
       \textsc{Trimming}(\bar\bF,\bar{\bm b}, S_k,\tau)$
       \tcp*{prune terms below relative tolerance $\tau$}
}
\BlankLine
 
\tcp{5. Candidate selection}
$k^\star \leftarrow \arg\min_{k}\;
   \textsc{RRScore}\big(\bar\bF,\bar{\bm b}, S_k;\,\rho,L\big)$\;
$ S^\star \leftarrow  S_{k^\star}$\;
\BlankLine
 
\tcp{6. Refit on the unnormalized system}
$\bm c^\star \leftarrow \arg\min_{\bm c}\;
   \lVert \bm b - \bF_{:, S^\star}\,\bm c\rVert_2$\;
\BlankLine
 
\Return $\;u_t = \sum_{\ell\in S^\star} c^\star_\ell\,f_\ell(u)$\;
\caption{Proposed identification method: S-IDENT}
\label{alg:proposed}
\end{algorithm}

\section{Intrinsic Difficulties with Type-S Dictionaries}\label{sec:S-correlation}

\begin{figure}
\centering
\begin{tabular}{cc}
(a) & (b) \\
\includegraphics[width=0.475\textwidth]{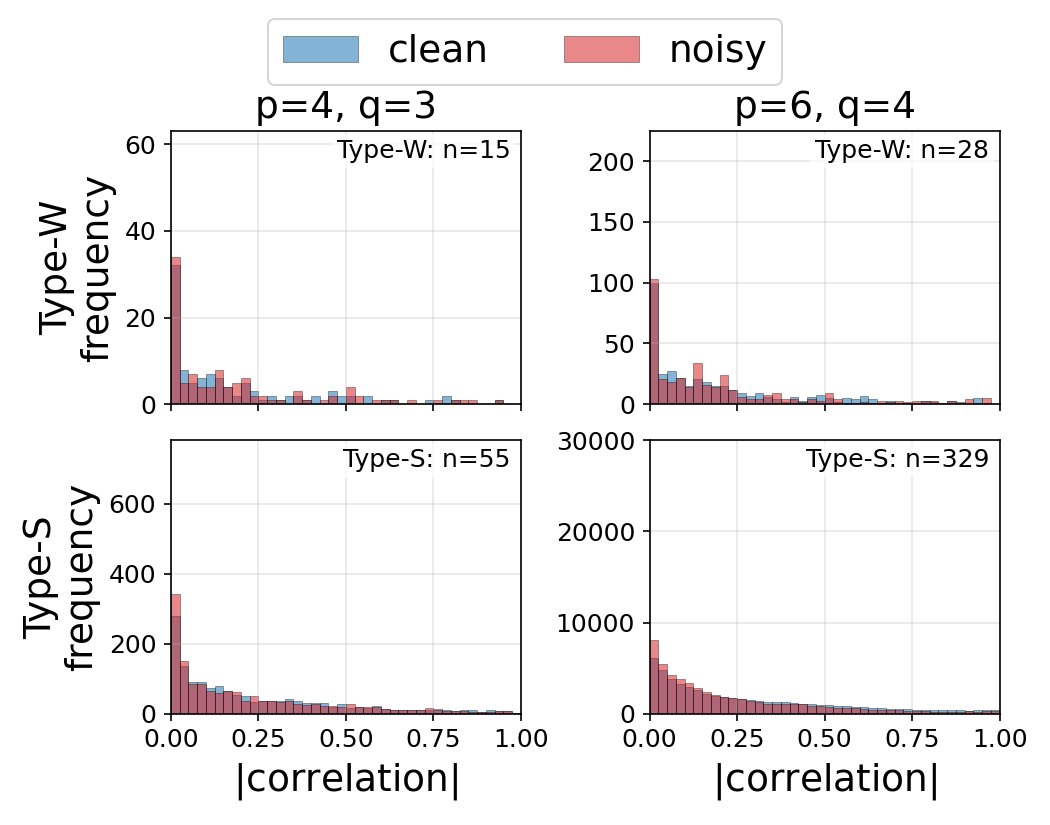} & \includegraphics[width=0.475\textwidth]{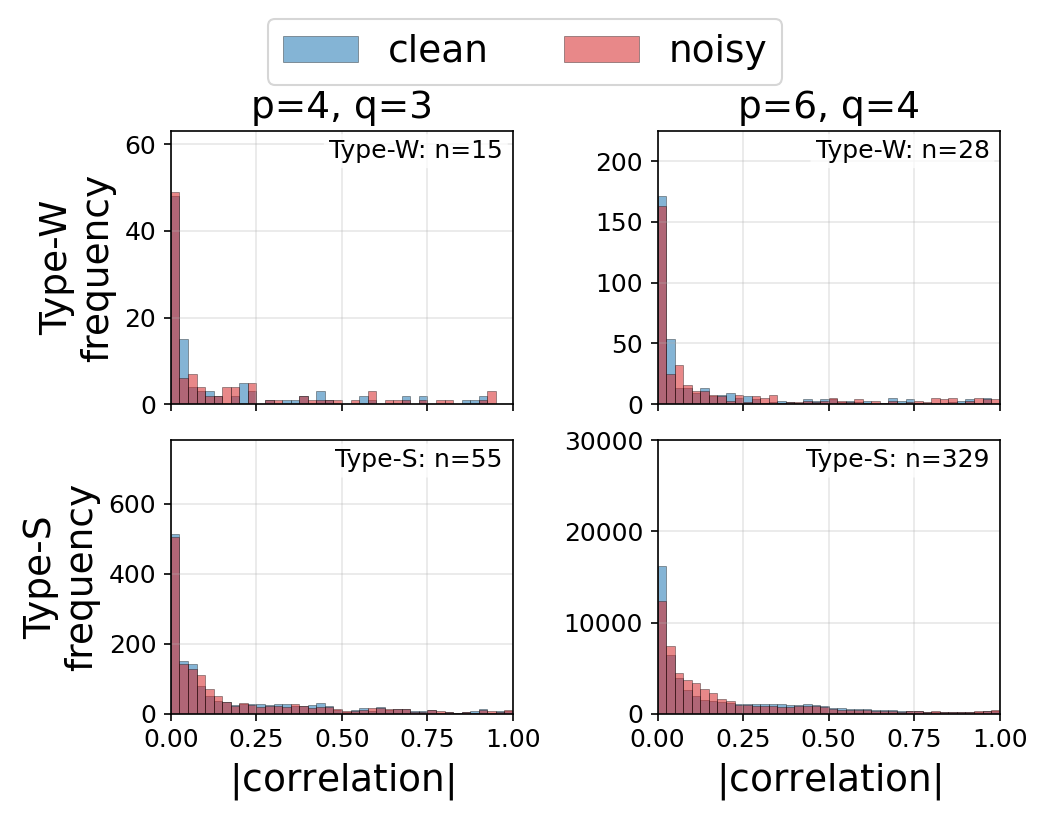}\\
\end{tabular}
\caption{Correlation among feature pairs in the Type-W and Type-S dictionaries induced by (a) viscous Burgers~\eqref{eq:vb} and (b) KdV~\eqref{eq:kdv}: 
Histograms of absolute correlations among feature pairs for both the Type-W and Type-S dictionaries with $(p,q)=(4,3)$ and $(6,4)$, on clean data (blue bars) and noisy data (red bars) with $\mathrm{NSR} = 20\%$.}\label{fig:corr_hist}
\end{figure}

There are intrinsic difficulties in using Type-S dictionaries compared to Type-W ones. First, as shown in Table~\ref{tab_dictionary_sizes} and in Subsection~\ref{sec:comparison_typeW}, there is a difference in the dictionary size needed to cover the same degree of differentiation.  

Second, the features are more correlated in the Type-S dictionaries. We investigate the pairwise correlation among features in the Type-W and Type-S dictionaries. For the Type-W dictionary, we also approximate the features with SURE-SG, as used in S-IDENT, for a fair comparison. Figure~\ref{fig:corr_hist} presents two examples showing the histograms of the absolute pairwise correlation among features in the Type-W and Type-S dictionaries. Panels (a) and (b) correspond to the viscous Burgers~\eqref{eq:vb} and KdV~\eqref{eq:kdv} equations, respectively; within each panel, the first row gives the histograms of correlations for Type-W and the second row those for Type-S. 
Omitting the constant feature, we consider the differentiation order $p$ and multiplication order $q$ with $(p,q)=(4,3)$, giving $15$ terms for Type-W and $55$ for Type-S, and $(p,q)=(6,4)$, giving $28$ terms for Type-W and $329$ for Type-S.
The majority of feature pairs are uncorrelated, yet compared to the Type-W dictionary, relatively more feature pairs in the Type-S dictionaries exhibit higher correlation when the dictionary becomes larger. Some phenomena are equation-dependent. For example, we find that dictionaries of the same type and size induced from the viscous Burgers equation exhibit more correlated pairs than those induced from KdV.
We also observe that the noise (NSR $=20\%$) does not strongly affect the correlation distribution of the Type-W feature pairs. In contrast, the noise generally gives more Type-S feature pairs with larger correlations, as shown by the slightly higher tails. Meanwhile, for smoother trajectories from the viscous Burgers equation, noise can produce more uncorrelated pairs, as shown by the higher first red bar compared to the blue one, whereas for KdV, noise tends to cause more correlated ones, as shown by the lower first red bar compared to the blue one.

\begin{figure}
\centering
\begin{tabular}{ccc}
(a)&(b)&(c)\\
\includegraphics[width=0.3\textwidth]{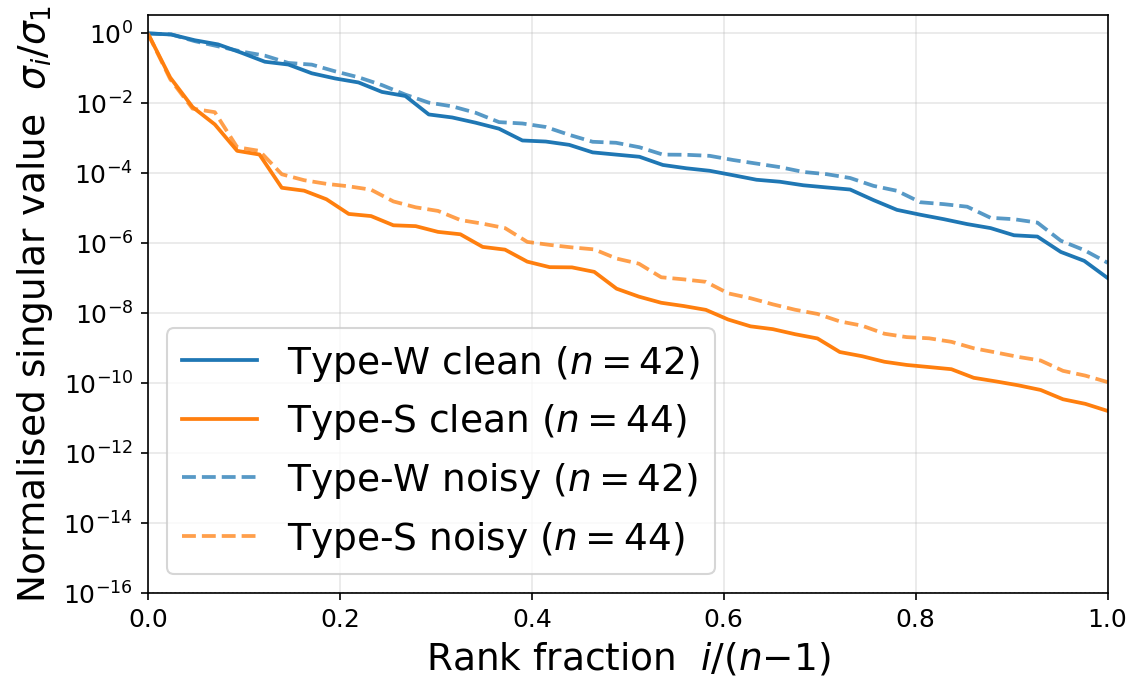}&
\includegraphics[width=0.3\textwidth]{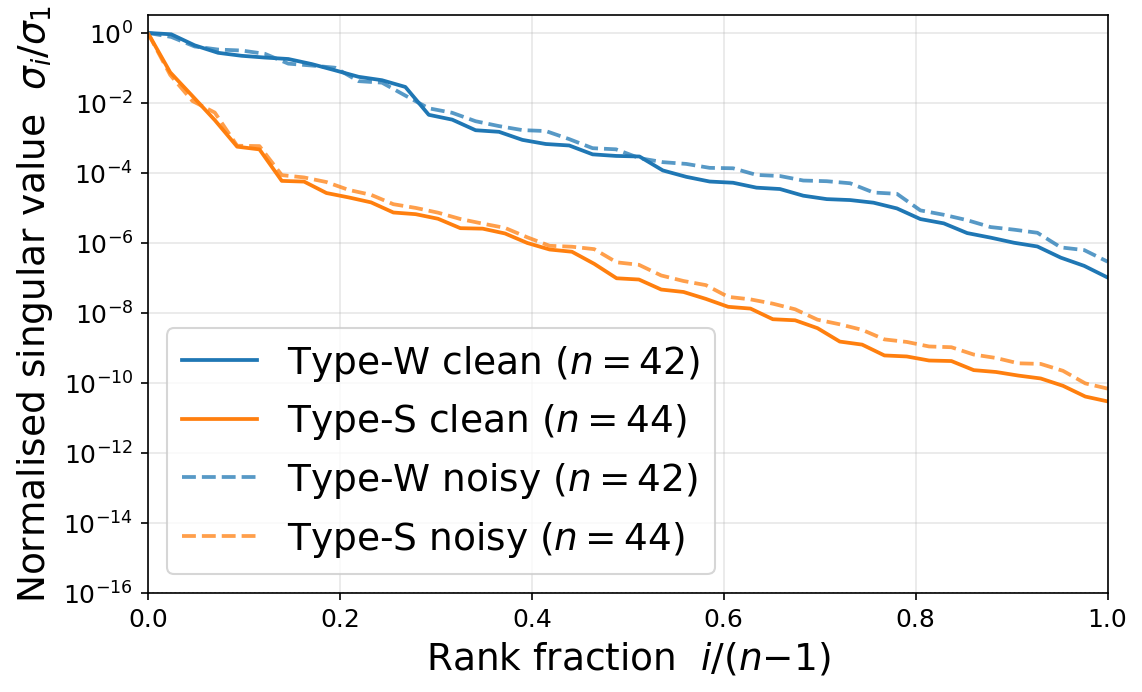}&
\includegraphics[width=0.3\textwidth]{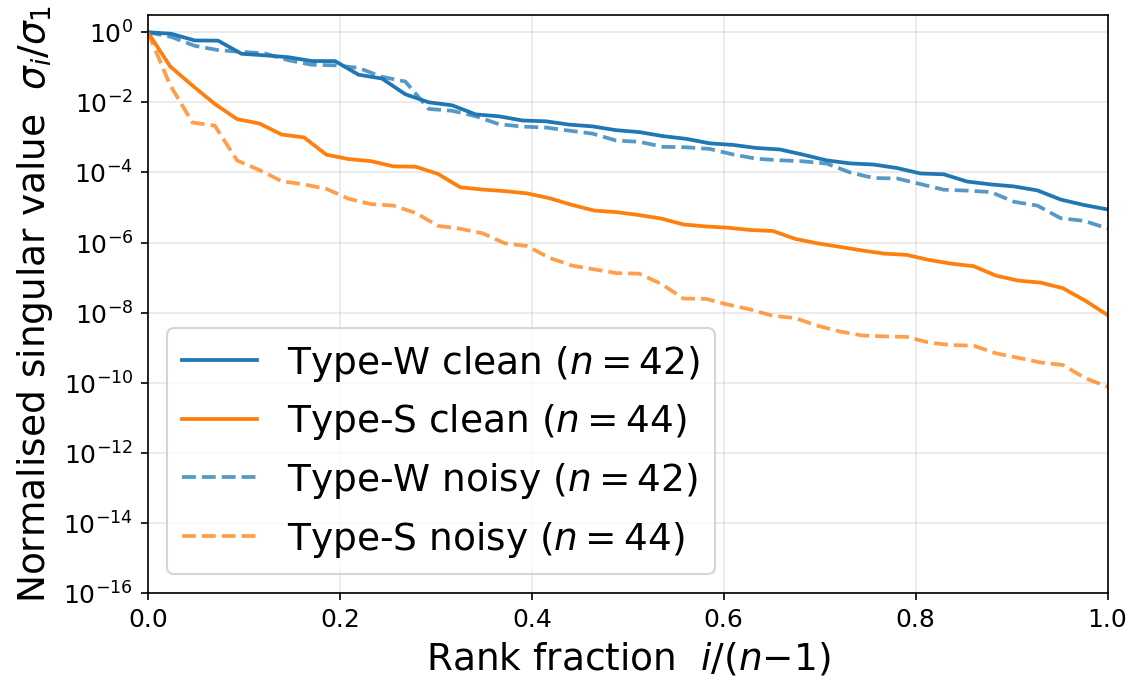}
\end{tabular}
\caption{Singular values of the feature matrices corresponding to the Type-W and Type-S dictionaries for the (a) viscous Burgers~\eqref{eq:vb}, (b) KdV~\eqref{eq:kdv}, and (c) Allen--Cahn~\eqref{eq:allencahn} equations. The Type-S dictionary uses $p=7$ and $q=2$, while the Type-W dictionary uses $p=6$ and $q=6$, so that they have comparable sizes. In general, the Type-S feature matrices are more ill-conditioned than the Type-W ones, making identification more challenging.}\label{fig:svd}
\end{figure}

Third, 
the Type-S feature matrices are generally more ill-conditioned than the Type-W feature matrices. In Figure~\ref{fig:svd}, we plot the spectra, i.e., the singular values, of the feature matrices derived from the solution data of the (a) viscous Burgers~\eqref{eq:vb}, (b) KdV~\eqref{eq:kdv}, and (c) Allen--Cahn~\eqref{eq:allencahn} equations, with and without noise. We choose Type-S with $p=7$ and $q=2$ ($44$ terms) and Type-W with $p=6$ and $q=6$ ($42$ terms) so that their sizes are comparable. For all three PDEs, we observe that the spectrum of the Type-S feature matrices decays faster than that of the Type-W feature matrices. From the classical perspective, the rapid decay of the singular values reflects ill-conditioning, which makes sparse regression more challenging, particularly under noise.

These observations highlight the challenges of identification with Type-S dictionaries from the perspective of sparse linear regression. However, it is surprising that S-IDENT performs well for Type-S dictionaries, as demonstrated in Subsections~\ref{sec:DPspecific} and~\ref{sec:comparison_Stype}, and comparably to Type-W, as seen in Subsection~\ref{sec:comparison_typeW}. The identifiability theory established in~\cite{he2024much} provides a partial answer to this when the dictionaries contain only linear features; however, the general question of how much correlation PDE identification can tolerate remains open.

\section{Noise Models} \label{noise discussion}

Most existing works in the literature assume additive Gaussian random noise.
In~\cite{he2022robust, messenger2021weak}, the noise is determined by the root mean square (RMS) of the observations 
\begin{equation*}
\sigma_E  := \NSR_{E} \sqrt{\frac{1}{N}\sum_{n=1}^N |U_n|^2},
\end{equation*}
which is identical to the noise model~\eqref{eq_nsr_energy} we use in this work.
In~\cite{brunton2016discovering}, the authors show identification results with explicit $\sigma$ values. There are also three distinct strategies for scaling the noise relative to the given data. Suppose the perturbations are i.i.d. samples from a normal distribution $\mathcal{N}(0,\sigma^2)$ with standard deviation $\sigma$. In~\cite{rudy2017data}, the noise is added relative to the standard deviation of the clean data $\{U_n\}_{n=1}^N$\footnote{Typically, the standard deviation of discrete data is computed with a deficit of one degree of freedom; that is, $\frac{1}{N-1}\sum_{n=1}^N |U_n - \frac{1}{N}\sum_{m=1}^NU_m|^2$, so that it is unbiased. We report~\eqref{eq_nsr_variance} as used in~\cite{rudy2017data}.}:
\begin{equation}\label{eq_nsr_variance}
\sigma_S := \NSR_{S} \sqrt{\frac{1}{N}\sum_{n=1}^N |U_n - \frac{1}{N}\sum_{m=1}^NU_m|^2}.
\end{equation}
In~\cite{tang2023weakident}, the noise depends on a centralized RMS of the observations:
\begin{equation}\label{eq_nsr_central}
\sigma_C:= \NSR_{C}  \sqrt{\frac{1}{N}\sum_{n=1}^N |U_n - (\min_{m} U_m +\max_{m}U_m)/2|^2}.
\end{equation}
In~\eqref{eq_nsr_variance}--\eqref{eq_nsr_central}, $\NSR_S$, $\NSR_E$, and $\NSR_C$ are non-negative scaling factors that quantify the noise-to-signal ratios (NSR). The authors in~\cite{tang2023weakident} compared $\NSR_E$ with $\NSR_C$ through examples. 

We note that for any data, the inequality
\[\sigma_S\leq \min(\sigma_E,\sigma_C)\]
always holds. This suggests that with the same $\NSR_{S}=\NSR_{E}=\NSR_{C}=p\%$, the standard deviation of $p\%$ noise using model~\eqref{eq_nsr_variance} is smaller than those given by~\eqref{eq_nsr_energy} and~\eqref{eq_nsr_central}. Meanwhile, the relation between $\sigma_E$ and $\sigma_C$ is undetermined.
The actual standard deviation of the noise added by specifying the $\NSR$ in the above models depends on the underlying PDEs. For example, if the PDE is such that whenever $u$ is a solution, $u+C$ for any constant $C$ is also a solution (e.g., the advection-diffusion, heat, and linear Schr\"{o}dinger equations), then the standard deviation of $p\%$ noise via~\eqref{eq_nsr_energy} will vary as $C$ changes, whereas the standard deviations of noise via the other two models remain unchanged.  

The above discussion highlights the importance of clearly stating the NSR model and exercising caution when interpreting the noise level. This is particularly important when developing benchmarks for comparing the robustness of different identification methods.

\section{Type-S Partial Differential Equations in Subsection \ref{sec:DPspecific}}\label{Asec:StypePDE}

Figure~\ref{fig:solution_group1} shows the equations (A)--(F) used in Subsection~\ref{sec:DPspecific}. It shows clean trajectory data of (A) the Harry Dym equation~\eqref{eq:harry-dym}, (B) the thin film equation~\eqref{eq:thin-film}, (C) the viscous Hamilton--Jacobi equation~\eqref{eq:viscousHJ}, (D) 2D nonlinear advection (linear gradient)~\eqref{eq:2dnladv_v1}, (E) 2D nonlinear advection (squared gradient)~\eqref{eq:2dnladv_v2}, and (F) the DSW equation~\eqref{eq:dsw}. These are examples of PDEs representable by Type-S features but not by Type-W features alone. The identification results are presented in Table~\ref{tab:group1_all}.

\newcommand{\figw}{0.23\textwidth}
\newcommand{\figww}{0.46\textwidth}
\newcommand{\figwww}{0.92\textwidth}
\begin{figure}
\centering
\begin{tabular}{@{}cccc@{}}
(A) & (B) &
 \multicolumn{2}{c}{(F)}  \\
\includegraphics[width=\figw]{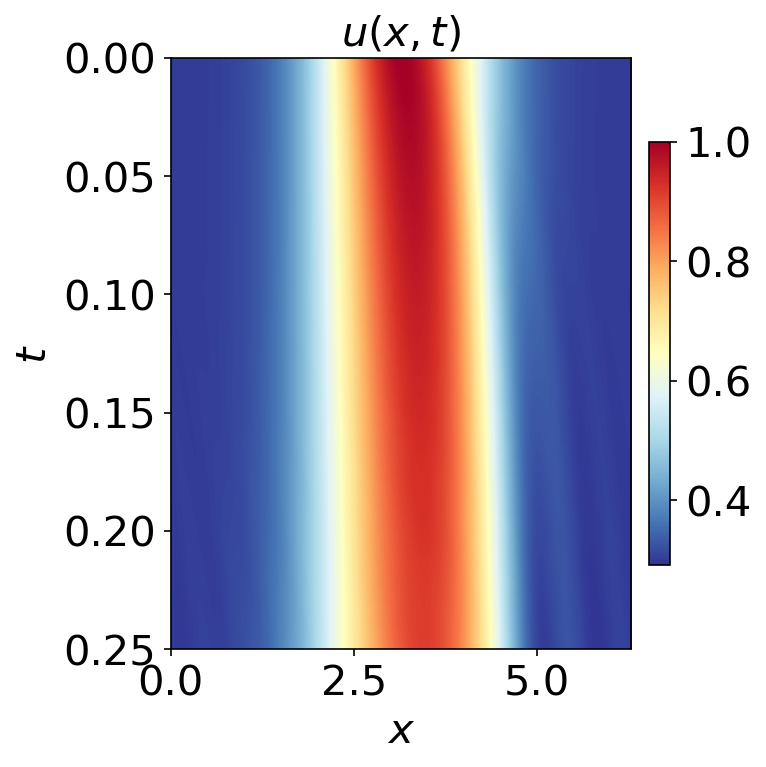} & 
\includegraphics[width=\figw]{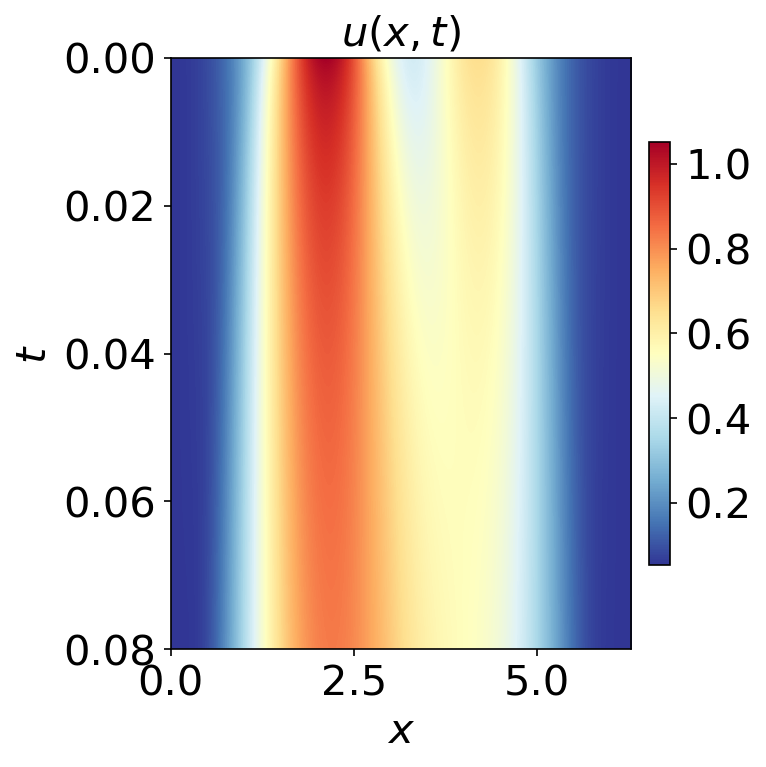} &
\multicolumn{2}{c}{\includegraphics[width=\figww]{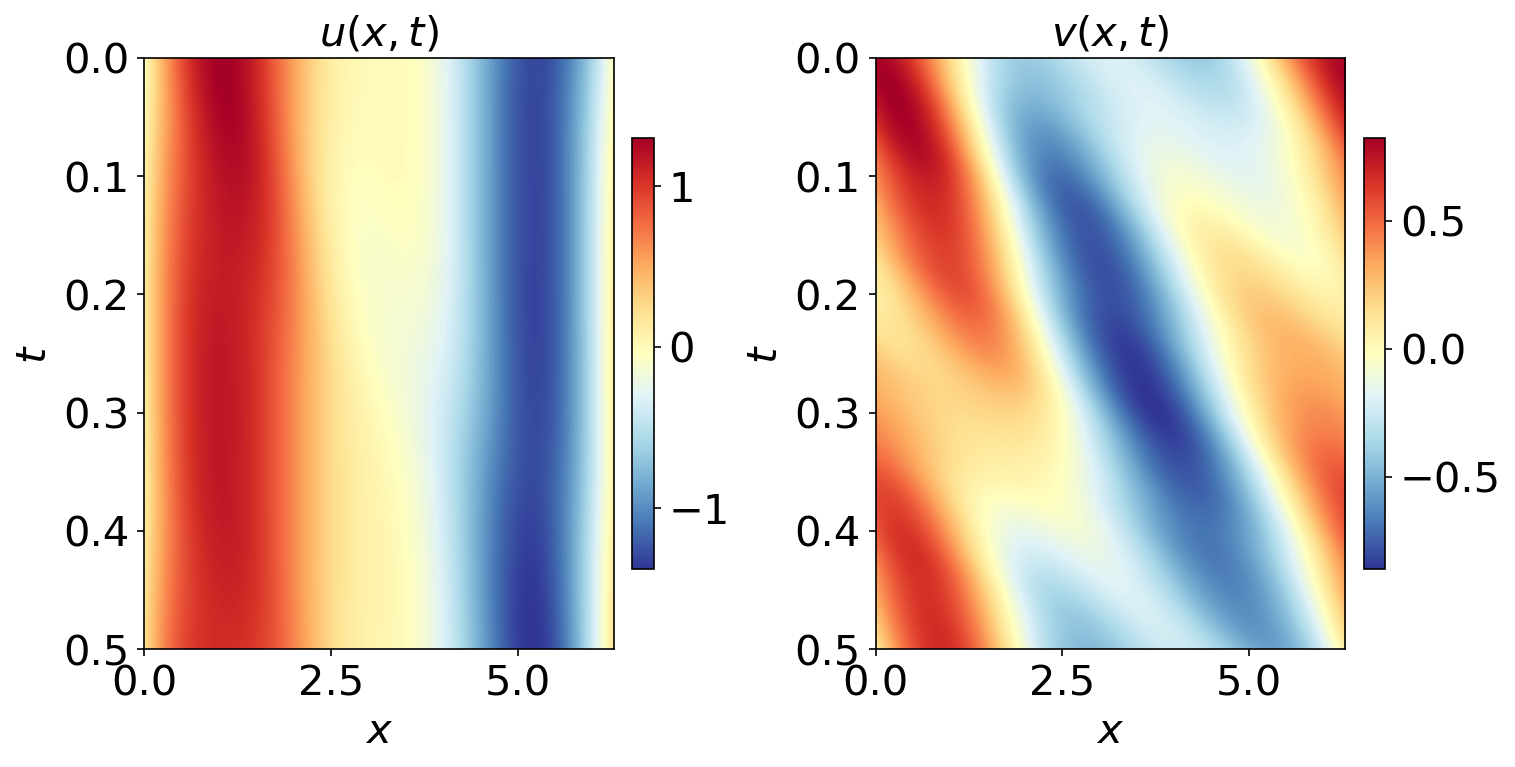}}
 \\
\multicolumn{4}{c}{(C)} \\
\multicolumn{4}{c}{\includegraphics[width=\figwww]{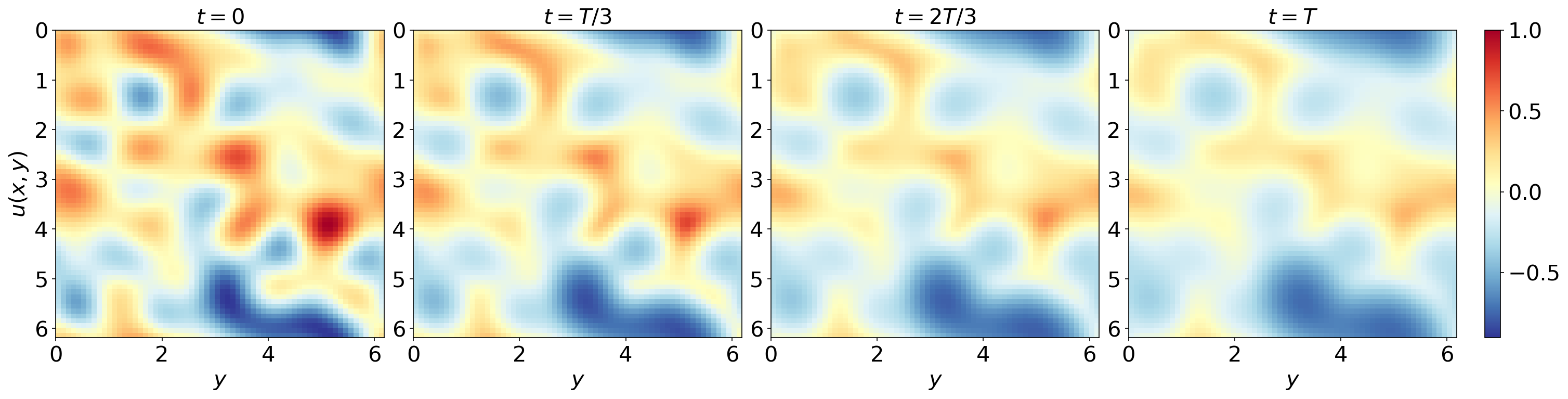}}\\
 \multicolumn{4}{c}{(D)} \\
\multicolumn{4}{c}{\includegraphics[width=\figwww]{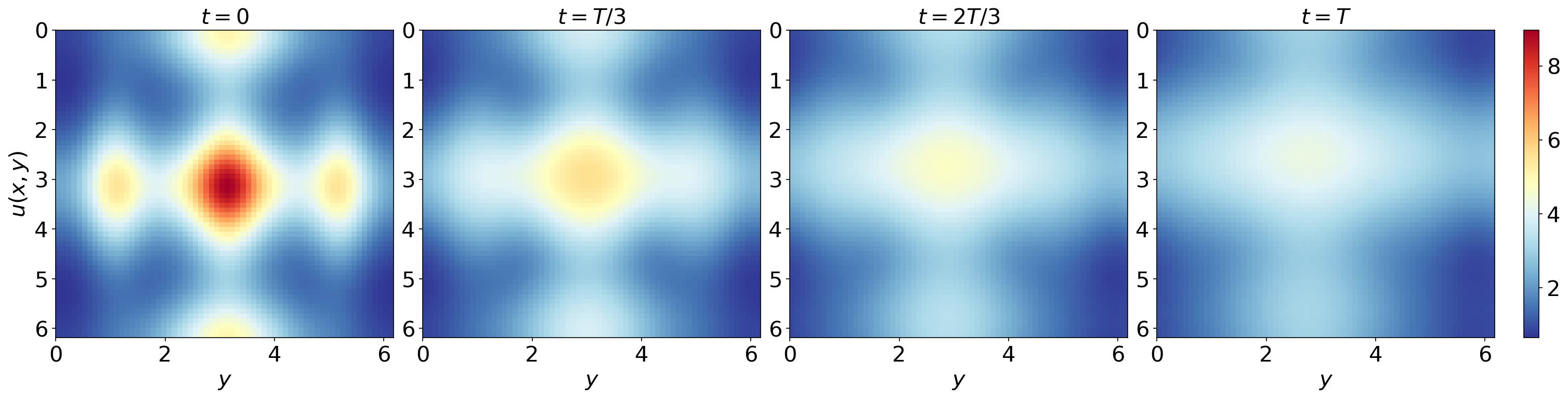}} \\
 \multicolumn{4}{c}{(E)} \\
\multicolumn{4}{c}{\includegraphics[width=\figwww]{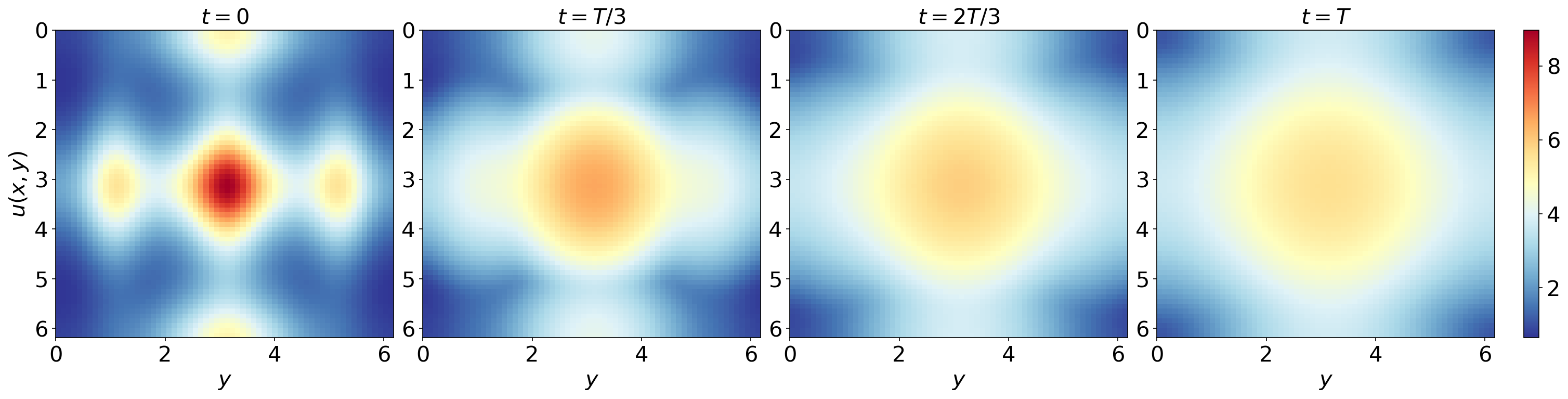}}
\end{tabular}
\caption{Clean trajectory data of (A) the Harry Dym equation~\eqref{eq:harry-dym}, (B) the thin film equation~\eqref{eq:thin-film}, (C) the viscous Hamilton--Jacobi equation~\eqref{eq:viscousHJ}, (D) 2D nonlinear advection (linear gradient)~\eqref{eq:2dnladv_v1}, (E) 2D nonlinear advection (squared gradient)~\eqref{eq:2dnladv_v2}, and (F) the DSW equation~\eqref{eq:dsw}. These are examples of PDEs representable by Type-S features but not by Type-W features alone. The identification results obtained by S-IDENT for the clean and noisy data are reported in Table~\ref{tab:group1_all}.}\label{fig:solution_group1}
\end{figure}

\section{Parameter Choices for the Compared Methods}
\label{sec:compare_details}

We present the parameter choices for the methods compared in Subsections~\ref{sec:comparison_Stype} and~\ref{sec:comparison_typeW}.
SINDy-PDE~\cite{rudy2017data} has a number of parameters to choose. We follow the default setup and employ SG differentiation with a fixed window size of $N/10$, where $N$ is the number of grid points along each dimension. 
There are also an $L^2$ regularization parameter $\lambda$ for conditioning and a tolerance step $\tau$ for hard-threshold sweeping. While fixing $\lambda=1\times10^{-5}$, we find that in order to correctly identify the tested PDEs~\eqref{eq:vb}--\eqref{eq:allencahn} from noiseless data, we need to adjust $\tau$ case by case. For the viscous Burgers equation~\eqref{eq:vb}, we use  $\tau=1.0$; for the KdV equation~\eqref{eq:kdv}, we use   $\tau=4.0$; and for the Allen--Cahn equation~\eqref{eq:allencahn}, we use  $\tau=0.5$. 

For Robust-IDENT, there are parameters for SDD and the model validation: in SDD, for the MLS smooting, we set polynomial degree to be $2$, window length to be $25$, standard deviation of the Gaussian weight to be $7$; for the differentiation, we use  FD with $5$-points; and for the cross-validation model selection, we set the associated parameter to be $5\times 10^{-3}$. With larger dictionaries and higher levels of noise as we address in this paper, Robust-IDENT finds many excessive terms even with clean data; we thus omit identified features with coefficients whose absolute values are smaller than $5\times 10^{-4}$.   

For Weak-SINDy, Weak-IDENT, and S-IDENT, we employ the default parameters, and the denoising parameters for all these methods adapt to the noise level.
For Weak-IDENT, we follow the default setup, downsampling the trajectory data by $5$ in both time and space and setting the trimming parameter to $0.05$.

\section{More Results with Type-W Dictionaries}\label{sec:w-type-more}
\begin{figure}
\centering
\setlength{\tabcolsep}{2pt}
\renewcommand{\arraystretch}{1.2}
\begin{tabular}{cccc}
& \textbf{Viscous Burgers} & \textbf{KdV} & \textbf{Allen-Cahn} \\
\adjustbox{valign=c}{\rotatebox{90}{\textbf{TPR}}} &
\adjustbox{valign=c}{\includegraphics[width=0.28\textwidth]{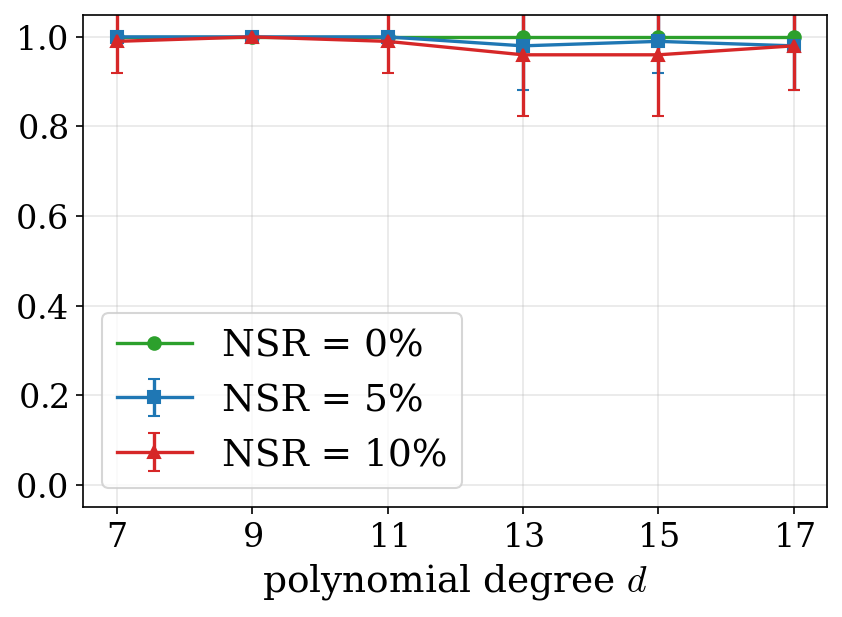}} &
\adjustbox{valign=c}{\includegraphics[width=0.28\textwidth]{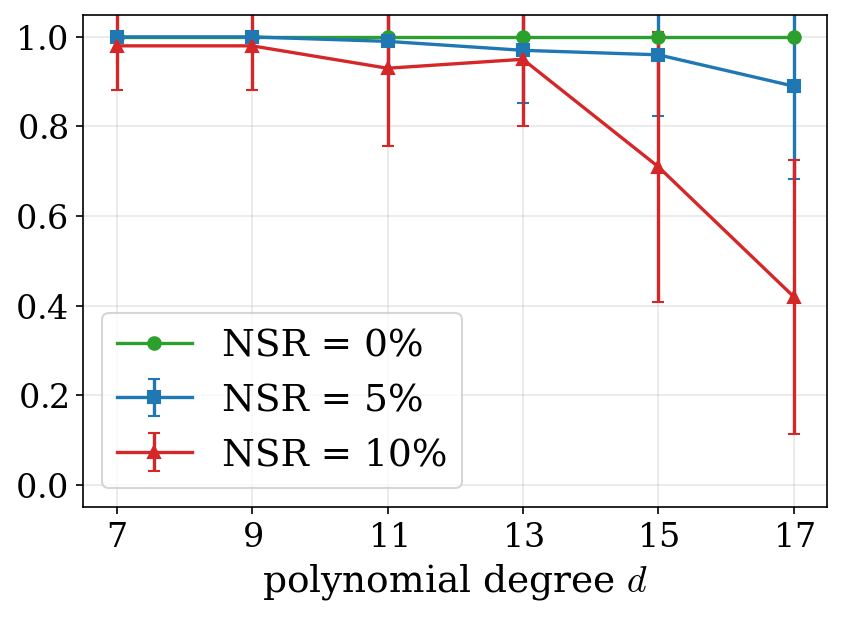}} &
\adjustbox{valign=c}{\includegraphics[width=0.28\textwidth]{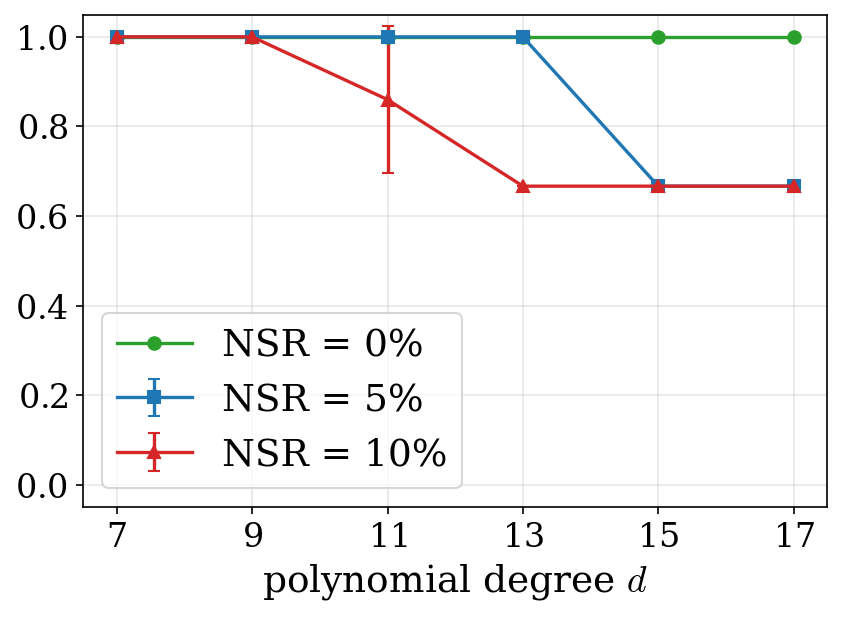}} \\
\adjustbox{valign=c}{\rotatebox{90}{\textbf{PPV}}} &
\adjustbox{valign=c}{\includegraphics[width=0.28\textwidth]{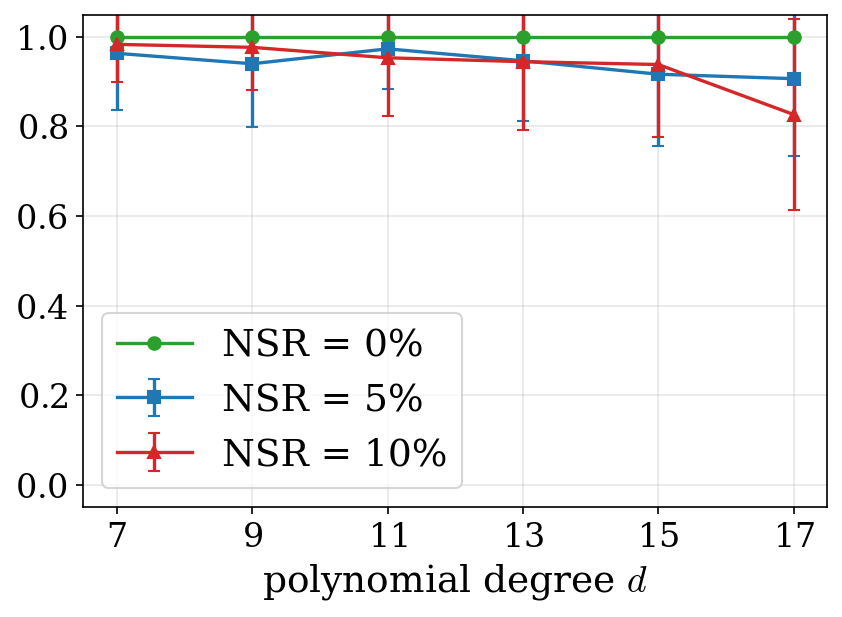}} &
\adjustbox{valign=c}{\includegraphics[width=0.28\textwidth]{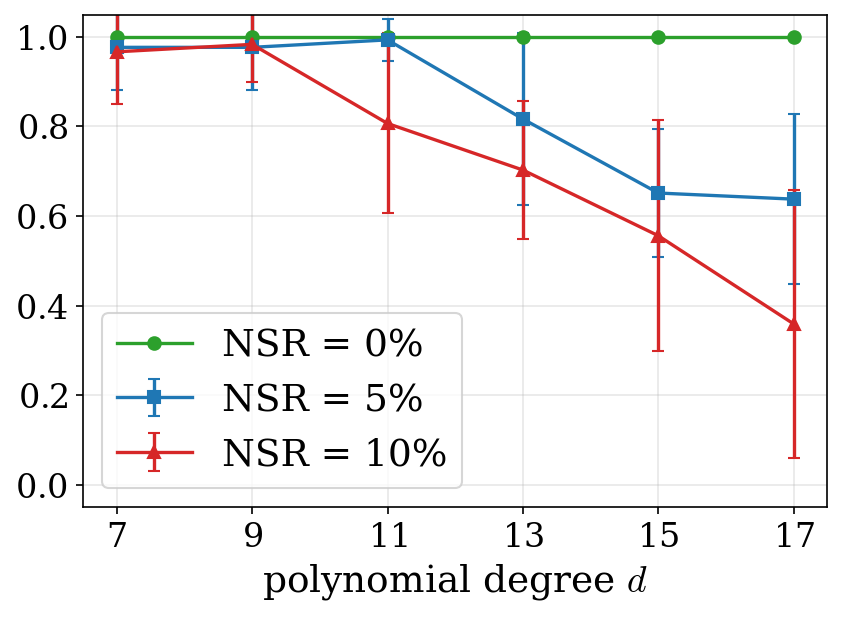}} &
\adjustbox{valign=c}{\includegraphics[width=0.28\textwidth]{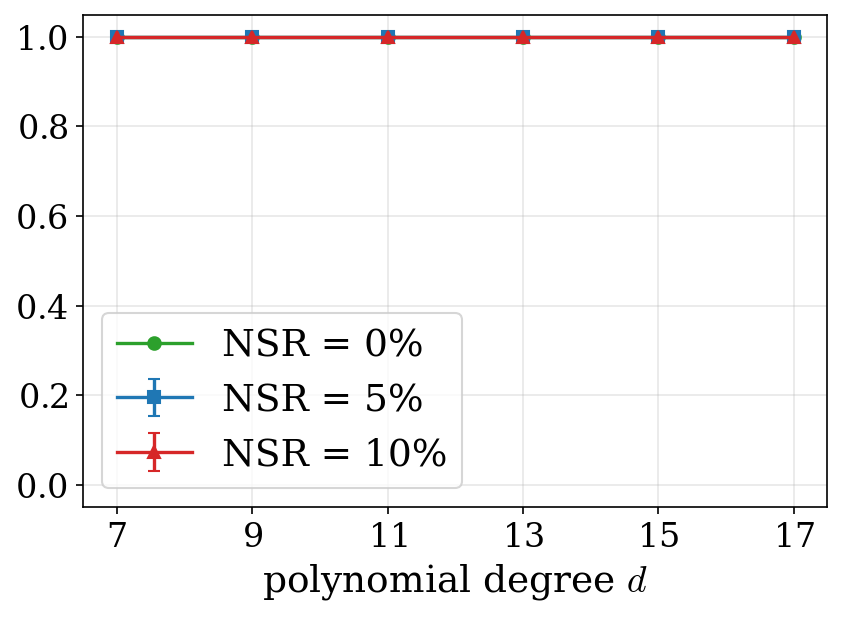}} \\
\adjustbox{valign=c}{\rotatebox{90}{$\boldsymbol{E_{\mathrm{in}}}$}} &
\adjustbox{valign=c}{\includegraphics[width=0.28\textwidth]{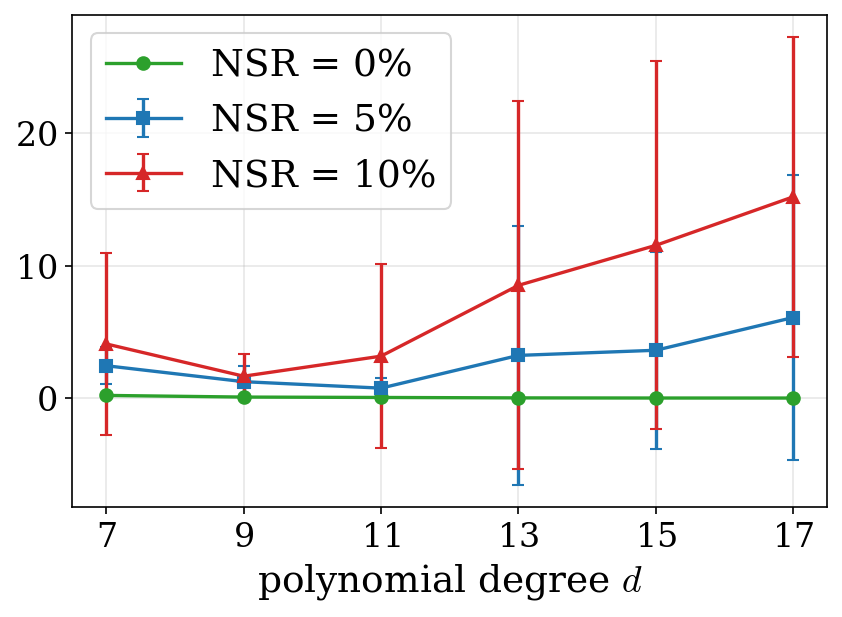}} &
\adjustbox{valign=c}{\includegraphics[width=0.28\textwidth]{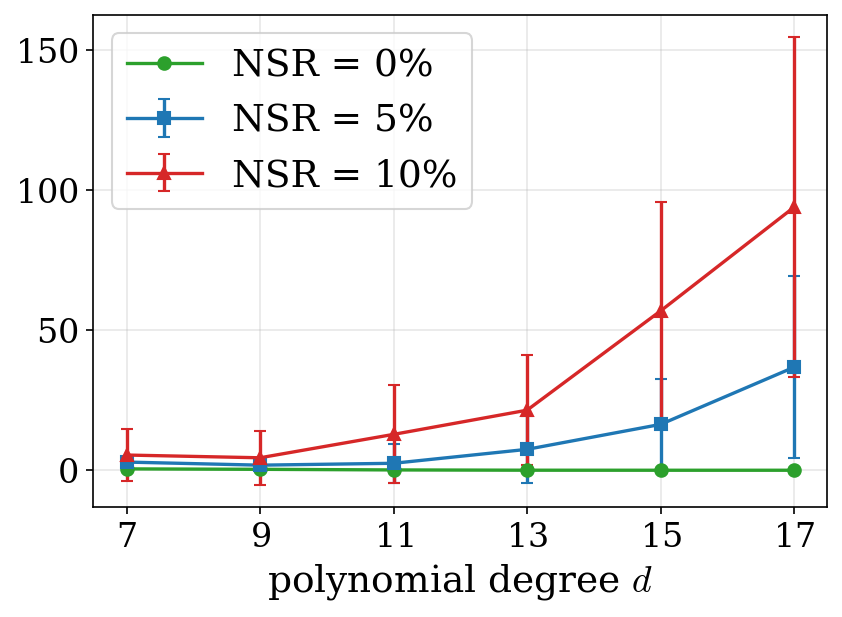}} &
\adjustbox{valign=c}{\includegraphics[width=0.28\textwidth]{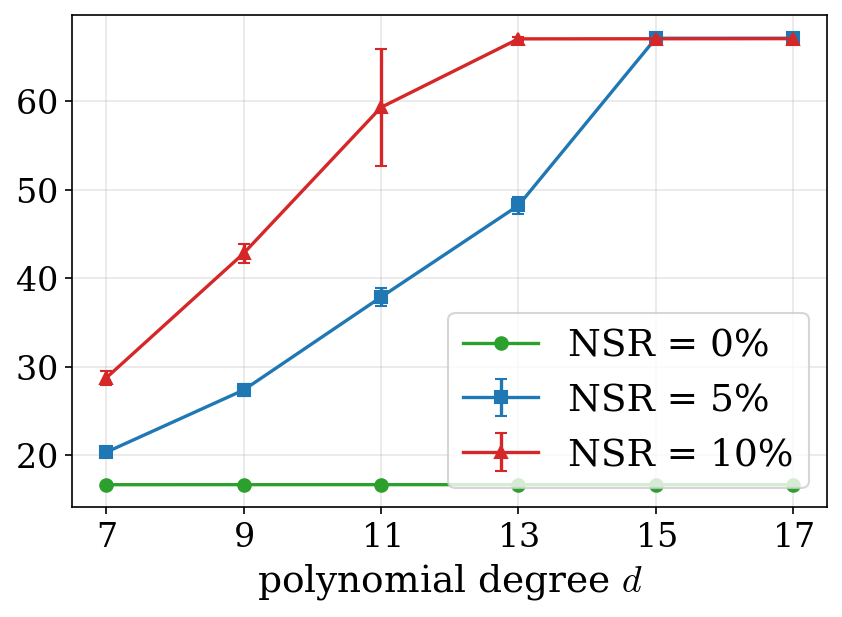}} \\
\adjustbox{valign=c}{\rotatebox{90}{$\boldsymbol{E_{\mathrm{out}}}$}} &
\adjustbox{valign=c}{\includegraphics[width=0.28\textwidth]{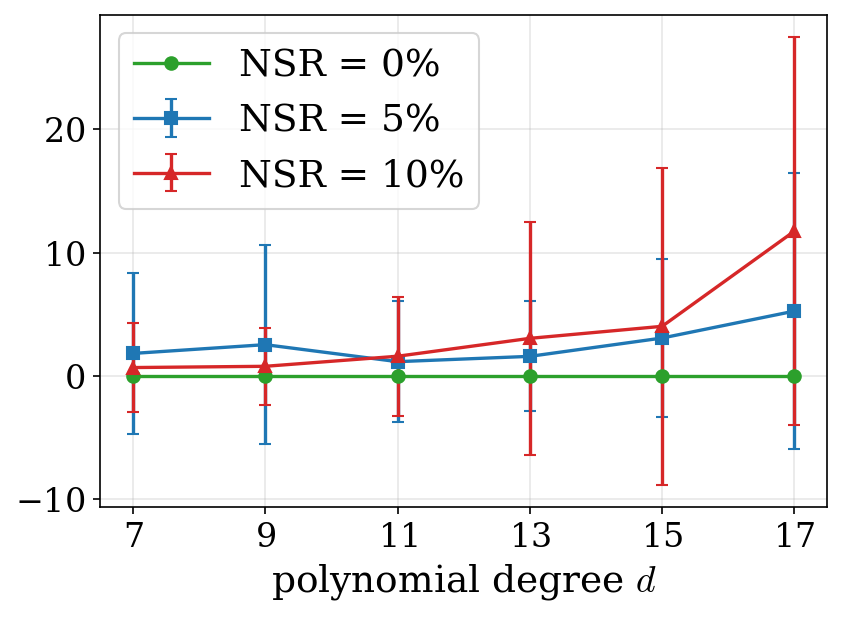}} &
\adjustbox{valign=c}{\includegraphics[width=0.28\textwidth]{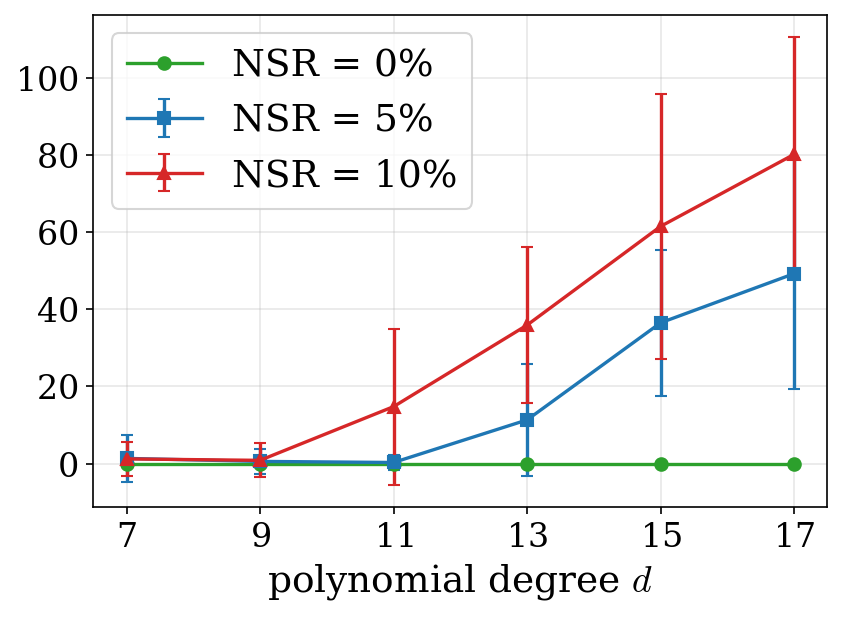}} &
\adjustbox{valign=c}{\includegraphics[width=0.28\textwidth]{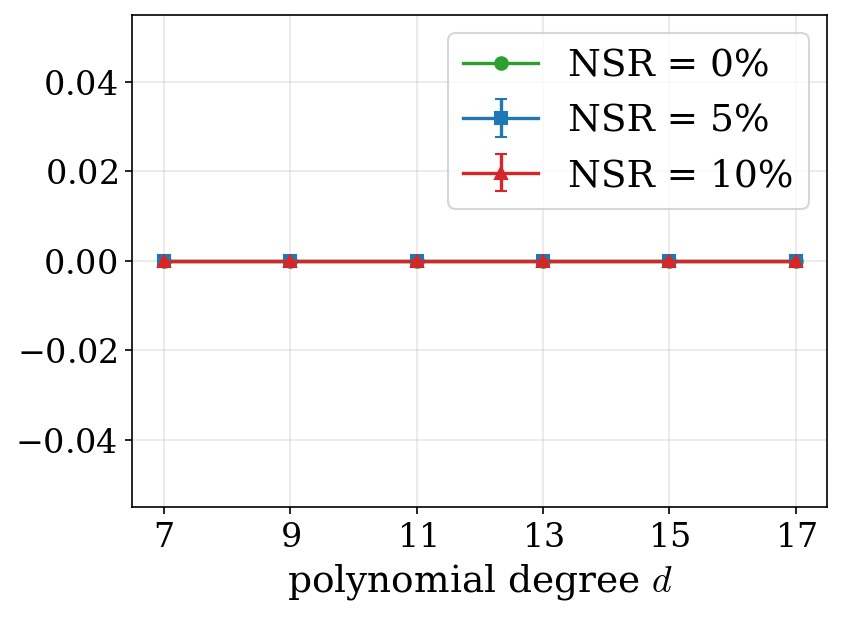}} \\
\adjustbox{valign=c}{\rotatebox{90}{\textbf{Recovery}}} &
\adjustbox{valign=c}{\includegraphics[width=0.28\textwidth]{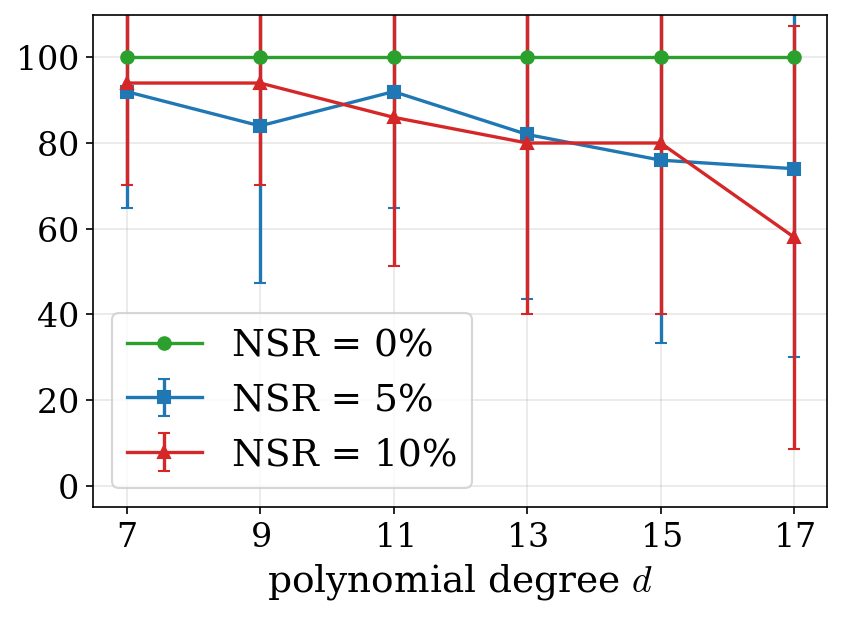}} &
\adjustbox{valign=c}{\includegraphics[width=0.28\textwidth]{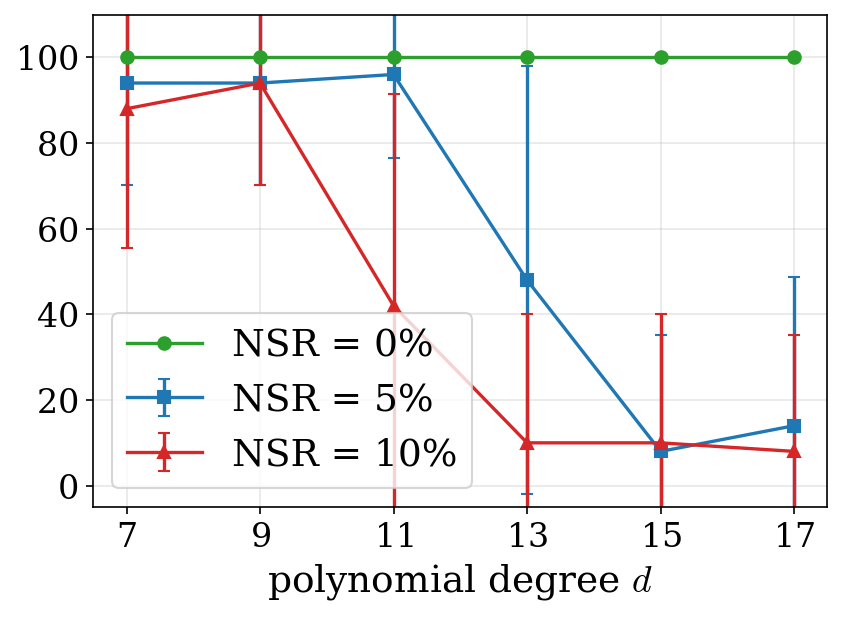}} &
\adjustbox{valign=c}{\includegraphics[width=0.28\textwidth]{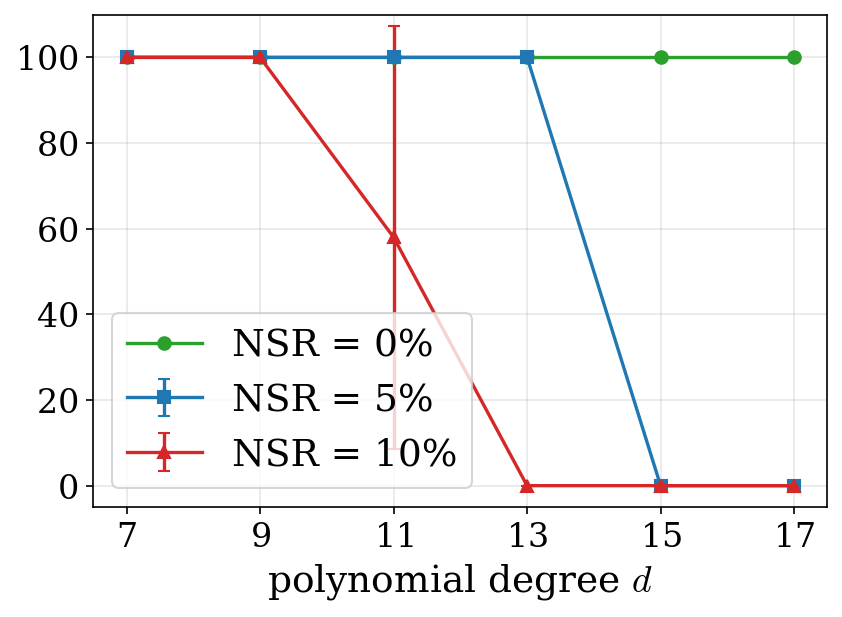}} \\
\end{tabular}
\caption{Type-W dictionary, SG polynomial order vs.\ S-IDENT results: Identification performance of the Type-W dictionary for the viscous Burgers equation~\eqref{eq:vb}, the KdV equation~\eqref{eq:kdv}, and the Allen--Cahn equation~\eqref{eq:allencahn}. We fix the SG polynomial degree to $7$ as the default in this paper. The corresponding Type-S results are presented in Figure~\ref{fig:poly-order-DP}.}
\label{fig:poly-order-PD}
\end{figure}

We present additional results for Subsections~\ref{sec:order_accuracy} and~\ref{sec:ablation_feature_approx}.   
Figure~\ref{fig:poly-order-PD} shows the identification performance of S-IDENT on noisy data collected from viscous Burgers~\eqref{eq:vb}, KdV~\eqref{eq:kdv}, and Allen--Cahn~\eqref{eq:allencahn} using SURE-SG with different orders of accuracy. As with the results shown in Figure~\ref{fig:poly-order-DP} for Type-S dictionaries, we observe that $d=7$ gives the best results in most cases, thus justifying our default parameter setting. A higher order of accuracy often leads to worse identification, especially when the noise level is high.
Table~\ref{tab:diff_compare_pd} compares the performance of the three differentiation strategies (Direct, Repeated, and Adaptive) described in Subsection~\ref{sec:ablation_feature_approx}, applied to Type-W $(6,4)$ dictionaries. We also observe that the Direct approach proposed for S-IDENT yields the best results in most cases. One slight difference from the Type-S results in Table~\ref{tab:diff_compare_dp} is that the overall accuracy is generally higher with Type-W, and this phenomenon is numerically investigated in Appendix~\ref{sec:S-correlation} from the perspective of pairwise correlation among features.

\begin{table}
\centering
\small
\begin{tabular}{llccccc}
\toprule
NSR & Method & TPR & PPV & $E_{\text{in}}$ & $E_{\text{out}}$ & E.R. (\%) \\
\midrule
\multicolumn{7}{c}{(a) \textit{Viscous Burgers} \eqref{eq:vb}}\\
\midrule
& Direct & \textbf{1.000 $\pm$ 0.000} & 0.987 $\pm$ 0.065 & 0.831 $\pm$ 0.247 & 0.444 $\pm$ 2.209 & 96.0 $\pm$ 19.6 \\
1\%  & Repeated & \textbf{1.000 $\pm$ 0.000} & \textbf{0.993 $\pm$ 0.047} & \textbf{0.216 $\pm$ 0.150} & \textbf{0.235 $\pm$ 1.648} & \textbf{98.0 $\pm$ 14.0} \\
 & Adaptive & \textbf{1.000 $\pm$ 0.000} & 0.960 $\pm$ 0.123 & 2.724 $\pm$ 1.195 & 3.485 $\pm$ 13.250 & 90.0 $\pm$ 30.0 \\
\cmidrule(lr){2-7}
& Direct & 0.990 $\pm$ 0.070 & \textbf{0.970 $\pm$ 0.104} & 3.748 $\pm$ 6.842 & \textbf{1.045 $\pm$ 3.709} & \textbf{90.0 $\pm$ 30.0} \\
10\%  & Repeated & 0.990 $\pm$ 0.070 & 0.953 $\pm$ 0.129 & 3.950 $\pm$ 6.996 & 1.684 $\pm$ 5.429 & 86.0 $\pm$ 34.7 \\
 & Adaptive & \textbf{1.000 $\pm$ 0.000} & 0.933 $\pm$ 0.145 & \textbf{3.662 $\pm$ 1.597} & 2.804 $\pm$ 6.582 & 82.0 $\pm$ 38.4 \\
\midrule
\multicolumn{7}{c}{ (b) \textit{KdV} \eqref{eq:kdv} } \\
\midrule
 & Direct & \textbf{1.000 $\pm$ 0.000} & \textbf{0.987 $\pm$ 0.065} & 1.539 $\pm$ 0.902 & \textbf{0.550 $\pm$ 2.848} & \textbf{96.0 $\pm$ 19.6} \\
1\% & Repeated & \textbf{1.000 $\pm$ 0.000} & \textbf{0.987 $\pm$ 0.065} & \textbf{0.927 $\pm$ 1.354} & 0.836 $\pm$ 4.653 & \textbf{96.0 $\pm$ 19.6} \\
 & Adaptive & \textbf{1.000 $\pm$ 0.000} & 0.717 $\pm$ 0.230 & 20.682 $\pm$ 5.319 & 26.561 $\pm$ 23.744 & 36.0 $\pm$ 48.0 \\
\cmidrule(lr){2-7}
 & Direct & \textbf{0.980 $\pm$ 0.098} & \textbf{0.967 $\pm$ 0.115} & \textbf{5.428 $\pm$ 9.439} & \textbf{1.219 $\pm$ 4.335} & \textbf{88.0 $\pm$ 32.5} \\
10\% & Repeated & 0.970 $\pm$ 0.119 & 0.953 $\pm$ 0.129 & 5.836 $\pm$ 11.685 & 1.755 $\pm$ 6.436 & 84.0 $\pm$ 36.7 \\
 & Adaptive & 0.970 $\pm$ 0.119 & 0.625 $\pm$ 0.170 & 35.290 $\pm$ 20.371 & 38.237 $\pm$ 20.007 & 8.0 $\pm$ 27.1 \\
\midrule
\multicolumn{7}{c}{(c) \textit{Allen-Cahn} \eqref{eq:allencahn}} \\
\midrule
 & Direct & \textbf{1.000 $\pm$ 0.000} & \textbf{1.000 $\pm$ 0.000} & 16.798 $\pm$ 0.044 & \textbf{0.000 $\pm$ 0.000} & \textbf{100.0 $\pm$ 0.0} \\
1\% & Repeated & \textbf{1.000 $\pm$ 0.000} & \textbf{1.000 $\pm$ 0.000} & 16.865 $\pm$ 0.071 & \textbf{0.000 $\pm$ 0.000} & \textbf{100.0 $\pm$ 0.0} \\
 & Adaptive & \textbf{1.000 $\pm$ 0.000} & \textbf{1.000 $\pm$ 0.000} & \textbf{16.276 $\pm$ 0.094} & \textbf{0.000 $\pm$ 0.000} & \textbf{100.0 $\pm$ 0.0} \\
\cmidrule(lr){2-7}
 & Direct & \textbf{1.000 $\pm$ 0.000} & \textbf{1.000 $\pm$ 0.000} & 28.731 $\pm$ 0.788 & \textbf{0.000 $\pm$ 0.000} & \textbf{100.0 $\pm$ 0.0} \\
10\% & Repeated & \textbf{1.000 $\pm$ 0.000} & \textbf{1.000 $\pm$ 0.000} & 30.577 $\pm$ 2.366 & \textbf{0.000 $\pm$ 0.000} & \textbf{100.0 $\pm$ 0.0} \\
 & Adaptive & \textbf{1.000 $\pm$ 0.000} & \textbf{1.000 $\pm$ 0.000} & \textbf{17.428 $\pm$ 2.098} & \textbf{0.000 $\pm$ 0.000} & \textbf{100.0 $\pm$ 0.0} \\
\bottomrule
\end{tabular}
\caption{Type-W dictionary, variant study of Direct (Proposed), Repeated, and Adaptive for S-IDENT in Subsection~\ref{sec:ablation_feature_approx}:
For each (PDE, NSR) and each metric, the best of the three SURE-SG strategies is shown in \textbf{bold}. Mean $\pm$ std over $50$ independent trials. The Direct strategy shows the best overall performance. The corresponding Type-S results are presented in Table~\ref{tab:diff_compare_dp}.}
\label{tab:diff_compare_pd}
\end{table}

\end{document}